\documentclass{amsart}

\usepackage{hyperref}

\usepackage{enumerate}
\usepackage{enumitem}
\usepackage{moreenum}
\usepackage{amsmath} 
\usepackage{amssymb}
\usepackage{mathrsfs}

\usepackage{enumitem}
\usepackage{wasysym}

\newtheorem{theorem}{Theorem}[section] 
\newtheorem{claim}[theorem]{Claim}
\newtheorem{sclaim}[theorem]{Subclaim}
\newtheorem{lemma}[theorem]{Lemma} 
 
\newtheorem{observation}[theorem]{Observation} 
 
\newtheorem{fact}[theorem]{Fact} 

\newcommand{\thistheoremname}{}
\newtheorem*{genericthm*}{\thistheoremname}
\newenvironment{namedthm*}[1]
{\renewcommand{\thistheoremname}{#1}%
	\begin{genericthm*}}
	{\end{genericthm*}}

\theoremstyle{definition}
\newtheorem{definition}[theorem]{Definition}

\newtheorem{problem}[theorem]{Problem}

\newtheorem{convention}[theorem]{Convention}

\theoremstyle{remark}

\newcommand{\lh}{{\ell g}}
\newcommand{\supp}{{\rm supp}}

\newcommand{\rest}{{\restriction}}
\newcommand{\dom}{{\rm dom}} 

\newcommand{\cov}{{\rm cov}}
\newcommand{\ran}{{\rm ran}}

\newcommand{\stem}{{\rm stem}}

\newcommand{\then}{{\underline{then}}}

\newcommand{\Then}{{\underline{Then}}}

\newcommand{\varp}{{\varepsilon}}

\newcommand{\bbS}{{\mathbb S}}

\newcommand{\sP}{{\mathscr P}}

\newcommand{\bbG}{{\mathbb G}}

\newcommand{\cM}{{\mathcal M}}

\newcommand{\bbP}{{\mathbb P}}

\newcommand{\bbQ}{{\mathbb Q}}

\newcommand{\bbN}{{\mathbb N}}

\newcommand{\tieconcat}{%
	\mathbin{\mathpalette\dotieconcat\relax}%
}
\newcommand{\dotieconcat}[2]{
	\text{\raisebox{.8ex}{$\smallfrown$}}%
}

\newcount\skewfactor
\def\mathunderaccent#1#2 {\let\theaccent#1\skewfactor#2
	\mathpalette\putaccentunder}
\def\putaccentunder#1#2{\oalign{$#1#2$\crcr\hidewidth
		\vbox to.2ex{\hbox{$#1\skew\skewfactor\theaccent{}$}\vss}\hidewidth}}
\def\name{\mathunderaccent\tilde-3 }

\newenvironment{PROOF}[2][\proofname.]
{\begin{proof}[#1]\renewcommand{\qedsymbol}{\textsquare$_{\rm #2}$}}
	{\end{proof}}

\usepackage{hyperref}

\begin{document}
	
	\title {On the weak Borel chromatic number and cardinal invariants of the continuum}
	\author {Márk Poór}
	\address{Einstein Institute of Mathematics\\
		Edmond J. Safra Campus, Givat Ram\\
		The Hebrew University of Jerusalem\\
		Jerusalem, 91904, Israel\\}
	\email{sokmark@gmail.com}

	\author {Saharon Shelah}
	\address{Einstein Institute of Mathematics\\
		Edmond J. Safra Campus, Givat Ram\\
		The Hebrew University of Jerusalem\\
		Jerusalem, 91904, Israel\\
		and \\
		Department of Mathematics\\
		Hill Center - Busch Campus \\ 
		Rutgers, The State University of New Jersey \\
		110 Frelinghuysen Road \\
		Piscataway, NJ 08854-8019 USA}
	\email{shelah@math.huji.ac.il}
	\urladdr{http://shelah.logic.at}
	\thanks{$^\dag$The first author was supported by the Excellence Fellowship Program for International Postdoctoral Researchers of The Israel Academy of Sciences and Humanities, and by the National Research, Development and Innovation Office
		– NKFIH, grants no. 124749, 129211. \\
		$^\ast$The second author was supported	by the Israel Science Foundation grant
		1838/19. \\
		 Research of both authors partially supported by NSF grant no: DMS 1833363.  \\ Paper 1226 on Shelah's list.
		References like \cite[Th0.2=Ly5]{Sh:950} mean the label
		of Th.0.2 is y5.  The reader should note that the version on the second author's website is usually more updated than the one in the mathematical archive.}
	
	
	
	
	
	\begin{abstract}
		We prove that consistently, cov($\mathcal{M})< \lambda_\mathbf{0} < \lambda_\mathbf{1} < \lambda_\mathbf{\infty} < 2^{\aleph_0}$, where $\lambda_\mathbf{0}$ denotes the weak Borel chromatic number of the Kechris-Solecki-Todor{\v{c}}evi{\'c} graph $\mathbb{G}_0$, that is, the minimal cardinality of a $\mathbb{G}_0$-independent Borel covering of $2^\omega$, while $\lambda_\mathbf{1}$ and $\lambda_\infty$ are the corresponding invariants of the graph $\mathbb{G}_1$ and the simple graph associated with the equivalence relation $E_0$.
	\end{abstract}
	
	\maketitle
	\numberwithin{equation}{section}
	\setcounter{section}{-1}

	\section{Introduction}

	Borel graphs and their combinatorial properties have become a growing area of research in the last two decades and it has interesting connections with other areas such as the theory of graph limits, countable group actions, paradoxical
	decompositions, as well as ergodic theory.
	
	The Borel chromatic number was studied and defined in \cite{KST}\cite{lecmil}.
	 For a graph $G = (X,E)$ on a Polish space $X$ its Borel chromatic number $\chi_\text{B}(G)$ is the least cardinal $\kappa$, such that for some Polish space $Y$ there exists a Borel coloring $c:X \to Y$ of $G$ with $|\ran(c)|= \kappa$ (i.e.\ for each $y \in Y$ the preimage $c^{-1}(y)$ is $G$-independent). It is clear by the theory of Polish spaces and the \emph{Perfect Set Property} of analytic sets that this number is an element of the set $\{0,1,2, \dots\} \cup \{\aleph_0, 2^{\aleph_0} \}$.
	 
	 The theory was extended by S. Geschke, who showed that if $X$ is Polish, then for each closed graph $G = (X,E)$ without perfect cliques, as well as for each locally countable $F_\sigma$ graph $G = (X,E)$ (i.e.\ each node has degree at most $\aleph_0$) there exists some ccc.\ forcing making the continuum large, while $X$ can be covered by $\aleph_1$-many Borel (in fact, closed) $G$-independent sets \cite{geschke}. Later M. Gaspar- S. Geschke \cite{weabo} have defined the weak Borel chromatic number of a fixed graph $G$ as the least possible cardinal $\kappa$, such that the underlying space can be covered by $\kappa$-many $G$-independent Borel sets.  Note that if either chromatic number is at most countable, then they coincide.
	 (Here we remark that S. Geschke had defined the weak Borel chromatic number (of the graph $G= (X,E)$) as the smallest cardinal $\kappa$, such that there exists a coloring $c: X \to \kappa$ with Borel fibers (i.e.\ for each $\alpha < \kappa$ the preimage $c^{-1}(\alpha)$ is Borel) \cite{geschke}. Note that if for a fixed graph either variant of the weak Borel chromatic number is at most $\aleph_1$, then they coincide.)

	 In the celebrated paper of A. Kechris-S. Solecki- S. Todor{\v{c}}evi{\'c} \cite{KST} the graph $\mathbb{G}_0$ (with $X = 2^\omega$) was constructed, and proved to be minimal among analytic graphs of uncountable Borel chromatic number in the sense that for each simple graph $G = (Y,F)$, where $Y$ is Polish, $E \in \mathbf{\Sigma}^1_1(Y^2)$ either $\chi_\text{B}(G) \leq \aleph_0$, or there exists a continuous homomorphism $f:2^\omega \to Y$ from $\mathbb{G}_0$ to $G$ (i.e.\ whenever $(x,x') \in E_{\mathbb{G}_0}$, then $(f(x), f(x')) \in F$ holds necessarily).  This also implies that whenever $G = (Y,F)$ is an analytic graph on a Polish space with uncountable weak Borel chromatic number, it is at least as the weak Borel chromatic number of $\mathbb{G}_0$.
	 
	 While the graph $\mathbb{G}_0$ is acyclic, so it can be colored by two colors,  B. Miller showed that the measurable chromatic number of it is $3$ \cite{miller2008measurable}. He asked whether anything can be said about the weak Borel chromatic number of $\mathbb{G}_0$ compared to other cardinal characteristics of the continuum.
	 In \cite{KST} not only has the authors verified that $\chi_\text{B}(G) > \aleph_0$, but it also followed from their argument, that each $\mathbb{G}_0$-independent Baire-measurable set $S \subseteq 2^\omega$ must be meager. This immediately implies that cov($\mathcal{M}$) is a lower bound for 
	 the weak Borel chromatic number of $\mathbb{G}_0$ as well.

	Due to 	M. Gaspar and S. Geschke, independently of this work, various Borel chromatic numbers of graphs were computed in models of set theory obtained by forcing with countable support iteration of uniform tree-forcing notions \cite{weabo}, or see further results by R. Banerjee, M. Gaspar \cite{banerjee2022borel}.	
	 Earlier F. Adams and J. Zapletal had studied cardinal invariants of closed graphs \cite{adams2018cardinal}.

	
	
	\section{Preliminaries, notations}
	Under ordinals we always mean von Neumann ordinals, and for a set $X$ the symbol $|X|$ always refers to the smallest ordinal with the same cardinality.  For a set $X$ the symbol $\sP(X)$ denotes the power set of $X$, while if $\kappa$ is an cardinal we use the standard notation $[X]^\kappa$ for $\{Y \in \sP(X): \ |Y| = \kappa\}$, similarly for $[X]^{<\kappa}$, $[X]^{\leq \kappa}$, etc.\ By a sequence we mean a function on an ordinal, where for a sequence $ \overline{s}= \langle s_\alpha: \ \alpha < \dom(\overline{s}) \rangle$ the length of $\overline{s}$ (in symbols $\lh(\overline{s})$) denotes $\dom(\overline{s})$. We denote the empty sequence by $\langle\rangle$. Moreover, for sequences $\overline{s}$, $\overline{t}$, we  let $\overline{s} \tieconcat \overline{t}$ denote the natural concatenation of them (of length $\lh(\overline{s}) + \lh(\overline{t})$).
	For a set $X$, and ordinal $\alpha$ we use	$^\alpha X = \{ \overline{s}: \ \lh(\overline{s}) =  \alpha, \ \ran(\overline{s}) \subseteq X\}$, and for cardinals $\lambda$, $\kappa$ we use the symbol  $\lambda^\kappa = |^\kappa \lambda|$ (that is, the least ordinal equivalent to it).
	
	For a finite sequence $\bar s \in \ ^{\omega>}2$ the symbol $[\bar s]$ stands for the basic open set in $2^\omega$ that $\bar s$ determines, i.e.
	$$[\bar s] = \{ x \in \ ^\omega 2: \ x \supseteq \bar s \}.$$
	A tree $T$ is a downward closed set consisting of finite sequences.
	
	If $\varphi(n,x)$ is a formula, then $\forall^\infty n \varphi(n,x)$ is true, if for all but finitely many $n \in \omega$ $\varphi(n,x)$ is true, $\exists^\infty n$ stands for there exists infinitely many $n$, and we use the quantifier $\exists!y$ as ``there exists a unique $y$".
	For $r,r' \in \omega^\omega$ under $r \leq^* r'$ we mean that $\forall^\infty n$ $r_n \leq r'_n$.
	
	Concerning forcing, $q \leq p$ means that $q$ is stronger, and for a notion of forcing $\bbP$ the term $1_\bbP$ stands for the unique largest element of $\bbP$. 
	
	\section{The forcing construction}
	
	Let $\bar{\mathfrak{s}} = \langle \mathfrak{s}_n: \ n \in \omega \rangle$ be fixed, such that
	\begin{enumerate}[label = $(\text{x}_{\arabic*})$, ref = $(\text{x}_{\arabic*})$]
		\item \label{x1} for each $n \in \omega$: $\mathfrak{s}_n \in \ ^n \omega$,
		\item \label{x2} the set $\bigcup_{n \in \omega} [\mathfrak{s}_n]$ is dense in $2^\omega$.
	\end{enumerate}

	Recall the definition of the graph $\mathbb{G}_0(\bar s)$ on the Cantor space \cite{KST}:
	\begin{definition}
		The graph $\mathbb{G}_{0}(\bar{\mathfrak{s}})$ is defined as follows: 
		$$\mathbb{G}_{0}(\bar{\mathfrak{s}}) = \{ (\mathfrak{s}_n\tieconcat \langle 0 \rangle \tieconcat x, \mathfrak{s}_n\tieconcat \langle 1 \rangle \tieconcat x): \ n \in \omega, \ x \in 2^\omega \} \subseteq [2^\omega]^2.$$
	\end{definition}
	
	\begin{theorem}\label{KST} (\cite{KST})
		For any sequence $\bar t$ satisfying \ref{x1}, \ref{x2} the graph $\bbG_0(\bar t) \subseteq [2^\omega]^2$ is a closed acyclic graph such that whenever $H \subseteq 2^\omega$ has the Baire property and $\bbG_0(\bar t)$-independent, then it must be meager.
		
		Moreover, if $G= (X,E)$ is an analytic graph on the Polish space $X$ and $\chi_{\text{B}}(G) > \aleph_0$, then there exists a continuous  map $f: 2^\omega \to X$, which is a homomorphism from $\bbG_0(\bar t)$ into $G$.
	\end{theorem}

	From now on we will only write $\mathbb{G}_0$ instead of $\mathbb{G}_0(\bar{\mathfrak{s}})$. Note that the graph $\bbG_0(\bar s)$ enjoys the expected properties if \ref{x1}, \ref{x2} holds without any regard to the specific sequence $\bar{\mathfrak{s}}$, justifying the use of the terminology $\chi_{\text{B}}(\bbG_0)$, $\chi_{\text{wB}}(\bbG_0)$, $\cov(I_{\bbG_0})$.
	
	\begin{definition}
		The graph $\mathbb{G}_{1}$ is defined as follows: 
		$$\mathbb{G}_{1} = \{ (x,y): \ x,y \in 2^\omega, \ \exists! n \in \omega \ x_n \neq y_n \} \subseteq [2^\omega]^2.$$
	\end{definition}

	\begin{definition}
		The Vitali relation $E_0$ is defined as follows: 
		$$E_{0} = \{ (x,y): \ x\neq y \in 2^\omega, \ \forall^\infty n \in \omega \ x_n = y_n \} \subseteq [2^\omega]^2.$$
	\end{definition}
	Note that this is not the standard definition of the Vitali relation, as we interpret it as a subset of $[2^\omega]^2$, while in the literature $E_0 \subseteq 2^\omega \times 2^\omega$ is an equivalence relation.
	
	\begin{definition}
		If $X$ is a topological space, and $G$ is a graph on it, then we let $I_G \subseteq \sP(X)$ denote the $\sigma$-ideal generated by Borel $G$-independent sets.
	\end{definition}
	Now we are ready to state our main theorem.
	\begin{theorem} \label{main}
		Assume $\mathbf{CH}$, and  let $\lambda_\mathbf{0} \leq \lambda_{\mathbf{1}} \leq \lambda_{\infty} \leq \lambda_\bbS = \kappa$ be infinite cardinals such that $ \lambda_\iota = \lambda_\iota^{\aleph_0}$ for each $\iota \in \{ \mathbf{0}, \mathbf{1}, \infty, \bbS\}$.
		Then in some cardinal preserving forcing extension we have
		\[ \begin{array}{rl}  \cov(\mathcal{M}) = \mathfrak{d} & = \aleph_1, \\ \cov(I_{\bbG_0}) & = \lambda_\mathbf{0}, \\
		 \cov(I_{\bbG_1}) & = \lambda_{\mathbf{1}}, \\
		   \cov(I_{E_0}) & = \lambda_{\mathbf{\infty}}, \\
		    \ 2^{\aleph_0} &  =  \kappa =\lambda_\bbS. \end{array} \]
	\end{theorem}
\begin{PROOF}{Theorem \ref{main}}{\ }
	We define our forcing posets in the following steps.
	
\begin{definition}\label{def1}{\ }
	\newcounter{enumD} \setcounter{enumD}{0}
\begin{enumerate}[label = $(\text{D}{\arabic*})$, ref = $\text{D}{\arabic*})$]
	\item \label{D1} For each $n$ we let $C_n \subseteq 2^{<\omega}$ be a finite, non-empty set such that for each $n$ 
			\begin{enumerate}
			\item $x \neq y \in C_n$ $\to$ $ x \nsubseteq y$ or $y \nsubseteq x$,
			\item  $C_{2n+1} = \{ \langle  0 \rangle , \langle 1 \rangle \}$, and
			\item for each $\bar t \in \prod_{i < 2n} C_i$ there exists $\bar t' \in C_{2n}$ such that $\bar t \tieconcat \bar t' = \mathfrak{s}_k$ (from \ref{x1}-\ref{x2}) for some $k \in \omega$.
		\end{enumerate}
	(This can be achieved by induction, e.g.  before constructing $C_{2n}$ imposing the auxiliary demand  that for some $r \in 2^\omega$ $\bar t \in C_{2n} \ \to \ r \notin [\bar t]$.)
	\item \label{P} Let $p \in \mathbb{P}^\mathbf{0}$, if 
	\begin{enumerate}[label = ($\roman*$)]
		\item $p = \langle p_i: \ i \in \omega \rangle$, where $\forall i: \  \emptyset \neq p_i \subseteq C_i$, and
		\item $\exists^\infty j$:  $p_{2j} = C_{2j}$ $\wedge$  $p_{2j+1} = C_{2j+1}$,
	\end{enumerate}  
	with $q$ stronger than $p$, (in symbols, $q \leq p$) iff $q_i \subseteq p_i$ for each $i$.
	\item \label{P1} For $\bbP^\mathbf{1}$ just recall the definition of the Silver real forcing: we let $p \in \mathbb{P}^\mathbf{1}$, if 
\begin{enumerate}[label = ($\roman*$)]
	\item $p = \langle p_i: \ i \in \omega \rangle$, where $\forall i: \  \emptyset \neq p_i \subseteq \{0,1\}$, and
	\item $\exists^\infty j$:  $p_{j} = \{0,1\}$,
\end{enumerate}  
with $q$ stronger than $p$, (in symbols, $q \leq p$) iff $q_i \subseteq p_i$ for each $i$.
\item  \label{Pinft} We let $p \in \mathbb{P}^\mathbf{\infty}$, if 
\begin{enumerate}[label = ($\roman*$)]
	\item $p = \langle p_i: \ i \in \omega \rangle$, where $\forall i \in \omega: \  \emptyset \neq p_i \subseteq \ ^{k_{p_i}} 2$ for some $k_{p_i} >0$, and
	\item $\exists^\infty j$:  $|p_{j}| = 2$,
\end{enumerate}  
with $q$ stronger than $p$, (in symbols, $q \leq p$) iff
\begin{itemize}
	\item there exists a strictly increasing infinite sequence $j_0 < j_1 < \dots$ of finite ordinals, for which
			\[(\forall i < \omega): \ k_{q_i} =  k_{p_{j_{i-1}+1}} + k_{p_{j_{i-1}+2}} + \dots + k_{p_{j_{i}}},\]
			where under $j_{-1}$ we mean $-1$,
	\item for each  $i < \omega$, $\bar t \in q_i$ there exists
	 $$(\bar t^*_{j_{i-1}+1} \in p_{j_{i-1}+1}) \ \& \ (\bar t^*_{j_{i-1}+2} \in p_{j_{i-1}+2}) \ \& \dots \& \ (\bar t^*_{j_{i}} \in p_{j_{i}}),$$
	  such that
	$$\bar t =  \bar t^*_{j_{i-1}+1} \tieconcat t^*_{j_{i-1}+2} \tieconcat \dots \tieconcat t^*_{j_{i}}.$$
	(So the biggest element $p$ can be described as the condition satisfying for each $j$ $p_j = \{ \langle 0 \rangle, \langle 1 \rangle\}$.)
\end{itemize}

	\item \label{cson} If $p \in \bbP^\iota$ ($\iota \in \{ \mathbf{0}, \mathbf{1}, \infty\}$), $\bar t \in \prod_{j<i} p_j$ for some $i <\omega$, then let $p^{[\bar t]}$
	denote the condition defined as $p^{[\bar t]}_j =  \{t_j\}$ for $j <i$, and $p^{[\bar t]} \rest [i,\omega) =  p \rest [i,\omega)$.
	\item  \label{csonk} Moreover, let $\bbP^{-m}$ denote the subforcing $\{ p \rest [m, \omega): \ p \in \bbP\}$ of $\bbP$ with the natural order.
	\stepcounter{enumD} \stepcounter{enumD}  \stepcounter{enumD} \stepcounter{enumD} \stepcounter{enumD} \stepcounter{enumD}
\end{enumerate}

Recall the definition of the Sacks forcing, and so let $\bbP^\bbS = \{T \subseteq 2^{<\omega}: \ T \text{ is a perfect tree} \}$ with  $ T \leq T'$ iff $T \subseteq T'$. 
\begin{definition}\label{def2} {\ }
	\begin{enumerate}[label = $(\text{D}{\arabic*})$, ref = $\text{D}{\arabic*})$]
		\setcounter{enumi}{\value{enumD}}
		\item \label{de}
	For $n \in \omega$  we define the partial order $\leq_n$ on $\bbP^\mathbf{0}$
	as $q \leq_n p$, iff
	\begin{itemize}
		\item $q \leq p$, and
		\item for some $\ell_0< \ell_1 < \dots < \ell_{n-1} < \omega$ we have
		 \[ \forall j <n: \ (p_{2\ell_j} = q_{2\ell_j} = C_{2\ell_j}) \&  (p_{2\ell_j+1} = q_{2\ell_j+1} = C_{2\ell_j+1}), \]
		 and
				\[ (\forall k < 2\ell_{n-1}) \	p_{k} = q_{k}.	 \]
 	\end{itemize}
 			\item \label{de1}
 	For $n \in \omega$  we define the partial order $\leq_n$ on $\bbP^\mathbf{1}$
 	as $q \leq_n p$, iff
 	\begin{itemize}
 		\item $q \leq p$, and
 		\item for some $l_0< l_1 < \dots <l_{n-1} < \omega$ we have
 		\[ (\forall k < l_{n-1}) \	p_{k} = q_{k}.	 \]
 		and
 			\[ \forall j <n: \ p_{l_j} =  q_{l_j} =  \{0,1\}, \]
 		
 	\end{itemize}
 	\stepcounter{enumD}
 	\item \label{deinf}
 	For $n \in \omega$  we define the partial order $\leq_n$ on $\bbP^\mathbf{\infty}$
 	as $q \leq_n p$, iff
 	\begin{itemize}
 		\item $q \leq p$, and
 		\item for some $l_0< l_1 < \dots l_{n-1} < \omega$ we have
 			\[ \forall j <n: \ |p_{l_j}| = |q_{l_j}| =2, \]
 			and
 		\[ (\forall k < l_{n-1}) \	p_{k} = q_{k},	 \]
 	
 	\end{itemize}
 	\stepcounter{enumD}
 	
 	\item \label{stemla} For $q \in \bbP^\bbS$ (so $q = T_q \subseteq \ ^{\omega>}2$ is a perfect tree) we define $\stem(q)$ to be the minimal branching node of $q$ (i.e. $q$ and stem($q$) satisfy $\stem(q) \tieconcat \langle 0 \rangle$, $\stem(q) \tieconcat \langle 1 \rangle \in  q$, but each proper initial segment of $\stem(q)$ has a unique immediate successor).
 
 	\item \label{dd} We define the partial order $\leq_n$ (for every $n \in \omega$) on $\bbS$
 			as 
 				\begin{itemize}
 				\item $q \leq_0 p$, iff $q \leq p$, 
 				\item $q \leq_{n+1} p$, iff  $q \leq p$,  $\stem(q) = \stem(p)$, and for this common stem $s$:
 				 $$q^{[s \tieconcat \langle 0 \rangle]} \leq_n p^{[s \tieconcat \langle 0 \rangle]}, \text{ and}$$
 						$$ q^{[s \tieconcat \langle 1 \rangle]} \leq_n p^{[s \tieconcat \langle 1 \rangle]}.$$ 
 				
  			\end{itemize}
  		
	\stepcounter{enumD}\stepcounter{enumD}\stepcounter{enumD}
\end{enumerate}
\end{definition}

A standard argument yields the following.
\begin{observation} \label{vanveg}
	$\bbP^\mathbf{0}$, $\bbP^\mathbf{1}$, $\bbP^\mathbf{\infty}$, $\bbP^\bbS$ satisfy Baumgartner's Axiom A with the partial orders defined above in \ref{de} and \ref{dd}, in particular if we are given the sequence 
	$$p_0 \geq_0 p_1 \geq_1 p_2 \geq_2 \dots \geq_{n-1} p_n \geq_n p_{n+1} \geq_{n+1} \dots, $$ 
	then there exists a common lower bound $p'$ w.r.t. $\geq$ (in fact even $p' \leq_n p_n$ can be assumed for each $n$).
\end{observation}

\begin{enumerate}[label = $(\text{D}{\arabic*})$, ref = $\text{D}{\arabic*})$]
		\setcounter{enumi}{\value{enumD}}
	\item For $I \subseteq \lambda_\mathbf{0}$ we let 
	$$\bbQ^\mathbf{0}_I = \{ f \in \ ^I(\bbP^\mathbf{0}): \ f(i) = 1_{\bbP^\mathbf{0}} \text{ for all, but countable } i's \}$$ be the countable support product of $\bbP$'s.
	\item For $I \subseteq [\lambda_\mathbf{0}, \lambda_\mathbf{1})$ let
	$$\bbQ^\mathbf{1}_I = \{ f \in \ ^I(\bbP^\mathbf{1}): \ f(i) = 1_{\bbP^\mathbf{1}} \text{ for all, but countable } i's \}$$ be the countable support product of $\bbP^\mathbf{1}$'s.
	\item Similarly for $I \subseteq [\lambda_\mathbf{1}, \lambda_\mathbf{\infty})$ let
	$$\bbQ^\mathbf{\infty}_I = \{ f \in \ ^I(\bbP^\mathbf{\infty}): \ f(i) = 1_{\bbP^\mathbf{\infty}} \text{ for all, but countable } i's \}$$ be the countable support product of $\bbP^\mathbf{\infty}$'s,
	\item and for $I \subseteq [\lambda_\mathbf{\infty}, \lambda_\mathbf{\bbS})$ let
	$$\bbQ^\mathbf{\bbS}_I = \{ f \in \ ^I(\bbP^\mathbf{\bbS}): \ f(i) = 1_{\bbP^\mathbf{\bbS}} \text{ for all, but countable } i's \}$$ be the countable support product of $\bbP^\mathbf{\bbS}$'s.	
	
	\item We let $\bbQ$ be the following countable support product:
	\[ \bbQ = \bbQ^0_{\lambda_\mathbf{0}} \times \bbQ^\mathbf{1}_{\lambda_\mathbf{1} \setminus \lambda_\mathbf{0}} \times \bbQ^\mathbf{\infty}_{\lambda_\mathbf{\infty} \setminus \lambda_\mathbf{1}} \times \bbQ^{\bbS}_{\lambda_\bbS \setminus \lambda_{\infty}} .\]
	\stepcounter{enumD}\stepcounter{enumD}\stepcounter{enumD}
\end{enumerate}
\end{definition}
We have to check that the forcing $\bbQ$ is indeed cardinal preserving, forcing $\mathfrak{d}$ to be $\aleph_1$, the continuum to be $\kappa = \lambda_\bbS$, there exists a system of $\lambda_\mathbf{0}$-many $\bbG_0$-independent Borel sets covering $2^\omega$, but any system of smaller cardinality is not sufficient (and similarly for $\bbG_1$, and $E_0$).
For this we will prove the following:
\begin{enumerate}[label = $(\circledast)_{\arabic*}$, ref =  $(\circledast)_{\arabic*}$]
	\item \label{e0} $\bbQ$ is proper, and has the $\aleph_2$-cc,
	\item \label{bou} $\bbQ$ is $\omega^\omega$-bounding, i.e. for each $r \in \ ^\omega \omega \cap V^\bbQ$ there exists $r' \in \ ^\omega \omega \cap V$ such that $r' \geq^* r$,
	\item  \label{ee2} for each $r \in \ ^\omega 2 \cap V^\bbQ$ there exists a 
	\begin{itemize}
		\item  tree $T_\mathbf{0} \in V^{\bbQ^\mathbf{0}_{\lambda_\mathbf{0}}} \cap \mathscr{P}(2^{<\omega})$  such that $r \in [T_\mathbf{0}]$, and  $[T_\mathbf{0}]$ is $\bbG_0$-independent, 
		\item tree $T_\mathbf{1} \in V^{\bbQ^\mathbf{0}_{\lambda_\mathbf{0}} \times \bbQ^\mathbf{1}_{\lambda_\mathbf{1} \setminus \lambda_\mathbf{0}}} \cap \mathscr{P}(2^{<\omega})$  such that  $r \in [T_\mathbf{1}]$, and $[T_\mathbf{1}]$ is $\bbG_1$-independent, 
		\item tree $T_\mathbf{\infty} \in V^{\bbQ^\mathbf{0}_{\lambda_\mathbf{0}} \times \bbQ^\mathbf{1}_{\lambda_\mathbf{1} \setminus \lambda_\mathbf{0}} \times \bbQ^\mathbf{\infty}_{\lambda_\mathbf{\infty} \setminus \lambda_\mathbf{1}}} \cap \mathscr{P}(2^{<\omega})$  such that  $r \in [T_\mathbf{\infty}]$, and $[T_\mathbf{\infty}]$ is $E_0$-independent, 
	\end{itemize}

\item \label{e3} If $\alpha < \lambda_\bbS$, $T  \in V^{\bbQ \rest \lambda_\bbS \setminus \{ \alpha\}}  \cap \sP(^{\omega>}\{0,1\})$ is a tree, such that 
	\begin{itemize}
	\item either $\alpha < \lambda_\mathbf{0}$, and $[T]$ is $\bbG_0$-independent,
	\item or $\alpha \in [\lambda_\mathbf{0}, \lambda_\mathbf{1})$, and $[T]$ is $\bbG_1$-independent,
	\item  $\alpha \in [\lambda_\mathbf{1}, \lambda_\mathbf{\infty})$, and $[T]$ is $E_0$-independent,
\end{itemize}
 then for $r_\alpha$, the generic real given by the $\alpha$'th coordinate ($\bbQ \rest \{ \alpha \}$), we have:
	\[ r_\alpha \notin [T].\]
\end{enumerate}
Similar statements to \ref{ee2} are proved independently in \cite[§3]{weabo} (for an extension adding a single real, and generalizing it to CS iterations), and see also  \cite[Thm 3.47- Corollary 3.49]{zapletal2019hypergraphs}, a more general result, albeit only for a single step extension, which is independent of both.

Observe that the properness of $\bbQ$ (together with our assumptions on the ground model) would imply that $(2^{\aleph_0})^{V^{\bbQ \rest \lambda_\mathbf{\iota}}} = \lambda_\iota$ ($\iota \in \{\mathbf{0}, \mathbf{1}, \infty, \bbS\}$), and so  $(2^{\aleph_0})^{V^\bbQ} = \kappa^{\aleph_0} = \kappa$ ($=\lambda_\bbS$), while it follows from \ref{ee2} that cov($I_{\bbG_0})\leq \lambda_\mathbf{0}$ (and the respective inequalities similarly hold for $\bbG_1$, and $E_0$).
 
 By the $\aleph_2$-cc of $\bbQ$, if $\mu$ is uncountable, then each system of Borel sets of size $\mu$ is in 
 $V^{\bbQ \rest M}$ for some $M \in V$ of size at most $\mu \cdot \aleph_1 = \mu$.
 
Moreover, \ref{bou} clearly implies $\mathfrak{d}^{V^\bbQ} = \aleph_1$, and the inequality $\cov(\mathcal{M}) \leq \mathfrak{d}$ holds in $\mathbf{ZFC}$ (since each compact set in $\bbN^\bbN$ is meager as well, or see \cite{BaJu95}).
Finally, since each uncountable Borel subset of a Polish space is a continuous image of (a closed subspace of) $\bbN^\bbN$ \cite[Theorems 7.9, 13.1]{kekrisz}, each Borel set is the union of $\mathfrak{d}$-many compact sets. Thus, if a Polish space can be covered by $\aleph_1 \leq \mu$-many  Borel sets, then we can replace each Borel set $B$ with a system of $\aleph_1$-many compact sets $\langle K_\alpha: \ \alpha < \omega_1 \rangle$ with $B = \bigcup_{\alpha < \omega_1} K_\alpha$, so (in the extension) there exists $\langle K'_\alpha: \ \alpha < $cov($I_{\bbG_0}$)$\rangle$ with each $K'_\alpha$  compact, covering $2^\omega$. This together with \ref{e3} implies cov($I_{\bbG_0}$)$\geq \lambda_\mathbf{0}$, cov($I_{\bbG_1}$)$\geq \lambda_\mathbf{1}$, cov($I_{E_0}$)$\geq \lambda_\mathbf{\infty}$.
 Therefore it is indeed enough to verify clauses \ref{e0}- \ref{e3}.
 
 \begin{claim}\label{cc}
 	$\bbQ$ has the $\aleph_2$-cc.
 \end{claim}
\begin{PROOF}{Claim \ref{cc}}
	Suppose that $\langle a_i: \ i < \omega_2 \rangle$ is an antichain. Since $|\supp(a_i)| \leq \aleph_0$ for $i < \omega_2$, by $\mathbf{CH}$ we can assume that $\{\supp(a_i): \ i < \omega_2 \}$ forms a $\Delta$-system with kernel $K$. But $|\bbP^\mathbf{0}| = |\bbP^\mathbf{1}| = |\bbP^\mathbf{\infty}| = |\bbP^\bbS| = |(2^{\aleph_0})^V| = \aleph_1$ (by Definition \ref{def1}), so $|^K(\bbP^\mathbf{0} \cup \bbP^\mathbf{1} \cup \bbP^\mathbf{\infty} \cup \bbP^\bbS)| = \aleph_1$, we are done.
\end{PROOF}

	\begin{convention}
		By passing down  to a dense subset of $\bbP^\mathbf{0}$, from now on we can assume that whenever $p \in \bbP^\mathbf{0}$, $k \in \omega$,
		$$ \neg (p_{2k} = C_{2k} \ \wedge \ p_{2k+1} = C_{2k+1}) \ \rightarrow \ |p_{2k}| = |p_{2k+1}| = 1.$$
	\end{convention}
\begin{definition} \label{megs} {\ } 
	\begin{enumerate}[label = $(\text{D}^*{\arabic*})$, ref = $\text{D}^*{\arabic*})$]
		\item  if $p \in \bbP^\mathbf{0}$, $n \in \omega$, and $d$ is the smallest integer for which

			 for some $l_0< l_1 < \dots < l_{n-1} < d$ we have
			\[ \forall j <n: \ (p_{2l_j} = q_{2l_j} = C_{2l_j}) \&  (p_{2l_j+1} = q_{2l_j+1} = C_{2l_j+1}), \]
		then we let
		$$ T_n(p) = \prod_{j<2d} p_{j}, $$ 
		and $$  T(p) = \bigcup_{n \in \omega} T_n(p), $$ 
		\item if $p \in \bbP^\mathbf{0}$, $\bar u \in \bigcup_{n \in \omega} \prod_{j<n} p_j$, (e.g. $\bar u \in T(p)$), then we define $p^{[\bar u]}= p^{(\bar u)}$ , to be a condition in $\bbP^0$ satisfying
		\begin{itemize}
			\item $p^{[\bar u]}_j = p^{(\bar u)}_j = \{u_j\}$, if $j <\lh(\bar u)$, 
			\item $p^{[\bar u]}_j =p^{(\bar u)}_j = p_j$, if $j \geq \lh(\bar u)$.
		\end{itemize}
		
	\end{enumerate}

	\begin{enumerate}[label = $(\text{D}^*{\arabic*})$, ref = $\text{D}^*{\arabic*})$]
		\setcounter{enumi}{2}
		\item \label{d2''} for $n \in  \omega$, $ p \in \bbP^\mathbf{1}$, if $d < \omega$ is the minimal natural number such that for some $l_0 < l_1 < \dots < l_{n-1}  <d$ 
		\[ \forall j <n: \ \left(p_{l_j} =  \{0,1\}\right), \]
 then for  each $\bar v = \langle v_0, v_1, \dots, v_{n-1} \rangle\in \ ^n2$ we define the sequence $\bar{\mathfrak{t}}_{\bar v}(p) = \langle \mathfrak{t}_{\bar v}(p)_j: \ j <d \rangle  \in \ ^d2$ as
		\begin{itemize}
			\item $\mathfrak{t}_{\bar v}(p)_{l_j} = v_j$,
			\item $\mathfrak{t}_{\bar v}(p)_{k} = a_k$, where $p_k = \{a_k\}$, $k \in (l_j, l_{j+1})$ for some $j<n-1$.
		\end{itemize}

		\item if $p \in \bbP^\mathbf{1}$, $\bar t \in \bigcup_{n \in \omega} \prod_{j<n} p_j$, (e.g. $\bar t = \bar{\mathfrak{t}}_{\bar u}(p)$ for some $\bar u \in \ ^{\omega>}2$), then we define $p^{[\bar t]}$, to be a condition in $\bbP^{\mathbf{1}}$ satisfying
	\begin{itemize}
		\item $p^{[\bar t]}_j = \{t_j\}$, if $j <\lh(\bar t)$, 
		\item $p^{[\bar t]}_j = p_j$, if $j \geq \lh(\bar t)$.
	\end{itemize}
		Moreover, if $\bar u \in \ ^n2$ for some $n$, then we let 
		$$p^{(\bar u)}:= p^{[ \bar{\mathfrak{t}}_{\bar u}(p)]}.$$
	
	\item \label{frakting}
	for $n \in  \omega$, $ p \in \bbP^\mathbf{\infty}$ if $d < \omega$ is the minimal natural number such that for some $l_0 < l_1 < \dots < l_{n-1}  <d$ 
	\[ \forall j <n: \ \left(|p_{l_j}| =  2 \right), \]
	then for  each $\bar v = \langle v_0, v_1, \dots, v_{n-1} \rangle\in \ ^n2$ we define the sequence $\bar{\bar{\mathfrak{t}}}_{\bar v}(p) = \langle \bar{\mathfrak{t}}_{\bar v}(p)_j: \ j <d \rangle  \in \ ^d2$ as
	\begin{itemize}
		\item $\bar{\mathfrak{t}}_{\bar v}(p)_{l_j} =  \bar t^{p_k}_{l_j, w_j}$, where $p_{l_j} = \{ \bar t^{p_k}_{l_j, 0} <_{\textrm{lex}} \bar t^{p_k}_{l_j, 1}\}$,
		\item $\bar{\mathfrak{t}}_{\bar v}(p)_{k} =  \bar{t}^{p_k}$, where $p_k = \{\bar t^{p_k}\}$, $k \in (l_j, l_{j+1})$ for some $j<n-1$.
	\end{itemize}

	\item if $p \in \bbP^\mathbf{\infty}$, $\bar{\bar t} \in \bigcup_{n \in \omega} \prod_{j<n} p_j$, (e.g. $\bar{\bar{t}} \in T(p)$), then we define $p^{[\bar{\bar t}]}$, to be a condition in $\bbP^\infty$ satisfying
	\begin{itemize}
		\item $p^{[\bar{\bar{t}}]}_j = \{\bar t_j\}$, if $j <\lh(\bar{\bar{t}})$, 
		\item $p^{[\bar{\bar{t}}]}_j = p_j$, if $j \geq \lh(\bar{\bar{t}})$.
	\end{itemize}
	Furthermore, if $\bar w \in \ ^n2$ for some $n$, then we let 
	$$p^{(\bar w)}:= p^{[ \bar{\bar{\mathfrak{t}}}_{\bar w}(p)]}.$$

	\item \label{frakt}
	For $p \in \bbP^\bbS$, $\bar s \in \ ^n2$ we define the node $\bar{\mathfrak{t}}_{\bar s}(p) \in p$ by induction on $\lh(\bar s)$ as follows.
	Let 
	\begin{equation} \bar{\mathfrak{t}}_{\langle\rangle}(p) = \langle \rangle, \end{equation} and
	\begin{equation} 
		\begin{array}{rl} \bar{\mathfrak{t}}_{\bar s \tieconcat \langle 0 \rangle}(p) = & \stem\left(p^{[\bar{\mathfrak{t}}_{\bar s}(p)]}\right) \tieconcat \langle 0 \rangle, \\ 
			\bar{\mathfrak{t}}_{\bar s \tieconcat \langle 1 \rangle}(p) = &  \stem\left(p^{[\bar{\mathfrak{t}}_{\bar s}(p)]}\right) \tieconcat \langle 1 \rangle, \end{array}
	\end{equation} 
	(recall \ref{stemla},  for the stem of a tree $T$ is the unique largest element in 
	$$\{t \in T: \ \forall t' \in T: \ t \subseteq t' \ \vee \ t' \subseteq t\}).$$
	Note that $q \leq_n p$, iff $q \leq p$, and for each $s \in \ ^n2$ we have $\bar{\mathfrak{t}}_s(p) = \bar{\mathfrak{t}}_s(q)$.
	
		\item if $p \in \bbP^\mathbf{\bbS}$ and $\bar s \in p$, then we define $p^{[\bar s]}$, to be the condition $p^{[\bar s]} = \{\bar t \in p: \ \bar t \supseteq \bar s\}$, and for arbitrary $\bar s \in \ ^{\omega>}2$ we let 
		$$p^{(\bar s)} = p^{[\bar{\mathfrak{t}}_{\bar s}(p)]}.$$

\end{enumerate} 
\end{definition} 

\begin{observation} \label{ossze}
	If $p \in \bbP^\infty$, $i \in \omega$,
	\begin{equation} \label{p} p_0 = \{\bar t^p_0\}, p_1 = \{\bar t^p_1\},  p_{i-1}= \{\bar t^p_{i-1}\}, \end{equation} and we consider the condition $p' \leq p$ defined as
	$$\begin{array}{l}
		p'_0 = \{ \bar t^{p'}_0 = \bar t^p_0 \tieconcat \bar t^p_1 \tieconcat \dots \tieconcat  \bar t^p_{i-1} \}, \\
		p'_{j+1} = p_{j+i} \ (\text{for }j \in \omega)
		\end{array} $$
		then $p  \Vdash_{\bbP^\infty}  p' \in \name{\mathbf{G}}$.
\end{observation}
\begin{PROOF}{Observation \ref{ossze}}
	Suppose that $p^* \leq p$, and $p^* \perp p'$. By further strengthening $p^*$ w.l.o.g. we can assume that $|p^*_0| = 1$, and for the unique $\bar t^* \in p^*_0$ we have 
	$$\lh(\bar t^*) > |\bar t^p_0| + |\bar  t^p_1| + \dots + |\bar t^p_{i-1}|.$$ 
	But then this bound together with
	$p^* \leq p$ and \eqref{p} implies that
	$$\bar t^* = (\bar t^p_0 \tieconcat \bar t^p_1 \tieconcat \dots \tieconcat  \bar t^p_{i-1}) \tieconcat t^p_{i} \tieconcat \dots t^p_{i+k} = \bar t^{p'}_0 \tieconcat \bar t^p_{i} \tieconcat \dots \bar t^p_{i+k}$$  
	for some $\bar t^p_{i} \in p_i = p'_1$, $\bar t^p_{i+1} \in p_{i+1} = p'_2$, $\dots$, $\bar t^p_{i+k} \in p_{i+k} = p'_{k+1}$, so not only is $p'$ compatible with $p^*$, but $p^* \leq p'$.
	
\end{PROOF}

		\newcounter{elokcc} \setcounter{elokcc}{0}
		
		\begin{definition} \label{ELL}
		We let $\bar{\bar{\ell}}\in \mathcal{L}$, iff
		$\bar{\bar{\ell}} = (\bar \ell^\mathbf{0}, \bar \ell^\mathbf{1}, \bar \ell^\mathbf{\infty}, \bar \ell^\bbS)$, where 
		\begin{itemize}
			\item $\bar \ell^\bbS = \langle \ell^\bbS_i: \ i< \omega \rangle \in \ ^{\omega}\omega$, with  $\sum_{i< \omega} \ell^\bbS_i < \infty$, 
			\item $\bar \ell^\mathbf{0} = \langle \ell^\mathbf{0}_i: \ i \in \omega \rangle \in \ ^{\omega}\omega$, with $\sum_{i< \omega} \ell^\mathbf{0}_i < \infty$, 
			\item $\bar \ell^\mathbf{1} = \langle \ell^\mathbf{1}_i: \ i \in \omega \rangle \in \ ^{\omega}\omega$, with $\sum_{i< \omega} \ell^\mathbf{1}_i < \infty$, 
			\item $\bar \ell^\mathbf{\infty} = \langle \ell^\mathbf{\infty}_i: \ i \in \omega  \rangle \in \ ^{\omega}\omega$, with  $\sum_{i< \omega} \ell^\mathbf{\infty}_i < \infty$, 
		\end{itemize}
	\end{definition}	
	\begin{definition} \label{D()}
	Fix $\langle \varp^\mathbf{0}_j: \ j \in \omega  \rangle\in \ ^\omega \lambda_\mathbf{0}$, $\langle \varp^{\mathbf{1}}_j: \ j \in \omega  \rangle \in \ ^\omega [\lambda_\mathbf{0}, \lambda_\mathbf{1})$,  $\langle \varp^{\mathbf{\infty}}_j: \ j \in \omega \rangle \in  \ ^\omega [\lambda_\mathbf{1}, \lambda_\infty)$,  $\langle \varp^\bbS_j: \ j \in \omega \rangle \in \ ^\omega [\lambda_\infty, \lambda_\bbS)$.
	The following symbols depend on the particular fixed $\bar \varp^\iota$'s, but we omit it as it will be always clear what those sequences are.

	Now 
	\begin{enumerate}[label = $(\blacktriangleright_{\arabic*})$, ref = $(\blacktriangleright_{\arabic*})$]
	

		\item \label{d2'} for $j \in  \omega$, $n\in \omega$ let
		 $$T^j_n(q) = T_n(q(\varp^\mathbf{0}_j))$$
		 	and
		 $$ T^j(q) = \bigcup_{n \in \omega} T^j_n(q), $$
		 \item 	Suppose that $\bar{\bar{\ell}} \in \mathcal{L}$, $q \in \bbQ'$. Then $(\bar u, \bar v, \bar w, \bar s) \in \mathbf{seq}_{\bar{\bar{\ell}}}(q)$, iff 
		 \begin{itemize} 
		 	\item $\bar s \in \ ^{\omega}(^{<\omega}2)$ with $s_j \in \ ^{\ell^\bbS_j}2$ ($j < \omega$),
		 	\item  $\bar u  \in \ \prod_{j < \omega}T^j(q)$ with $u_j \in T^j_{\ell^\mathbf{0}_j}(q)$ ($j < \omega$),
		 	\item  $\bar v \in \ ^{\omega}(^{<\omega}2)$ with $v_j \in \ ^{\ell^\mathbf{1}_j}2$ ($j < \omega$),
		 	\item  $\bar w \in \ ^{\omega}(^{<\omega}2)$ with $w_j \in \ ^{\ell^\mathbf{\infty}_j}2$ ($j < \omega$).
		 \end{itemize}


		\item  \label{()} For each
		$$\begin{array}{ll}
			\bar u = \langle u_j: \ j \in \omega  \rangle & \in \prod_{j \in \omega} T^{j}(q)), \\
			\bar v = \langle v_j: \ j \in \omega \rangle & \in \prod_{j \in \omega} \left( ^{\omega>}2 \right), \\
			\bar w = \langle v_j: \ j \in \omega \rangle & \in \prod_{j \in \omega} \left( ^{\omega>}2 \right),	\\
			\bar s = \langle s_j: \ j \in \omega \rangle & \in \prod_{j \in \omega} \left( ^{\omega>}2 \right),
			
		\end{array}$$
		we let $q^{(\bar u, \bar v, \bar w, \bar s)} \in \bbQ$ be defined as
		\begin{itemize}
			\item $q^{(\bar u, \bar v, \bar w, \bar s)}(\varp^\mathbf{0}_j) = q(\varp^\mathbf{0}_j)^{(u_j)}$ ($j \in \omega$),
			\item $q^{(\bar u, \bar v, \bar w, \bar s)}(\varp^\mathbf{1}_j) = q(\varp^\mathbf{1}_j)^{(v_j)}$  ($j \in \omega$),
			\item  $q^{(\bar u, \bar v, \bar w, \bar s)}(\varp^\mathbf{\infty}_j) = q(\varp^\infty_j)^{(w_j)}$  ($j \in \omega$),
			\item  $q^{(\bar u, \bar v, \bar w, \bar s)}(\varp^\bbS_j) = q(\varp_j^\bbS)^{(s_j)}$  ($j \in \omega$),
		\end{itemize}
		\stepcounter{elokcc}
	\end{enumerate}
	

\end{definition}

\begin{definition} \label{seq}

\end{definition}
Observe that 
\begin{equation} \label{seqfi} \text{whenever } \bar{\bar{\ell}} \in \mathcal{L}: \ |\mathbf{seq}_{\bar{\bar{\ell}}}(q)| < \aleph_0, \end{equation}
since $\bar{\bar{\ell}}$ has a finite sum.

	


	\begin{definition}
	Assuming the sequences $\bar \varp^\iota$ are as in Definition \ref{D()} ($\iota \in \{\mathbf{0}, \mathbf{1}, \infty, \bbS\}$, and $\bar{\bar{\ell}} \in \mathcal{L}$, we let the partial order $\leq_{\bar{\bar{\ell}}}$ defined by
	\[ \begin{array}{rl}
		p \leq_{\bar{\bar{\ell}}} q \ \iff & \bullet_1 \ (p \leq q) \ \wedge \\ & \bullet_2 \ \mathbf{seq}_{\bar{\bar{\ell}}}(p) = \mathbf{seq}_{\bar{\bar{\ell}}}(q) \ \wedge \\ 
		& \bullet_3 \ \forall (\bar u, \bar v, \bar w, \bar s) \in \mathbf{seq}_{\bar{\bar{\ell}}}(q): \ p^{(\bar u, \bar v, \bar w, \bar s)} \leq q^{(\bar u, \bar v, \bar w, \bar s)}.
	\end{array}
	\]
\end{definition}

		Note the following easy corollaries of our definitions:
	\begin{observation} \label{sv}
		Let $\bar \varp^\iota$  ($\iota \in \{\mathbf{0}, \mathbf{1}, \infty, \bbS\}$) be as in Definition \ref{D()}, and $\bar{\bar{\ell}} \in \mathcal{L}$ be given.
		
		Then, if $p \geq q \in \bbQ$ holds, then $p \geq_{\bar{\bar{\ell}}} q$, iff 
		\begin{itemize}
			\item for each $j \in \omega$: $T^j_{\ell^\mathbf{0}_j}(p) = T^j_{\ell^\mathbf{0}_j}(q)$, and
			\item  	for each $j\in \omega$ and $u \in \ ^{\ell^\mathbf{1}_j}2$ we have $\mathfrak{t}^{\mathbf{1},j}_u(p) = \mathfrak{t}^{\mathbf{1},j}_u(q)$, and
			\item  for each $j\in \omega$ and $v \in \ ^{\ell^\mathbf{\infty}_j}2$: $\mathfrak{t}^{\mathbf{\infty},j}_v(p) = \mathfrak{t}^{\mathbf{\infty},j}_v(q)$, and
			\item 	for each $j\in \omega$ and $s \in \ ^{\ell^\bbS_j}2$: $\mathfrak{t}^{\bbS,j}_s(p) = \mathfrak{t}^{\bbS,j}_s(q)$.
		\end{itemize}
	\end{observation}

		\begin{observation} \label{kit}  \label{obskite}
		If $p, q \in \bbQ$, $\bar{\bar{\ell}} \in \mathcal{L}$, $\bar \varp^\iota$ ($\iota \in \{\mathbf{0}, \mathbf{1}, \infty, \bbS\}$) are as in Definition \ref{D()},  $(\bar u, \bar v, \bar w, \bar s)  \in \mathbf{seq}_{\bar{\bar{\ell}}}(p)$, and $q \leq p^{(\bar u, \bar v, \bar w, \bar s)}$ then for some $q_* \leq_{\bar{\bar{\ell}}} p$, 
		$$q_*^{(\bar u, \bar v, \bar w, \bar s)} \Vdash  q \in \mathbf{G},$$
		even $q_*^{(\bar u, \bar v, \bar w, \bar s)} = q$, if in addition $p^{(\bar u, \bar v, \bar w, \bar s)}(\varp^\infty_j) \geq_{\ell^\infty_j} q(\varp^\infty_j)$ holds for each $j$.

		Moreover, we can assume that whenever $m \in \omega$, and $s' \in \ ^{\omega>}2$ is not comparable with $s_m$ (i.e.\ $s' \nsubseteq s_m$, $s'\nsupseteq s_m$), then $q_*(\varp^{\mathbf{0}}_m) \in \bbP^{\mathbf{0}}$ satisfies $q_*(\varp^{\mathbf{0}}_m)^{(s')} = p(\varp^{\mathbf{0}}_m)^{(s')}$.
	\end{observation}
		Note that we cannot expect above $q_*^{(\bar u, \bar v, \bar w, \bar s)} =  q$ to hold in general, since on coordinates of the form $\varp^\infty_j$ possibly $q(\varp^\infty_j)_0 = \{\bar t\}$, where $\bar t =  (\bar t^*)^0 \tieconcat  (\bar t^*)^1$ with $ (\bar t^*)^0 \in q_*^{(\bar u, \bar v, \bar w, \bar s)}(\varp^\infty_j)_0$,  and $ (\bar t^*)^1 \in q_*^{(\bar u, \bar v, \bar w, \bar s)}(\varp^\infty_j)_1$.
		
		
	\begin{observation}\label{obski}
		If $p \geq q \in \bbQ$, $\bar{\bar{\ell}} \in \mathcal{L}$, $(\bar u, \bar v, \bar w, \bar s)  \in \mathbf{seq}_{\bar{\bar{\ell}}}(q)$, then $q^{(\bar u, \bar v, \bar w, \bar s)} \leq p^{(\bar u', \bar v', \bar w', \bar s')}$ for some $(\bar u', \bar v', \bar w', \bar s')  \in \mathbf{seq}_{\bar{\bar{\ell}}}(p)$.
	\end{observation}
	
	The next claim verifies the properness part of \ref{e0}, and \ref{bou}.
	\begin{claim} \label{fuzio}
		Let $q \in \bbQ$, $D_0,D_1, \dots, D_i, \dots$ be a countable sequence of maximal antichains of $\bbQ$.
		Then for a suitable extension $q' \leq q$ we have that for each $i \in \omega$ $q'$ is compatible with only finitely many elements of $D_i$.
	\end{claim}
	\begin{PROOF}{Claim \ref{fuzio}}
		Assume that $q \in \bbQ$, and the $D_j$'s are fixed. In what follows we will 
		sketch a standard fusion argument for Baumgartner's Axiom A. 
		

		The following is a trivial application of Observation $\ref{vanveg}$:
		\begin{observation} \label{ob2}
			Suppose that the sequence $\langle q_n: \ n \in \omega \rangle \in \ ^\omega\bbQ$ is decreasing, the sequences $\bar \varp^\iota$ ($\iota \in \{\mathbf{0}, \mathbf{1}, \infty, \bbS \}$) are as in Definition \ref{D()}, and for each $k$ there is $\bar{\bar{\ell}}^k = ((\bar \ell^k)^\mathbf{0}, (\bar \ell^k)^\mathbf{1},  (\bar \ell^k)^\mathbf{\infty},  (\bar \ell^k)^\bbS) \in \mathcal{L}$ such that 
			\begin{itemize} 
				\item for each $\alpha \in \bigcup_{n \in \omega} \supp(q_n)$ we have that $\alpha = \varp^\iota_j$ for some $\iota$ and $j$,
				and  
				$$ \langle (\ell^k)^\iota_j: \ k \in \omega\rangle  \text{ is nondecreasing, converging to } \infty,$$
				\item $q_{n+1} \leq_{\bar{\bar{\ell}}^n} q_n$ holds for each $n$.
			\end{itemize} 
			Then there exists a common lower bound $q_\omega \in \bbQ$ of the sequence $\langle q_n: \ n \in \omega \rangle$.
			Moreover, there exists $q_\omega$ for which for each $n$ $q_\omega \leq_{\bar{\bar{\ell}}^n} q_n$.
		\end{observation}
		We will define the sequences $\langle q_i: \ i < \omega \rangle$, $\langle \varp^\mathbf{0}_i: \ i \in \omega \rangle$, $\langle \varp^\mathbf{1}_i: \ i \in \omega \rangle$, $\langle \varp^\infty_i: \ i < \omega \rangle$ $\langle \varp^\bbS_i: \ i < \omega \rangle$, $\langle \bar{\bar{\ell}}^i: \ i \in \omega \rangle$  satisfying the following:
		\begin{enumerate}[label = $(\XBox_{\arabic*})$, ref = $(\XBox_{\arabic*})$]
			\item \label{h1} $q_0 = q$, and for each $i$ we have $q_i \in \bbQ$,
			\item $\{ \varp^\mathbf{0}_i: \ i \in \omega \} \subseteq \lambda_{\mathbf{0}}$, $\{ \varp^\mathbf{1}_i: \ i \in \omega \} \subseteq \lambda_{\mathbf{1}} \setminus \lambda_\mathbf{0}$, $\{\varp^\infty_i: \ i < \omega \} \subseteq \lambda_\infty \setminus \lambda_\mathbf{1}$,  $\{\varp^\bbS_i: \ i < \omega \} \subseteq \lambda_\bbS \setminus \lambda_\infty$,
			\item \label{Sn} for each $n$ $\supp(q_{n}) \subseteq \{\varp^\iota_j: \ \iota \in \{\mathbf{0}, \mathbf{1}, \infty, \bbS\}, \ j \in \omega\}$,
			\item for each  $\iota \in \{\mathbf{0}, \mathbf{1}, \infty, \bbS\}$ and $j \in \omega$ the sequence 
			$$\langle (\ell^n)^{\iota}_j: \ n \in \omega \rangle \ \text{ is nondecreasing, and tends to } \infty,$$ 
			\item $\forall n$ $ q_n \geq_{\bar{\bar{\ell}}^n} q_{n+1}$,
			\item $\forall n$ the condition $q_{n+1}$ is  compatible with only finitely many conditions in $D_n$.
		\end{enumerate}
		Provided that such sequences exist we can appeal to Observation $\ref{ob2}$, which will complete the proof of Claim $\ref{fuzio}$.

	We can clearly define a sequence of $\bar{\bar{\ell}}^i$'s as in Observation \ref{ob2}.
	Now by Observation $\ref{sv}$, and some standard bookkeeping arguments it is easy to see that the entire induction can be done once we specify how to define the condition $q_{n+1} \in \bbQ$ from $q_n$ and the adequate fragment of $\bar \varp^\iota$'s.
This $q_{n+1}$ will satisfy that
		\newcounter{elokc} \setcounter{elokc}{0}
	\begin{enumerate}[label = ($\blacktriangle_{\arabic*}$), ref = ($\blacktriangle_{\arabic*}$)]
		\item \label{e1} $q_{n+1} \leq_{\bar{\bar{\ell}}^n} q_n$,
		\item \label{e2} whenever $(\bar u, \bar v, \bar w, \bar s) \in \mathbf{seq}_{\bar{\bar{\ell}}^n}(q_n)$,
		 then $(q_{n+1})^{(\bar u, \bar v, \bar w, \bar s)}$ is compatible with exactly one element of $D_n$.
	\end{enumerate}
	For this	\begin{enumerate}[label = $(\circledcirc_{\arabic*})$, ref = $(\circledcirc_{\arabic*})$]
		\setcounter{enumi}{\value{elokc}}
		\item let
	$$M = |\mathbf{seq}_{\bar{\bar{\ell}}^n}(q_n)|,$$ 
	and 
	fix an enumeration 
	\begin{equation}\label{vv}
		\langle (\bar u^i, \bar v^i, \bar w^i, \bar s^i):  \ i <M \rangle
	\end{equation} of $\mathbf{seq}_{\bar{\bar{\ell}}^n}(q_n)$.
	\stepcounter{elokc}
	\end{enumerate}
	Note that \ref{e2} includes $M$-many different objectives, each one is corresponding to some $(\bar u^i, \bar v^i, \bar w^i, \bar s^i)$  from $\eqref{vv}$.
	So
		\begin{enumerate}[label = $(\circledcirc_{\arabic*})$, ref = $(\circledcirc_{\arabic*})$]
		\setcounter{enumi}{\value{elokc}}
		\item \label{Mth} we construct the sequence $\langle q^*_i: \ i \leq M \rangle$
		satisfying
		\[ q^*_0 = q_n \geq_{\bar{\bar{\ell}}^n} q^*_1 \geq_{\bar{\bar{\ell}}^n} \dots \geq_{\bar{\bar{\ell}}^n} q^*_M,\]
		and 
		\[ (\forall i<M): \ (q^*_{i+1})^{(\bar u^i, \bar v^i, \bar w^i, \bar s^i)} \Vdash p_* \in \mathbf{G}, \ \text{ for some } p_* \in D_n,\]
		 thus $q^*_M$ will work (i.e.\ \ref{e1}, \ref{e2} hold).
		\stepcounter{elokc}
	\end{enumerate}

	Assuming that $i <M$ and $q^*_i$ is defined, 
	pick $q' \leq (q^{*}_i)^{(\bar u^i, \bar v^i, \bar w^i, \bar s^i)}$, such that $q' \leq p_*$ for some $p_* \in D_n$. Let 
	$q^*_{i+1} \leq_{\bar{\bar{\ell}}^n} q^*_i$,  $(q^*_{i+1})^{(\bar u^i, \bar v^i, \bar w^i, \bar s^i)} \Vdash p_* \in \mathbf{G}$ (guaranteed by Observation $\ref{obskite}$).	Clause \ref{Mth} clearly holds, so we are done.
		
\end{PROOF}

We can turn to the proof of \ref{ee2}:
\begin{claim}\label{ee2p}
	For the forcing $\bbQ$ defined above clause \ref{ee2} holds.
\end{claim}
\begin{PROOF}{Claim \ref{ee2p}}
	Fix a  $\bbQ$-name $\name z$ with $q \Vdash_{\bbQ} \name z \in 2^\omega$.
	By Claim $\ref{fuzio}$ (and a standard density argument) we can assume, that 
	\begin{enumerate}[label = $(\XBox_{\arabic*})$, ref = $(\XBox_{\arabic*})$]

		\item	$\name z$ is a $\bbQ' = \prod_{\iota \in \{\mathbf{0},\mathbf{1}, \infty, \bbS\}} \bbQ^{\iota}_{X_\iota}$-name for some $X_{\mathbf{0}} \in [\lambda_\mathbf{0}]^{\aleph_0}$ , $X_{\mathbf{1}} \in [\lambda_\mathbf{1} \setminus \lambda_\mathbf{0}]^{\aleph_0}$, $X_{\mathbf{\infty}} \in [\lambda_\mathbf{\infty} \setminus \lambda_\mathbf{1}]^{\aleph_0}$,  $X_\bbS \in [\lambda_\bbS \setminus \lambda_\mathbf{\infty}]^{\aleph_0}$,
	\end{enumerate}
	moreover, w.l.o.g.\
		\begin{enumerate}[label = $(\XBox_{\arabic*})$, ref = $(\XBox_{\arabic*})$]
					\setcounter{enumi}{1}
		\item	$\Vdash_{\bbQ'} \name z \notin V$.
		\end{enumerate} 
	
		\newcounter{enubo} \setcounter{enubo}{2}
	\begin{enumerate}[label = $(\XBox_{\arabic*})$, ref = $(\XBox_{\arabic*})$]
		\setcounter{enumi}{2}
		\item \label{ab} Fix  enumerations $X_\mathbf{0} = \{ \varp^\mathbf{0}_j: \ j \in \omega \}$, $X_\mathbf{1} = \{\varp^{\mathbf{1}}_j: \ j \in \omega \}$,  $X_\mathbf{\infty} = \{\varp^{\mathbf{\infty}}_j: \ j \in \omega \}$,  $X_\bbS = \{\varp^\bbS_j: \ j \in \omega \}$.
	 
	 \item If $\varp \in X_\mathbf{1} \cup X_\infty \cup X_\bbS$, $q \in \bbQ'$, $\bar s \in \ ^{\omega>} 2$ we let $q^{\{ \varp \},(\bar s )} \in \bbQ'$ be defined as
	 $$\begin{array}{rl} q^{\{ \varp \},(\bar s )} \upharpoonright X_\mathbf{0} \cup X_\mathbf{1} \cup X_\infty \cup X_\bbS \setminus \{\varp \} = &  (q \upharpoonright X_\mathbf{0} \cup X_\mathbf{1} \cup X_\infty \cup X_\bbS \setminus \{\varp \}),\\
	  q^{\{ \varp\}, (\bar s)}(\varp) = & (q(\varp)^{(\bar s)}). \end{array}$$
  \item If $\varp \in X_\mathbf{1} \cup X_\infty \cup X_\bbS$, $q,p \in \bbQ'$,  $n \in \omega$, then 
  $q \leq_{\{\varp\},n} p$, if $q \leq p$ and $q(\varp) \leq_n p(\varp)$.
			\stepcounter{enubo} 	\stepcounter{enubo}  	\stepcounter{enubo}  \stepcounter{enubo}  \stepcounter{enubo}
	\end{enumerate}

	We will again need the terminology introduced in Definition \ref{D()}.
	
		\begin{definition}
		If $\Vdash_{\bbQ'} \name z \in 2^\omega$, $n \in \omega$, $\iota \in \{\mathbf{1}, \infty\}$, then we let $q \in D^{\iota,\textbf{un}}_{n}(\name z)$ (where we mean \textbf{un} as an abbreviation for ``unique"), iff 
		\begin{enumerate}[label = $(\roman*)$, ref =  $(\roman*)$]
			\item $q \in \bbQ'$, 
			\item  there exist  $k \in \omega$, and $i_0 \neq i_1 \in \{0,1\}$ for which
			\begin{enumerate}[label = $(ii)_{\arabic*}$, ref =  $(ii)_{\arabic*})$]
				\item $q^{\{\varp^\iota_n\},(\langle 0\rangle)} \Vdash \name z_k = i_0$,
				\item $q^{\{\varp^\iota_n\},(\langle 1\rangle)} \Vdash \name z_k = i_1$, and
				\item whenever $q \geq_{\{\varp^\iota_n\}, 1} r$, and $r^{\{\varp^\iota_n\},\langle 0\rangle }$ or $r^{\{ \varp^\iota_n \},\langle 1 \rangle}$ decides $\name z_j$ for some $j \neq k$ then so does $r$. 
			\end{enumerate} 
		\end{enumerate}
	\end{definition}
	
		\begin{definition}
		If $\Vdash_{\bbQ'} \name z \in 2^\omega$, $n \in \omega$, $\iota \in \{\mathbf{1}, \infty\}$, then we let $q \in D^{\iota,\textbf{eq}}_{n}(\name z)$, iff 
		whenever $q \geq_{\{\varp^\iota_n\}, 1} r$, and $r^{\{\varp^\iota_n\},(\langle 0 \rangle)}$ or $r^{\{\varp^\iota_n \}, (\langle 1 \rangle)}$ decides $\name z_k$ for some $k \in \omega$ then so does $r$.
	\end{definition}

	\begin{definition} \label{defmul}
	If $\Vdash_{\bbQ'} \name z \in 2^\omega$, $n \in \omega$, $\iota \in \{\mathbf{1}, \infty\}$, then we let $q \in D^{\iota,\textbf{mul}}_{n}(\name z)$, iff there is no $q' \leq_{\{\varp^\iota_n\},1} q$ with $q' \in D^{\iota,\textbf{un}}_{n}(\name z) \cup D^{\iota,\textbf{eq}}_{n}(\name z)$.
	\end{definition}
	
		Note the following:
	\begin{fact} \label{elfa}
		If $\Vdash_{\bbQ'} \name z \in 2^\omega$, $n \in \omega$, $\iota \in \{\mathbf{1}, \infty\}$, then
		\begin{enumerate}
			\item  \label{elfa1} $D^{\iota,\textbf{un}}_{n}(\name z) \cup D^{\iota,\textbf{mul}}_{n}(\name z) \cup D^{\iota,\textbf{eq}}_{n}(\name z)$ is dense (in fact, even $\leq_{\{\varp^\iota_n\}, 1}$-dense) in $\bbQ'$.
			\item \label{elfaop} if $q \in D^{\iota,\textbf{un}}_{n}(\name z))$
			 ($D^{\iota,\textbf{mul}}_{n}(\name z)$, $D^{\iota,\textbf{eq}}_{n}(\name z)$, resp.), and $q' \leq_{\{\varp^\iota_n\},1} q$, then $q' \in D^{\iota,\textbf{un}}_{n}(\name z))$
			 ($D^{\iota,\textbf{mul}}_{n}(\name z)$, $D^{\iota,\textbf{eq}}_{n}(\name z)$, resp.).
			 \item \label{triv} $D^{\iota,\textbf{un}}_{n}(\name z)$, $D^{\iota,\textbf{mul}}_{n}(\name z)$, $D^{\iota,\textbf{eq}}_{n}(\name z)$ are pairwise disjoint.
		\end{enumerate}
	\end{fact}

	

	 The proof of the present claim is by clarifying Subclaims $\ref{guszkl}$ and $\ref{megegy}$:
	 

	\begin{claim} \label{guszkl}
			Let $q \in \bbQ'$, $\Vdash_{\bbQ'} \name z \in 2^\omega$, and $\bar{\bar{\ell}} = ( \bar \ell^\mathbf{0}, \bar \ell^\mathbf{1}, \bar \ell^\mathbf{\infty}, \bar \ell^\bbS) \in \mathcal{L}$ be given,			and let $m \in \omega$ be fixed.
		Then for some  $r \in \bbQ'$,  $r\leq_{\bar{\bar{\ell}}} q$, for each $(\bar u, \bar v, \bar w, \bar s) \in \mathbf{seq}_{\bar{\bar{\ell}}}(r)$ one of the following holds:
		\begin{enumerate}[label = $\odot^m_{\arabic*}(r^{(\bar u, \bar v, \bar w, \bar s)})$:, ref =  $\odot^m_{\arabic*}(r^{(\bar u, \bar v, \bar w, \bar s)})$]
			\item $r^{(\bar u, \bar v, \bar w, \bar s)}$ forces that $\name z$ does not depend on $\{ p(\varp^\bbS_m): \ p \in \mathbf{G} \}$, i.e. there is no $p \in \bbQ'$, $p \leq r^{(\bar u, \bar v, \bar w, \bar s)}$ for which there exists  $k \in \omega$ and $c \in \{0,1\}$:
			\[ \begin{array}{rcl} p^{\{\varp^\bbS_m \}, (\langle 0 \rangle)} & \Vdash_{\bbQ'} & \name z_k = c, \\
			 p^{\{\varp^\bbS_m \}, (\langle 1 \rangle)} & \Vdash_{\bbQ'} & \name z_k = 1-c.\end{array} \]
			\item for each $p \in \bbQ'$, $p \leq r^{(\bar u, \bar v, \bar w, \bar s)}$  there exist $q \leq p$, $k \in \omega$, and $c \in \{0,1\}$ such that:
			\[ \begin{array}{rcl} q^{\{\varp^\bbS_m \}, (\langle 0 \rangle)}& \Vdash_{\bbQ'} & \name z_k = c, \\ 
			 q^{\{\varp^\bbS_m \}, (\langle 1 \rangle)} & \Vdash_{\bbQ'} & \name z_k = 1-c. \end{array}  \]
		\end{enumerate}
	\end{claim}
	\begin{PROOF}{Subclaim \ref{guszkl}}
		 	\newcounter{sclcc} \setcounter{sclcc}{0}
		Observe that
		\begin{enumerate}[label = $(\blacksquare_{\arabic*})$, ref = $(\blacksquare_{\arabic*})$]
			\item \label{megf} if $p \geq r \in \bbQ'$,  then $\odot_1^m(p) \to \odot^m_1(r)$, and similarly, $\odot^m_2(p) \to \odot^m_2(r)$ for every $m \in \omega$.
			\stepcounter{sclcc}
		\end{enumerate}
		
	Note that
		\begin{enumerate}[label = $(\blacksquare_{\arabic*})$, ref = $(\blacksquare_{\arabic*})$]
		\setcounter{enumi}{\value{sclcc}}
		\item if for $p \in \bbQ$ there is no extension $p' \leq p$ with $\odot^m_1(p')$, then
			$\odot^m_2(p)$ holds (and conversely), 
		\stepcounter{sclcc}
	\end{enumerate}
	therefore,
	\begin{enumerate}[label = $(\blacksquare_{\arabic*})$, ref = $(\blacksquare_{\arabic*})$]
		\setcounter{enumi}{\value{sclcc}}
		\item \label{dpd} for $n \in \omega$ 
		the set 
		$$D^m_\odot = \{ p \in \bbQ': \ \odot^m_1(p) \ \vee \ \odot^m_2(p)\}$$
		is dense open (and the sets $\{ p \in \bbQ': \ \odot^m_1(p) \}$, $\{ p \in \bbQ': \ \odot^m_2(p) \}$ are open).
		\stepcounter{sclcc}
	\end{enumerate}
	
	For later reference we remark the following corollary of Observation $\ref{obski}$:
	\begin{observation} \label{kellob}
		If $r \in \bbQ'$ is given by Subclaim $\ref{guszkl}$ (for a fixed $m$ and $\bar{\bar{\ell}} \in \mathcal{L}$), and $r \geq_{\bar{\bar{\ell}}} r'$, then for each $(\bar u, \bar v, \bar w, \bar s) \in \mathbf{seq}_{\bar{\bar{\ell}}}(r)$ we have
		$(r')^{(\bar u, \bar v, \bar w, \bar s)} \in D^m_\odot$, i.e. either $\odot_1^m((r')^{(\bar u, \bar v,  \bar w, \bar s)})$, or $\odot^m_2((r')^{(\bar u,  \bar v, \bar w, \bar s)})$ holds.
	\end{observation}

	\begin{fact}\label{prf}
		For every $p \in \bbQ'$, $m \in \omega$, $(\bar u, \bar v, \bar w, \bar s) \in \mathbf{seq}_{\bar{\bar{\ell}}}(p)$ there exists $p' \leq_{\bar{\bar{\ell}}} p$ for which either $\odot^m_1((p')^{(\bar u, \bar v, \bar w, \bar s)})$, or $\odot^m_2((p')^{(\bar u, \bar v, \bar w, \bar s)})$ holds.
	\end{fact}
	\begin{PROOF}{Fact \ref{prf}}
		Using \ref{dpd} choose $p'' \leq p^{(\bar u, \bar v, \bar w, \bar s)}$ with $p'' \in D^m_\odot$. By an argument similar to that of Observation \ref{ossze} we can assume that $p''(\varp^\infty_j) \leq_{\ell^\infty_j} p^{(\bar u, \bar v, \bar w, \bar s)}(\varp^\infty_j)$ for each $j$, so by Observation $\ref{obskite}$ there exists a condition $p' \leq_{\bar{\bar{\ell}}} p$ such that $(p')^{(\bar u, \bar v, \bar w, \bar s)}) = p''$.
	\end{PROOF}	
	Since $\leq_{\bar{\bar{\ell}}}$ is a partial order, enumerating the finite set
	 $$\left\{(\bar u, \bar v, \bar w, \bar s) \in \ \mathbf{seq}_{\bar{\bar{\ell}}}(q) \right\}$$
	  as $\{(\bar u^i, \bar v^i, \bar w^i, \bar s^i): \ i <M\}$ we can choose a sequence 
	  $$q_0 = q^*_0 \geq_{\bar{\bar{\ell}}}  \dots \geq_{\bar{\bar{\ell}}}  q^*_{M-1} \geq_{\bar{\bar{\ell}}}  q^*_M$$
	  requiring $(q^*_{i+1})^{(\bar u^i, \bar v^i, \bar w^i, \bar s^i)} \in D^m_\odot$ ($i < M$) (recall Observation $\ref{sv}$). Thus $r = q^*_{M}$ works.
	\end{PROOF}

	Now we can turn back to the proof of Claim $\ref{ee2p}$.
	
	\begin{definition} \label{ll}
		Fix a sequence $\langle \xi_n: \ n \in \omega \rangle$ that lists 
		$$X_{\mathbf{0}} \cup X_{\mathbf{1}} \cup X_{\mathbf{\infty}} \cup X_{\bbS} = \{ \varp^\iota_j: \ \iota \in \{\mathbf{0}, \mathbf{1}, \infty, \bbS\}, \ j \in \omega\}$$
			with each such element occurring infinitely many times (where the $X^\iota$'s are from \ref{ab}). Then we define
		\begin{enumerate}[label =$\arabic*)$, ref = $\arabic*)$]
			\item  	 the sequence $\langle \bar{\bar{\ell}}^n: \ n \in \omega \rangle$ so that
			\begin{itemize}
				\item $\bar{\bar{\ell}}^n = ((\bar \ell^n)^\mathbf{0}, (\bar \ell^n)^\mathbf{1}, (\bar \ell^n)^\mathbf{\infty}, (\bar \ell^n)^\bbS) \in \mathcal{L}$ for each $n$,
				\item $\bar{\bar{\ell}}^0$ consists of constant zero sequences,
				\item  if $\xi_n = \varp^\iota_m$, then we define $\bar{\bar{\ell}}^{n+1}$ so that 
				$$(\ell^{n+1})^{\iota'}_k = \left\{ \begin{array}{ll}(\bar \ell^{n})^{\iota'}_{k} +1, & \text{ if } \iota' = \iota \ \wedge \ k = m \\  (\bar \ell^n)^{\iota'}_k, & \text{ otherwise.} \end{array} \right.$$
			\end{itemize}
			\item \label{llt} for $q \in \bbQ'$, $n \in \omega$ and $(\bar u, \bar v, \bar w, \bar s) \in \mathbf{seq}_{\bar{\bar{\ell}}^n}(q)$ we define the sequence 
			$$\bar{t}^{(\bar u, \bar v, \bar w, \bar s)} = \langle t^{(\bar u, \bar v, \bar w, \bar s)}_j: \ j <n \rangle$$
			inductively as follows:
			if $k <n$, $K = |\{ j<k: \ \xi_k = \xi_j\}|$ and if $\xi_k =$
			\begin{itemize} 
				\item  $=\varp^{\mathbf{0}}_m$, then set $t^{(\bar u, \bar v, \bar w, \bar s)}_k = u_m(K)$,
				\item $=\varp^{\mathbf{1}}_m$, then set $t^{(\bar u, \bar v, \bar w, \bar s)}_k = v_m(K)$,
				\item $=\varp^{\mathbf{\infty}}_m$, then set $t^{(\bar u, \bar v, \bar w, \bar s)}_k = w_m(K)$,
				\item $=\varp^{\bbS}_m$, then set $t^{(\bar u, \bar v, \bar w, \bar s)}_k = s_m(K)$,
			\end{itemize}
		\item \label{lltt} for $q \in \bbQ'$ and the finite sequence $\bar t'$ 
		we let 
		$$\mathbf{qp}(q,\bar t') = (\bar u_*, \bar v_*, \bar w_*, \bar s_*) , $$
		if
		$$\bar t' = \bar t^{(\bar u_*, \bar v_*, \bar w_*, \bar s_*)},$$
		(where $(\bar u_*, \bar v_*, \bar w_*, \bar s_*) \in \mathbf{seq}_{\bar{\bar{\ell}}^{|\bar t|}}(q)$, and $\bar t^{(\bar u_*, \bar v_*, \bar w_*, \bar s_*)}$ is defined as above),

		\item and (for $q \in \bbQ'$)$, (\bar u, \bar v, \bar w, \bar s)$, $(\bar u', \bar v', \bar w', \bar s') \in \bigcup_{j<\omega} \mathbf{seq}_{\bar{\bar{\ell}}^j}(q)$ we
		define $(\bar u, \bar v, \bar w, \bar s) \sqsubseteq(\bar u', \bar v', \bar w', \bar s')$ naturally, i.e.
		\[ (\bar u, \bar v, \bar w, \bar s) \sqsubseteq(\bar u', \bar v', \bar w', \bar s'), \text{ iff } \bar t^{(\bar u, \bar v, \bar w, \bar s)} \subseteq \bar t^{(\bar u', \bar v', \bar w', \bar s')}, \]
		as well as 
		$$(\bar u, \bar v, \bar w, \bar s) \sqsubset (\bar u', \bar v', \bar w', \bar s'), $$ $$ \text{ iff } ((\bar u, \bar v, \bar w, \bar s) \sqsubseteq (\bar u', \bar v', \bar w', \bar s') \ \wedge \ (\bar u, \bar v, \bar w, \bar s) \neq (\bar u', \bar v', \bar w', \bar s'),$$
		\item for ($q \in \bbQ'$)$, (\bar u, \bar v, \bar w, \bar s) \in  \mathbf{seq}_{\bar{\bar{\ell}}^n}(q)$, $k \leq n$ we let  $ (\bar u, \bar v, \bar w, \bar s) \upharpoonright k$ to be the (unique)  member $(\bar u', \bar v', \bar w', \bar s')$ of $ \mathbf{seq}_{\bar{\bar{\ell}}^k}(q)$ for that 
		\[(\bar u', \bar v', \bar w', \bar s') \sqsubseteq (\bar u, \bar v, \bar w, \bar s). \]
		\end{enumerate}
	
	\end{definition}

	\begin{observation} \label{obst} { \ }
		\begin{enumerate}[label = $\alph*)$, ref = $\alph*)$] 
			\item \label{obst1}If $\bar t^{*}$ and $\bar t^{(\bar u, \bar v, \bar w, \bar s)}$ (with $(\bar u, \bar v, \bar w, \bar s) \in \mathbf{seq}_{\bar{\bar{\ell}}^{|\bar t^*|}}(q)$) satisfies that whenever $k$ is such that $t_k \neq  t^{(\bar u, \bar v, \bar w, \bar s)}_k$, then 
				\begin{enumerate}[label = $(\ast)_k$, ref = $(\ast)_k$]
					\item \label{csk} $t_k \in \{0,1\}$ and $\xi_k$ is \emph{not} of the form $\varp^\mathbf{0}_m$ (for any $m$),
				\end{enumerate}
			 \then \ $\mathbf{qp}(q,\bar t^*)$ is defined.
			\item \label{obst2} If 
			$$(\bar u, \bar v, \bar w, \bar s) \neq (\bar u', \bar v', \bar w', \bar s') \in \mathbf{seq}_{\bar{\bar{\ell}}^{n}}(q)$$ are such that $\bar t^{(\bar u, \bar v, \bar w, \bar s)}$ and $\bar t^{(\bar u', \bar v', \bar w', \bar s')}$ differs on exactly one coordinate, the $k$'th for which \ref{csk} holds, then there exists a condition $p' \leq_{\{\xi_k\}, 1} p^{(\bar u, \bar v, \bar w, \bar s) \rest k}$ such that
			$$ \{ (p')^{\{ \xi_k\}, (\langle 0 \rangle)}, (p')^{\{ \xi_k\}, (\langle 1 \rangle)} \} = \{ p^{(\bar u, \bar v, \bar w, \bar s)}, p^{(\bar u', \bar v', \bar w', \bar s')} \}.$$
		\end{enumerate}
	\end{observation}

	\begin{sclaim} \label{megegy}
			 Let $r$, $\name z$ be as in Subclaim $\ref{guszkl}$, and $\langle \bar{\bar{\ell}}^n: \ n \in \omega \rangle$ defined in Definition \ref{ll}.
			 Suppose that $r \Vdash \name z \notin V$.
	Then there exists a condition $r_* \in \bbQ'$, $r \geq r_*$, and 
	 $$ \bar y = \left\langle y^{(\bar u, \bar v, \bar w, \bar s)}: \ n \in \omega, \ (\bar u, \bar v, \bar w, \bar s) \in \mathbf{seq}_{\bar{\bar{\ell}}^n}(r_*) \right\rangle,$$
	 $$ \bar{\mathbf{x}} = \left\langle \mathbf{x}_{(\bar u, \bar v, \bar w, \bar s)}: \ n \in \omega, \ (\bar u, \bar v, \bar w, \bar s) \in \mathbf{seq}_{\bar{\bar{\ell}}^n}(r_*) \right\rangle, $$
	 such that for each $n$
	 \begin{enumerate}[label = $\varphi_{\alph*}(r_*$, ref = $\varphi_{\alph*}(r_*$]
	 	\item[$\varphi_{a}(r_*, \bar y, \bar{\mathbf{x}})$:] \label{r*a} for each $n$ and $(\bar u, \bar v, \bar w, \bar s) \in \mathbf{seq}_{\bar{\bar{\ell}}^n}(r_*)$: 
	 	\begin{itemize}
	 		\item $y^{(\bar u, \bar v, \bar w, \bar s)} \in \ ^{\omega>}2$,
	 		\item $\mathbf{x}_{(\bar u, \bar v, \bar w, \bar s)} \in \{\mathbf{un}, \mathbf{eq}, \mathbf{mul}\}$,
	 	\end{itemize} 
	 	\item[$\varphi_{b}(r_*,\bar y)$:] \label{r*b} for each $m <n$:
	 	\begin{itemize}
	 		\item  $(\bar u, \bar v, \bar w, \bar s) \in \mathbf{seq}_{\bar{\bar{\ell}}^m}(r_*)$, 
	 		\item  $(\bar u', \bar v', \bar w', \bar s') \in \mathbf{seq}_{\bar{\bar{\ell}}^n}(r_*)$,
	 	\end{itemize}
 	 		we have 
 	 		$$(\bar u, \bar v, \bar w, \bar s) \sqsubset (\bar u', \bar v', \bar w', \bar s') 
	 	\Rightarrow \ y^{(\bar u, \bar v, \bar w, \bar s)} \subsetneq y^{(\bar u', \bar v', \bar w', \bar s')}, $$
	 	\item[$\varphi_{c}(r_*, \bar y)$:] for each $n$, $(\bar u, \bar v, \bar w, \bar s) \in \mathbf{seq}_{\bar{\bar{\ell}}^n}(r_*)$:
	 	\[ (r_{*})^{(\bar u, \bar v, \bar w, \bar s)}\Vdash \name z \in [y^{(\bar u, \bar v, \bar w, \bar  s)}],\]
	 	\item[$\varphi_{d}(r_*,\bar y,  \bar{\mathbf{x}})$:] \label{r*d-} if $(\bar u, \bar v, \bar w, \bar s)  \in \mathbf{seq}_{\bar{\bar{\ell}}^{n+1}}(r_*)$, then 
	 	\begin{enumerate}[label = $\arabic*)$, ref = $\arabic*)$]
	 			\item if $\xi_{n} = \varp^\bbS_m$ for some $m$, then
				 $\odot^m_{1}(r_*^{(\bar u, \bar v, \bar w, \bar s) \upharpoonright n})$  $\vee$ $ \odot^m_{2}(r_*^{(\bar u, \bar v, \bar w, \bar s) \upharpoonright n})$ (from Subclaim $\ref{guszkl}$), and
	 			$$\mathbf{x}_{(\bar u, \bar v, \bar w, \bar s) \upharpoonright n} = \mathbf{eq} \ \iff \ \odot^m_{1}(r_*^{(\bar u, \bar v, \bar w, \bar s) \upharpoonright n}), $$
	 			$$\mathbf{x}_{(\bar u, \bar v, \bar w, \bar s) \upharpoonright n} = \mathbf{mul} \ \iff \ \odot^m_{2}(r_*^{(\bar u, \bar v, \bar w, \bar s) \upharpoonright n}), $$
	 			\item if $\xi_{n} = \varp^\iota_m$, where $\iota \in \{\mathbf{1},\infty\}$ for some $m$, then for each $(\bar u, \bar v, \bar w, \bar s) \in \mathbf{seq}_{\bar{\bar{\ell}}^n}(r_*)$
	 			\begin{itemize}
	 				\item $\mathbf{x}_{(\bar u, \bar v, \bar w, \bar s)} = \mathbf{eq}$, iff $r_*^{(\bar u, \bar v, \bar w, \bar s)} \in D^{\iota, \mathbf{eq}}_m(\name z)$,
	 			\end{itemize}
	 	\end{enumerate}

	 	\item[$\varphi_{e}(r_*,\bar y,  \bar{\mathbf{x}})$:] \label{r*d}  if $(\bar u, \bar v, \bar w, \bar s) \neq (\bar u', \bar v', \bar w', \bar s') \in \mathbf{seq}_{\bar{\bar{\ell}}^{n+1}}(r_*)$  (for some $n$), are such that $\bar u = \bar u'$, then the following implications hold true:
	 	\begin{enumerate}[label = $e\arabic*)$, ref = $e\arabic*)$]
	 		\item \label{E1} if $\xi_{n} = \varp^\bbS_m$ for some $m$,  $\bar s_m \neq \bar s'_m$, and  $\odot^m_{2}((r_*)^{(\bar u, \bar v, \bar w, \bar s)})$, then there exists 
	 		$$j \geq \lh(y^{(\bar u, \bar v, \bar w, \bar s) \upharpoonright n}, \lh(y^{(\bar u', \bar v', \bar w', \bar s') \upharpoonright n}),$$
 			 such that ($j < \lh(y^{(\bar u, \bar v, \bar w, \bar s)}, \lh(y^{(\bar u', \bar v', \bar w', \bar s')})$) and
 			 $$ y^{(\bar u, \bar v, \bar w, \bar s)}_j \neq y^{(\bar u', \bar v', \bar w', \bar s')}_j.$$

		\item \label{E2} if $\psi((\bar u, \bar v, \bar w, \bar s),(\bar u', \bar v', \bar w', \bar s'))$, under which we mean that ($\bar u = \bar u'$, and) for each $k < n+1$ either
			\begin{itemize}
				\item $t^{(\bar u, \bar v, \bar w, \bar s)}_k = t^{(\bar u', \bar v', \bar w', \bar s')}_k$, or
				\item $\mathbf{x}_{(\bar u, \bar v, \bar w, \bar s) \upharpoonright k-1} = \mathbf{eq}$,
			\end{itemize}
			then $y^{(\bar u, \bar v, \bar w, \bar s)} =y^{(\bar u', \bar v', \bar w', \bar s')}$,
		\item \label{E3} if  $\xi_{n} = \varp^\mathbf{1}_m$ or $\varp^\mathbf{\infty}_m$ for some $m$, and 
		\begin{enumerate}
			\item $\psi((\bar u, \bar v, \bar w, \bar s) \upharpoonright n,(\bar u', \bar v', \bar w', \bar s') \upharpoonright n)$, but 
			\item $\mathbf{x}_{(\bar u, \bar v, \bar w, \bar s) \upharpoonright n} \neq \mathbf{eq}$ and  $t^{(\bar u, \bar v, \bar w, \bar s)}_{n} \neq t^{(\bar u', \bar v', \bar w', \bar s')}_{n}$,
			\item  and $\xi_{n}$ is of the form $\varp^\mathbf{1}_m$, or $\varp^\infty_m$ for some $m$,
		\end{enumerate}
		then
			\begin{itemize}
				\item if $\xi_n = \varp^\mathbf{1}_m$, then  
				$$\begin{array}{l} \exists i< \lh(y^{(\bar u, \bar v, \bar w, \bar s)}), \lh(y^{(\bar u', \bar v', \bar w', \bar s')}): \\  y^{(\bar u, \bar v, \bar w, \bar s)}_i \neq y^{(\bar u', \bar v', \bar w', \bar s')}_i, \end{array}$$
				moreover, if this $i$ is unique, then  
				$$y^{(\bar u, \bar v, \bar w, \bar s)} \rest i =  y^{(\bar u', \bar v', \bar w', \bar s')} \rest i \neq \mathfrak{s}_i$$
				(where $\bar{\mathfrak{s}}$ is from \ref{x1}, \ref{x2}),
				\item if $\xi_n = \varp^\mathbf{\infty}_m$, then 
					$$\begin{array}{l} \exists i < i' < \lh(y^{(\bar u, \bar v, \bar w, \bar s)}), \lh(y^{(\bar u', \bar v', \bar w', \bar s')}): \\ 
						 y^{(\bar u, \bar v, \bar w, \bar s)}_i \neq y^{(\bar u', \bar v', \bar w', \bar s')}_i \ \text{ and} \\  y^{(\bar u, \bar v, \bar w, \bar s)}_{i'} \neq y^{(\bar u', \bar v', \bar w', \bar s')}_{i'}, \end{array}$$
				
			\end{itemize}
		
		 \end{enumerate}
 	\end{enumerate}
 	\end{sclaim}
 	First we verify Claim $\ref{ee2p}$ provided the extension $r_*$ of $r$ and the $y^{(\bar u, \bar v, \bar w, \bar s)}$'s given by Subclaim $\ref{megegy}$, i.e.\ satisfying $\varphi_a(r_*,\bar y, \bar{\mathbf{x}})$--$\varphi_e(r_*,\bar y, \bar{\mathbf{x}})$. 
 	First define
 	
 	$$T_* = \{ y^{(\bar u, \bar v, \bar w, \bar s)}: \ (\bar u, \bar v, \bar w, \bar s) \in \bigcup_{n \in \omega} \mathbf{seq}_{\bar{\bar{\ell}}^n}(r^*)\},$$
 	and note that
 	$$ r_* \Vdash \name z \in [T_*].$$
 	As 
 	$$\bbQ'  \simeq \bbQ^{\mathbf{0}}_{X_{\mathbf{0}}} \times \bbQ^{\mathbf{1}}_{X_{\mathbf{1}}}  \times \bbQ^{\mathbf{\infty}}_{X_{\mathbf{\infty}}} \times \bbQ^{\bbS}_{X_{\bbS}},$$ we can
 	\begin{enumerate}
 		\item first add a $\bbQ^\mathbf{0}_{X_\mathbf{0}}$-generic filter $\mathbf{G}_{X_\mathbf{0}}$ to $V$ with $r_* \rest X_\mathbf{0} \in \mathbf{G}_{X_{\mathbf{0}}}$, and define $T_\mathbf{0} \in V[\mathbf{G}_{X_\mathbf{0}}]$, such that $[T_\mathbf{0}]$ is $\bbG_0(\bar{\mathfrak{s}})$-independent, and
 		\begin{equation}\label{stad0} V[\mathbf{G}_{X_\mathbf{0}}] \models \ ``r_* \rest (X_\mathbf{1} \cup X_{\infty} \cup X_{\bbS}) \Vdash_{\bbQ' \rest (X_{\mathbf{1}}  \cup X_{\mathbf{\infty}} \cup X_{\bbS})} \name z \in [T_\mathbf{0}]", \end{equation}
 	\end{enumerate}
 		or
 	\begin{enumerate}
 		\setcounter{enumi}{1}
 		\item  add a $\bbQ^\mathbf{0}_{X_\mathbf{0}} \times \bbQ^\mathbf{1}_{X_{\mathbf{1}}}$-generic filter $\mathbf{G}_{X_\mathbf{0} \cup X_{\mathbf{1}}}$ to $V$ with $r_* \rest (X_\mathbf{0} \cup X_{\mathbf{1}}) \in \mathbf{G}_{X_\mathbf{0} \cup X_{\mathbf{1}}}$, and define $T_\mathbf{1} \in V[\mathbf{G}_{X_\mathbf{0} \cup X_\mathbf{1}}]$, such that $[T_\mathbf{1}]$ is $\bbG_1$-independent, and
 		\begin{equation}\label{stad1} V[\mathbf{G}_{X_\mathbf{0} \cup X_\mathbf{1}}] \models \ ``r_* \rest (X_{\infty} \cup X_{\bbS}) \Vdash_{\bbQ' \rest (X_{\mathbf{\infty}} \cup X_{\bbS})} \name z \in [T_\mathbf{1}]", \end{equation}
 	\end{enumerate} 
 	or 
 	\begin{enumerate}
 		\setcounter{enumi}{2}
 		\item  add a $\bbQ^\mathbf{0}_{X_\mathbf{0}} \times \bbQ^\mathbf{1}_{X_{\mathbf{1}}} \times \bbQ^\infty_{X_{\mathbf{\infty}}}$-generic filter $\mathbf{G}_{X_\mathbf{0} \cup X_{\mathbf{1}} \cup X_{\mathbf{\infty}}}$ to $V$ with $r_* \rest (X_\mathbf{0} \cup X_{\mathbf{1}} \cup X_{\mathbf{\infty}}) \in \mathbf{G}_{X_\mathbf{0} \cup X_{\mathbf{1}} \cup X_{\mathbf{\infty}}}$, and define $T_\mathbf{\infty} \in V[\mathbf{G}_{X_\mathbf{0} \cup X_\mathbf{1} \cup X_{\mathbf{\infty}}}]$, such that $[T_\mathbf{\infty}]$ is $E_0$-independent, and
 		\begin{equation}\label{stad8} V[\mathbf{G}_{X_\mathbf{0} \cup X_\mathbf{1} \cup X_{\mathbf{\infty}}}] \models \ ``r_* \rest (X_{\bbS}) \Vdash_{\bbQ' \rest  X_{\bbS}} \name z \in [T_\mathbf{\infty}]". \end{equation}
 	\end{enumerate} 
 	So fix the mutually generic filters $\mathbf{G}_{X_\mathbf{0}}$, $\mathbf{G}_{X_\mathbf{1}}$, $\mathbf{G}_{X_{\mathbf{\infty}}}$ (containing $r_* \rest X_{\mathbf{0}}$, $r_* \rest X_{\mathbf{1}}$ and $r_* \rest X_{\mathbf{\infty}}$),  define $T_\mathbf{0} \in V[\mathbf{G}_{X_{\mathbf{0}}}]$, $T_\mathbf{1} \in V[\mathbf{G}_{X_\mathbf{0}} \times \mathbf{G}_{X_\mathbf{1}}]$, $T_\infty \in V[\mathbf{G}_{X_\mathbf{0}} \times \mathbf{G}_{X_\mathbf{1}} \times \mathbf{G}_{X_\mathbf{\infty}}]$ as follows:  
 	$$T_\mathbf{0} = \{ y^{(\bar u, \bar v, \bar w, \bar s)}: \ (\bar u, \bar v, \bar w, \bar s) \in \bigcup_{n \in \omega} \mathbf{seq}_{\bar{\bar{\ell}}^n}(r_*): \  r_*^{(\bar u, \bar v, \bar w, \bar s)} \rest X_{\mathbf{0}} \in \mathbf{G}_{X_\mathbf{0}} \},$$
 		$$T_\mathbf{1} = \{ y^{(\bar u, \bar v, \bar w, \bar s)}: \ (\bar u, \bar v, \bar w, \bar s) \in \bigcup_{n \in \omega} \mathbf{seq}_{\bar{\bar{\ell}}^n}(r_*): \  r_*^{(\bar u, \bar v, \bar w, \bar s)} \rest (X_{\mathbf{0}} \cup X_{\mathbf{1}}) \in \mathbf{G}_{X_\mathbf{0} \cup X_\mathbf{1}}\},$$
		$$T_\mathbf{\infty} = \{ y^{(\bar u, \bar v, \bar w, \bar s)}: \ (\bar u, \bar v, \bar w, \bar s) \in \bigcup_{n \in \omega} \mathbf{seq}_{\bar{\bar{\ell}}^n}(r_*): \  r_*^{(\bar u, \bar v, \bar w, \bar s)} \rest (X_{\mathbf{0}} \cup X_{\mathbf{1}} \cup X_{\infty}) \in \mathbf{G}_{X_\mathbf{0} \cup X_\mathbf{1} \cup X_{\infty}}\},$$

 	(Recalling Definition $\ref{def1}$), for each fixed $n$
 	$$\{ r_*^{(\bar u, \bar v, \bar w, \bar s)}: \ (\bar u, \bar v, \bar w, \bar s) \in  \mathbf{seq}_{\bar{\bar{\ell}}^n}(r_*) \} 	\text{ is predense below }r_*,$$
 hence a standard density argument implies $\eqref{stad0}$, $\eqref{stad1}$, and $\eqref{stad8}$. It remains to check that $[T_\mathbf{0}]$ ($[T_\mathbf{1}]$, $[T_\mathbf{\infty}]$, resp.) is indeed $\mathbb{G}_0(\bar{\mathfrak{s}})$- ($\mathbb{G}_1$-, $E_0$-, resp.)-independent.
 	For these one only needs to check the following assertions (using $\varphi_b(r^*,\bar y)$,  $\varphi_e(r^*, \bar y)$ from Subclaim $\ref{megegy}$), which is left to the reader:
 	\begin{itemize}
 		\item For every branch $\langle b_i: \ i \in \omega \rangle$ in $T_*$ there is an infinite sequence $\langle (\bar u^i, \bar v^i, \bar w^i, \bar s^i): \ i \in \omega \rangle$, such that $(\bar u^i, \bar v^i, \bar w^i, \bar s^i) \in \mathbf{seq}_{\bar{\bar{\ell}}^i}(r_*)$ with 
 		$$(\bar u^i, \bar v^i, \bar w^i, \bar s^i) \sqsubseteq (\bar u^{i+1}, \bar v^{i+1}, \bar w^{i+1}, \bar s^{i+1}),$$ and $b_i = y^{(\bar u^i, \bar v^i, \bar w^i, \bar s^i)}$ (for each $i \in \omega$) (recall that $\mathbf{seq}_{\bar{\bar{\ell}}^i}(r_*)$ is finite by \eqref{seqfi}, and use Kőnig's theorem).
 			\item If $\langle (\bar u^i, \bar v^i, \bar w^i, \bar s^i): \ i \in \omega \rangle$ and $\langle ( (\bar u')^i,  (\bar v')^i, (\bar  w')^i,  ( \bar s')^i): \ i \in \omega \rangle$ are different, $\sqsubseteq$-increasing,  $(\bar u^i, \bar v^i, \bar w^i, \bar s^i), ( ( \bar u')^i, (\bar v')^i, (\bar w')^i,  (s')^i) \in \mathbf{seq}_{\bar{\bar{\ell}}^i}(r_*)$,  and we have $\bar u^i = (\bar \bar u')^i$  for each $i$, then at least one of the following holds:
 			\begin{enumerate}[label = $(\wr)_\arabic*$, ref = $(\wr)_\arabic*$]
 				\item \label{wr1} for each $n$ the premise \ref{E2} (from $\varphi_{e}(r_*,\bar y)$) holds, and thus 
 				$$y^{(\bar u^n, \bar v^n, \bar w^n, \bar s^n)} = y^{(\bar (u')^n, \bar (v')^n, \bar (w')^n, \bar (s')^n)}.$$
 				\item \label{wr2} for some $n$ the premise in \ref{E1} holds, so it holds infinitely many often (since $\{q \in \bbQ': \ \odot^m_2(q)\}$ is open for arbitrary $m$) and so $\cup\{y^{(\bar u^n, \bar v^n, \bar w^n, \bar s^n)}: \ n \in \omega \}$ and $\cup\{ y^{((\bar u')^n, (\bar v')^n, (\bar w')^n, (\bar s')^n)}: \ n \in \omega\}$ differ on infinitely many digits.
 				\item \label{wr3}for some $n$ the premise in \ref{E3} holds, and so $\cup\{y^{(\bar u^n, \bar v^n, \bar w^n, \bar s^n)}: \ n \in \omega \}$ and $\cup\{ y^{((\bar u')^n, (\bar v')^n, (\bar w')^n, (\bar s')^n)}: \ n \in \omega\}$ are not connected in $\mathbb{G}_0(\bar{\mathfrak{s}})$.
 			\end{enumerate}
 			\item If $\langle (\bar u^i, \bar v^i, \bar w^i, \bar s^i): \ i \in \omega \rangle$ and $\langle ((\bar u')^i, (\bar v')^i, (\bar w')^i, (\bar s')^i): \ i \in \omega \rangle$ are different, $\sqsubseteq$-increasing,  $(\bar u^i, \bar v^i, \bar w^i, \bar s^i), ((\bar u')^i, (\bar v')^i, (\bar w')^i, (\bar s')^i) \in \mathbf{seq}_{\bar{\bar{\ell}}^i}(r_*)$,  and we have $\bar u^i = (\bar u')^i$,  $\bar v^i = (\bar v')^i$ for each $i$, then either \ref{wr1} or \ref{wr2} holds, or
 			\begin{enumerate}[label = $(\wr)'_\arabic*$, ref = $(\wr)'_\arabic*$]
 				\setcounter{enumi}{2}
 				\item for some $n$ the premise in \ref{E3} holds,where necessarily $\xi_n = \varp^\infty_m$ for some $m$ (since $\bar v^i = (\bar v')^i$ for each $i$, in particular $v^{n+1}_m = (v')^{n+1}_m$) and so $\cup\{y^{(\bar u^n, \bar v^n, \bar w^n, \bar s^n)}: \ n \in \omega \}$ and $\cup\{ y^{((\bar u')^n, (\bar v')^n, (\bar w')^n, (\bar s')^n)}: \ n \in \omega\}$ differ in at least two digits.
 			\end{enumerate}
 			\item If $\langle (\bar u^i, \bar v^i, \bar w^i, \bar s^i): \ i \in \omega \rangle$ and $\langle ((\bar u')^i, (\bar v')^i, (\bar w')^i, (\bar s')^i): \ i \in \omega \rangle$ are different, $\sqsubseteq$-increasing,  $(\bar u^i, \bar v^i, \bar w^i, \bar s^i), ((\bar u')^i, (\bar v')^i, (\bar w')^i, (\bar s')^i) \in \mathbf{seq}_{\bar{\bar{\ell}}^i}(r_*)$,  and we have $\bar u^i = (\bar u')^i$,  $\bar v^i = (\bar v')^i$, $\bar w^i = (\bar w')^i$  for each $i$, then either \ref{wr1} or \ref{wr2} holds. 		
 	\end{itemize}
 	Note that the assertions above are absolute between transitive models.
 	
	 \begin{PROOF}{Subclaim \ref{megegy}}(Subclaim \ref{megegy})
	Similarly to that in the proof of Claim $\ref{fuzio}$ 
	 \newcounter{bcou} \setcounter{bcou}{0}
	\begin{enumerate}[label = $(\blacklozenge_{\arabic*})$, ref = $(\blacklozenge_{\arabic*})$]
		\setcounter{enumi}{\value{bcou}}
		\item \label{nemkez}	we are going to define the sequences $$ \begin{array}{l}
		\langle r_i: \ i < \omega \rangle, \\
	\langle y^{(\bar u, \bar v, \bar w, \bar s)}: \ (\bar u, \bar v, \bar w, \bar s) \in \mathbf{seq}_{\bar{\bar{\ell}}^i}(r_i), i \in \omega \rangle \\
		\langle \mathbf{x}_{(\bar u, \bar v, \bar w, \bar s)}: \ (\bar u, \bar v, \bar w, \bar s) \in \mathbf{seq}_{\bar{\bar{\ell}}^i}(r_i), i \in \omega \rangle
 \end{array}$$ 
		satisfying the requirements of the following scheme:
	\newcounter{enuboo} \setcounter{enuboo}{0}
	\begin{enumerate}[label = $(\XBox_{\arabic*})$, ref = $(\XBox_{\arabic*})$]
		\item \label{r0} $r_0 = r$, and for each $i$ $r_i \in \bbQ'$,
			\item \label{qla} $\forall i$ $r_{i} \geq_{\bar{\bar{\ell}}^i} r_{i+1}$,
			\item \label{r2} for each $i$: 
			\begin{itemize}
				\item $\varphi_a(r_{i}, \bar y^{*(i)}, \bar{\mathbf{x}}^{*(i-1)})$,
				\item $\varphi_b(r_{i}, \bar y^{*(i)})$,
				\item $\varphi_c(r_{i}, \bar y^{*(i)})$,
				\item $\varphi_d(r_{i},\bar y^{*(i)}, \bar{\mathbf{x}}^{*(i-1)})$,
				\item $\varphi_e(r_{i},\bar y^{*(i)}, \bar{\mathbf{x}}^{*(i-1)})$,
			\end{itemize}
			where $\bar y^{*(i)}$ is a restriction of the sequence $\bar y$ defined as
			$$\bar y^{*(i)} = \langle y^{(\bar u, \bar v, \bar w, \bar s)}: \ (\bar u, \bar v, \bar w, \bar s) \in \mathbf{seq}_{\bar{\bar{\ell}}^k}(r_{i}), k \leq i \rangle,$$
			and $$\bar{\mathbf{x}}^{*(i-1)} = \langle \mathbf{x}_{(\bar u, \bar v, \bar w, \bar s)}: \ (\bar u, \bar v, \bar w, \bar s) \in \mathbf{seq}_{\bar{\bar{\ell}}^k}(r_{i}), \ k \in i \rangle.$$
				
		\stepcounter{enuboo} 	\stepcounter{enuboo}  	
	\end{enumerate} 
	\stepcounter{bcou}
	\end{enumerate}
	Again, once we have constructed the $r_i$'s and $y^{(\bar u,\bar v, \bar w, \bar s)}$'s, we can let $r_*$ be
	a common lower bound of the $r_i$'s such that for each $i$ $r_i \geq_{\bar{\bar{\ell}}^i} r_*$ holds. Then  for each $i<j$ $r_i \geq_{\bar{\bar{\ell}}^i} r_{j} \geq_{\bar{\bar{\ell}}^i} r_*$, so by Observation $\ref{sv}$
	for each $j>i$:
	 $$\mathbf{seq}_{\bar{\bar{\ell}}^i}(r_i)= \mathbf{seq}_{\bar{\bar{\ell}}^i}(r_{j}) = \mathbf{seq}_{\bar{\bar{\ell}}^i}(r_*),$$
	 and for each $(\bar u, \bar v, \bar w, \bar s) \in \ \mathbf{seq}_{\bar{\bar{\ell}}^i}(r_i)$ we have 
	 $$r_{i}^{(\bar u, \bar v, \bar w, \bar s)} \geq r_j^{(\bar u, \bar v, \bar w, \bar s)} \geq r_*^{(\bar u, \bar v, \bar w, \bar s)}.$$ 
	 Also note that by Observation $\ref{kellob}$ (and recalling \ref{dpd})
	 if $\xi_i = \varp^\bbS_m$ for some $m$, and $(\bar u, \bar v, \bar w, \bar s) \in \mathbf{seq}_{\bar{\bar{\ell}}^i}(r_{i+1})$, then
	  $$\begin{array}{l}
	 	\odot^m_1(r_{i+1}^{(\bar u, \bar v, \bar w, \bar s)}) \iff \odot^m_1(r_{*}^{(\bar u, \bar v, \bar w, \bar s)}), \\ 
	 	 \odot^m_2(r_{i+1}^{(\bar u, \bar v, \bar w, \bar s)}) \iff \odot^m_2(r_{*}^{(\bar u, \bar v, \bar w, \bar s)}). \end{array} $$
 	 Similarly, 
 	 for each $i$, if $\xi_i = \varp^\iota_m$ for some $m$, where $\iota \in \{\mathbf{1}, \infty\}$, and $(\bar u, \bar v, \bar w, \bar s) \in \mathbf{seq}_{\bar{\bar{\ell}}^i}(r_{i+1})$, then $r_{i+1}^{(\bar u, \bar v, \bar w, \bar s)} \in D_m^{\iota,\mathbf{x}_{(\bar u, \bar v, \bar w, \bar s)}}(\name z)$ (by \eqref{r2}, recalling the definition of $\varphi_d(.,.)$ in Subclaim \ref{megegy}). So
 	 we obtain (by Observation \ref{elfaop+}), that
 	 $$ \forall r \leq_{\bar{\bar{\ell}}^{i+1}} r_{i+1}: \ \ r \in  D_m^{\iota,\mathbf{x}_{(\bar u, \bar v, \bar w, \bar s)}}(\name z).$$
		This clearly implies that if for each $i$ and $g \in  \{a,b,c,d,e \}$  we have $\varphi_g(r_{i}, \bar y^{*(i)}, \bar{\mathbf{x}}^{*(i-1)} )$ (or $\varphi_g(r_{i}, \bar y^{*(i)})$) holds,  then for each $g \in \{a,b,c,d,e\}$ $\varphi_g(r_*, \bar y, \bar{\mathbf{x}})$ (or $\varphi_g(r_*, \bar y)$)  holds, too, as desired.

	\begin{enumerate}[label = $(\blacklozenge_{\arabic*})$, ref = $(\blacklozenge_{\arabic*})$]
		\setcounter{enumi}{\value{bcou}}
		\item	So suppose that we have set $r_0 = r$, and we have already defined 
		$$r_0 \geq_{\bar{\bar{\ell}}^0} r_1 \geq_{\bar{\bar{\ell}}^1} r_2 \geq_{\bar{\bar{\ell}}^2} \dots \geq_{\bar{\bar{\ell}}^{n-1}} r_n$$
		satisfying \ref{r0}-\ref{r2}.
		\stepcounter{bcou}
	\end{enumerate}

	Depending on the value of $\xi_n$ we do the following.
	\begin{enumerate}[label = $(\blacklozenge_{3})(\roman*)$, ref = $(\blacklozenge_{3})(\roman*)$]
		\item Case $i$: $\xi_n = \varp^\mathbf{0}_m$ for some $m$:

	\end{enumerate}
		Let	$N = | \mathbf{seq}_{\bar{\bar{\ell}}^{n+1}}(r_n)|$, and fix an enumeration
	$ \langle (\bar u^j, \bar v^j, \bar w^j, \bar s^j): \ j < N \rangle.$
	Then defining the $\leq_{\bar{\bar{\ell}}^{n+1}}$-decreasing sequence $\langle p_j: \ j \leq M \rangle$ with $p_0= r_n$, and $p_{j+1}^{(\bar u^j, \bar v^j, \bar w^j, \bar s^j)}$ deciding  the value of $\name z_k$ (recalling the first part of Observation \ref{kit}), where $k = \lh(y^{(\bar u^j, \bar v^j, \bar w^j, \bar s^j) \rest n})$, setting $r_{n+1} = p_N$ works.
		
		\begin{enumerate}[label = $(\blacklozenge_{3})(\roman*)$, ref = $(\blacklozenge_{3})(\roman*)$]
		\setcounter{enumi}{1}
		\item Case $ii$: $\xi_n = \varp^\bbS_m$ (for some $m$).
	\end{enumerate}
	\begin{enumerate}[label = $(\blacklozenge_{3})(ii)_\arabic*$, ref = $(\blacklozenge_{3})(ii)_\arabic*$]
		\item \label{ii1}	Let	$M = | \mathbf{seq}_{\bar{\bar{\ell}}^{n}}(r_n)|$, and fix an enumeration
		$ \langle (\bar u^j, \bar v^j, \bar w^j, \bar s^j): \ j < M \rangle.$
		Then (again by Observation \ref{kit}) defining the $\leq_{\bar{\bar{\ell}}^{n}}$-decreasing sequence $\langle p_j: \ j \leq M \rangle$ with $p_0= r_n$, and 
		$$\odot_1^m(p_{j+1}^{(\bar u^j, \bar v^j, \bar w^j, \bar s^j)}) \vee \odot_2^m(p_{j+1}^{(\bar u^j, \bar v^j, \bar w^j, \bar s^j)})$$ and set
		$r^*_{n} = p_M$
		(so easily 
		\begin{equation}
			\mathbf{seq}_{\bar{\bar{\ell}}^{n}}(r_n) = \mathbf{seq}_{\bar{\bar{\ell}}^{n}}(r^*_n),
		\end{equation}
		and
		\begin{equation}
			\forall j < M: \ \odot_1^m((r^*_{n})^{(\bar u^j, \bar v^j, \bar w^j, \bar s^j)}) \vee \odot_2^m((r_n^*)^{(\bar u^j, \bar v^j, \bar w^j, \bar s^j)})).
		\end{equation}
		
	\end{enumerate}

\begin{enumerate}[label = $(\blacklozenge_{3})(ii)_\arabic*$, ref = $(\blacklozenge_{3})(ii)_\arabic*$]
		\setcounter{enumi}{1}
		\item \label{ppro} 	 Now let 
		$$\begin{array}{rl} Y = \{ \langle (\bar u, \bar v, \bar w, \bar s), (\bar u', \bar v', \bar w', \bar s') \rangle: & (\bar u, \bar v, \bar w, \bar s), (\bar u', \bar v', \bar w', \bar s')\in \mathbf{seq}_{\bar{\bar{\ell}}^{n+1}}(r^*_n), \\ & \bar u = \bar u'\}|, \end{array} $$
		$$N = |Y|,$$
		and fix the enumeration
		$$\left\langle \langle (\bar u^j, \bar v^j, \bar w^j, \bar s^j), ((\bar{u}')^j, (\bar v')^j, (\bar w')^j, (\bar s')^j) \rangle: \ j <N \right\rangle \text{ of } Y. $$
		We are going to construct
		\begin{itemize} 
			\item  the sequence 
		$$q_0 = r^*_{n} \geq_{\bar{\bar{\ell}}^{n+1}} q_1 \geq_{\bar{\bar{\ell}}^{n+1}} q_2 \geq_{\bar{\bar{\ell}}^{n+1}} \dots \geq_{\bar{\bar{\ell}}^{n+1}} q_{N},$$
		\item and for each $ (\bar u, \bar v, \bar w, \bar s) \in \ \mathbf{seq}_{\bar{\bar{\ell}}^{n+1}}(r_n^*)$
		$$x^{(\bar u, \bar v, \bar w, \bar s)}_0 \subseteq x^{(\bar u, \bar v, \bar w, \bar s)}_1 \subseteq \dots \subseteq x^{(\bar u, \bar v, \bar w, \bar s)}_{N}$$
		with 
		\begin{enumerate}[label =$(\alph*)$, ref = $(\alph*)$]
			\item \label{fa}  $x^{(\bar u, \bar v, \bar w, \bar s)}_k \in \ ^{\omega>}2$ ($k\leq N$),
			\item 	$x^{(\bar u, \bar v, \bar w, \bar  s)}_0 = y^{(\bar u, \bar v, \bar w, \bar s) \rest n},$
			\item \label{cforsz} and $ q_k^{(\bar u, \bar v, \bar w, \bar s)} \Vdash \name z \in [x^{(\bar u, \bar v, \bar w, \bar s)}_k]$,
			\item \label{d} and for each $k <N$ 
			 $$x^{(\bar u^k, \bar v^k, \bar w^k, \bar s^k)}_k \subsetneq x^{(\bar u^k, \bar v^k, \bar w^k, \bar s^k)}_{k+1},$$
			  $$x^{((\bar u')^k, (\bar v')^k, (\bar w')^k, (\bar s')^k)}_k \subsetneq x^{((\bar u')^k, (\bar v')^k, (\bar w')^k, (\bar s')^k)}_{k+1},$$
			\item \label{ftln} for each $k<N$  if $s^k_m \neq (s')^k_m$ and $\odot^m_{2}((p_k)^{(\bar u^k, \bar v^k, \bar w^k, \bar s^k)})$ (from Subclaim $\ref{guszkl}$) hold,
				then for some 
				$$j \geq \lh(x^{(\bar u^k, \bar v^k, \bar w^k, \bar s^k)}_{k}), \lh(x^{((\bar u')^k, (\bar v')^k, (\bar w')^k, (\bar s')^k)}_k)$$
				  $$x^{(\bar u^k, \bar v^k, \bar w^k, \bar s^k)}_{k+1}(j) \neq x^{((\bar u')^k, (\bar v')^k, (\bar w')^k, (\bar s')^k)}_{k+1}(j).$$ 
				
	
		\end{enumerate}
	
		\end{itemize}
		
		\stepcounter{bcou}
	\end{enumerate}
	Before constructing the $q_k$'s ($k\leq N$) and $x_k^{(\bar u, \bar v, \bar w, \bar s)}$'s ($(\bar u, \bar v, \bar w, \bar s) \in \mathbf{seq}_{\bar{\bar{\ell}}^{n+1}}(r_n^*)$) 
	\begin{enumerate}[label = $(\blacklozenge_{3})(ii)_\arabic*$, ref = $(\blacklozenge_{3})(ii)_\arabic*$]
		\setcounter{enumi}{2}
		\item \label{pppro} we clarify why setting 
		$$\begin{array}{l} r_{n+1} = q_N,\\ 
			y^{(\bar u, \bar v, \bar w, \bar s)} = \text{ the maximal element in } \{x \in \ ^{\omega>}2 : \ q_N^{(\bar u, \bar v, \bar w, \bar s)} \Vdash x \in \name z \} (\supseteq x^{(\bar u, \bar v, \bar w, \bar s)}_N) \end{array} $$
		works for our purposes (i.e. satisfying \ref{qla}, \ref{r2}):
	\end{enumerate}
	First recall that $r_{n+1} \leq r_0 \Vdash \name z \notin V$, so the maximal elements above do really exist, and so the  $y^{(\bar u, \bar v, \bar w, \bar s)}$'s are well defined finite sequences.
	
	First,  \ref{d} clearly implies $\varphi_b(r_{n+1}, \bar y^{*(n+1)})$.	Second, as $p_0 = r_n \geq_{\bar{\bar{\ell}}^n} p_M =  r_n^* \geq_{\bar{\bar{\ell}}^{n+1}} q_N = r_{n+1}$ recalling Observation \ref{sv} clearly $\mathbf{seq}_{\bar{\bar{\ell}}^{n}}(r_n^*) = \mathbf{seq}_{\bar{\bar{\ell}}^{n}} (r_{n+1})$  holds. 
Recalling \ref{ii1} for each $j <M$ either 	
$$\odot^m_{1}((p_{j+1})^{(\bar u^j, \bar v^j, \bar w^j, \bar s^j)}), \text{ or } \odot^m_{2}((p_{j+1})^{(\bar u^j, \bar v^j, \bar w^j, \bar s^j)})$$
(where $(\bar u^j, \bar v^j, \bar w^j, \bar s^j)$ is meant as the $j$'th entry on the list in \ref{ii1}).
 Now (by $p_k \geq_{\bar{\bar{\ell}}^n} p_M = r_n^* \geq_{\bar{\bar{\ell}}^n} r_{n+1}$) clearly either $$\odot^m_{1}((r_{n+1})^{(\bar u^j, \bar v^j, \bar w^j, \bar s^j)})), \text{ or } \odot^m_{2}((r_{n+1})^{(\bar u^j, \bar v^j, \bar w^j, \bar s^j)})$$ holds, and $\varphi_d(r_{n+1}, \bar y^{*(n+1)})$ follows. 
	Finally, for $\varphi_e(r_{n+1}, \bar{y}^{*(n+1})$ we need to check clause \ref{E2}, i.e.\ when there is no reason for $y^{(\bar u, \bar v, \bar w, \bar s)}$ and $y^{(\bar u', \bar v', \bar w', \bar s')}$ to be different, then these two are the same, which follows from the next claim (applying to $r_{n+1} = q_N$), so \ref{r2} holds, indeed:
	\begin{sclaim} \label{kisza}
		Assume that $p \in \bbQ'$, $(\bar u, \bar v, \bar w, \bar s),(\bar u', \bar v', \bar w', \bar s') \in \mathbf{seq}_{\bar{\bar{\ell}}^{n+1}}(p)$, such that $\bar u = \bar u'$, and whenever $k \leq n$
		\begin{itemize} 
			\item if $\xi_k = \varp^\iota_m$, for some $\iota  \in \{\mathbf{1}, \infty\}$ and $m \in \omega$, then $p^{(\bar u, \bar v, \bar w, \bar s) \rest k} \in D^{\iota,\mathbf{eq}}_m(\name z)$,
			\item if $\xi_k = \varp^\bbS_m$, for some $m \in \omega$, then $\odot^m_1(p^{(\bar u, \bar v, \bar w, \bar s) \rest k})$,
		\end{itemize}
		Then for every $x \in \ ^{\omega>}2$
		\[ p^{(\bar u, \bar v, \bar w, \bar s)} \Vdash \name z \in [x] \ \iff \ p^{(\bar u', \bar v', \bar w', \bar s')} \Vdash \name z \in [x]. \]
	\end{sclaim}
	\begin{PROOF}{Subclaim \ref{kisza}}
		Fix $x \in \ ^{\omega>}2$, and let 
		\begin{itemize}
			\item $(\bar u_*^0, \bar v_*^0, \bar w_*^0, \bar s_*^0) = (\bar u', \bar v', \bar w', \bar s') \in \mathbf{seq}_{\bar{\bar{\ell}}^{n+1}}(p),$
			\item $t^0 = t^{(\bar u_*^0, \bar v_*^0, \bar w_*^0, \bar s_*^0)}$,
			\item $(\bar u_*^{n+1}, \bar v_*^{n+1}, \bar w_*^{n+1}, \bar s_*^{n+1}) = (\bar u, \bar v, \bar w, \bar s) \in \mathbf{seq}_{\bar{\bar{\ell}}^{n+1}}(p),$
			\item $ t^{n+1} = t^{(\bar u_*^{n+1}, \bar v_*^{n+1}, \bar w_*^{n+1}, \bar s_*^{n+1})},$
		\end{itemize}
		as defined in \ref{llt} from Definition \ref{ll}.
		Now set 
		$$\bar t^k = \bar t^{(\bar u_*^{n+1}, \bar v_*^{n+1}, \bar w_*^{n+1}, \bar s_*^{n+1})} \rest [0, k) \ \cup \  (\bar t^{(\bar u_*^{0}, \bar v_*^{0}, \bar w_*^{0}, \bar s_*^{0})}\rest [k,n+1)  \ \ (k \leq n+1),$$
		and let 
		$$ (\bar u_*^{k}, \bar v_*^{k}, \bar w_*^{k}, \bar s_*^{k}) = \mathbf{qp}(p,\bar t^k)$$
		(which exists by clause \ref{obst1} from Observation \ref{obst}).
		Let $p_k = p^{(\bar u_*^{k}, \bar v_*^{k}, \bar w_*^{k}, \bar s_*^{k})}$ for $k \leq n+1$, observe that $p^{(\bar u, \bar v, \bar w, \bar s)} = p_{n+1}$, $p^{(\bar u', \bar v', \bar w', \bar s')} = p_{0}$. 	We claim that
		\begin{equation} \label{ezkell} 
			\text{for each } k \leq n: \ p_k \Vdash \name z \in [x] \ \iff \ p_{k+1} \Vdash \name z \in [x], \end{equation} which will complete the proof of Subclaim \ref{kisza}.
		
		Fix $k \leq n+1$, observe that $\bar t^k$ and $\bar t^{k+1}$ differ by at most a single digit,
			and if $t^k_k \neq t^{k+1}_k$, then $\xi_k = \varp^\iota_m$ for some $\iota \in \{\mathbf{1}, \infty, \bbS\}$, and $m$, by the assumptions of the subclaim. Therefore, 	clause \ref{obst2} from Observation \ref{obst}	 implies that there exists
		$p_k' \leq_{\{ \varp^\iota_m\}, 1} p^{(\bar u_*^{k}, \bar v_*^{k}, \bar w_*^{k}, \bar s_*^{k}) \rest k}$, and 
			$$\{ (p'_k)^{\{ \varp^\iota_m\}, (\langle 0 \rangle)}, (p'_k)^{\{ \varp^\iota_m\}, (\langle 1 \rangle)}\} = \{p_k, p_{k+1} \}.$$
			But by our construction $(\bar u_*^{k}, \bar v_*^{k}, \bar w_*^{k}, \bar s_*^{k}) \rest k = (\bar u, \bar v, \bar w, \bar s) \rest k$, so by the conditions of the subclaim $p^{(\bar u_*^{k}, \bar v_*^{k}, \bar w_*^{k}, \bar s_*^{k}) \rest k} \in D^{\iota,\mathbf{eq}}_m(\name z)$. This yields that $p_k' \in  D^{\iota,\mathbf{eq}}_m(\name z)$ by \ref{elfaop} from Fact \ref{elfa}.
			Hence by the definition of being in $D^{\iota,\mathbf{eq}}_m(\name z)$ 
			we obtain
			$$(p'_k)^{\{ \varp^\iota_m\}, (\langle 0 \rangle)} \Vdash \name z \in [x] \ \iff \  (p'_k)^{\{ \varp^\iota_m\}, (\langle 1 \rangle)}\} \Vdash \name z \in [x],$$
			and so \eqref{ezkell} holds, indeed.

	\end{PROOF}
	
	Now it only remains to construct the sequence promised in \ref{ppro}.
	Assume $k <N$, and $r^*_n \geq_{\bar{\bar{\ell}}^{n+1}} q_k$ has been already chosen.
	Recall that $m \in \omega$ is defined so that $\xi_n = \varp^\bbS_m$.
	The properties of $r_n^*$ imply that
	 $$\odot_1^m(q_k^{(\bar u^k, \bar v^k, \bar w^k, \bar s^k) \rest n}) \ \vee \ \odot_2^m(q_k^{(\bar u^k, \bar v^k, \bar w^k, \bar s^k) \rest n}),$$
	 so $$\odot_1^m(q_k^{(\bar u^k, \bar v^k, \bar w^k, \bar s^k)}) \ \vee \ \odot_2^m(q_k^{(\bar u^k, \bar v^k, \bar w^k, \bar s^k)}).$$
	 First suppose that $\odot_1^m(q_k^{(\bar u^k, \bar v^k, \bar w^k, \bar s^k)})$.
	 Then it suffices to set $q_{k+1} \leq_{\bar{\bar{\ell}}^{n+1}} q_k$, such that both $q_{k+1}^{(\bar u^k, \bar v^k, \bar w^k, \bar s^k)}$, and $q_{k+1}^{((\bar u')^k, (\bar v')^k, (\bar w')^k, (\bar s')^k)})$ decide the first 
	 $$\lh(x^{(\bar u^k, \bar v^k, \bar w^k, \bar s^k)}_k) + \lh(x^{((\bar u')^k, (\bar v')^k, (\bar w')^k, (\bar s')^k)}_k)+1$$ digits of $\name z$, ensuring that the relevant parts in \ref{fa}-\ref{ftln} from \ref{ppro} hold.
	 
	 So we can turn to the case of $\odot_2^m(q_k^{(\bar u^k, \bar v^k, \bar w^k, \bar s^k)})$.
	 By this property, there exists $q_* \leq q_k^{(\bar u^k, \bar v^k, \bar w^k, \bar s^k)}$, and $j \in \omega$, such that
	 $q_*^{\{\varp^\bbS_m\}, (\langle 0 \rangle)} \Vdash \name z_j = i_0$,  $q_*^{\{\varp^\bbS_m\}, (\langle 1 \rangle)} \Vdash \name z_j = i_1$ for some $i_0 \neq i_1$
	 (w.l.o.g.\ we can assume that
	 $$j \geq \lh(x^{(\bar u^k, \bar v^k, \bar w^k, \bar s^k)}_k) + \lh(x^{((\bar u')^k, (\bar v')^k, (\bar w')^k, (\bar s')^k)}_k)+1).$$
	 
	  W.l.o.g.\ $q_*(\varp^\infty_l)^{(w^k_l)} = q_k(\varp^\infty_l)^{(w^k_l)}$ for each $l \in \omega$,  so by Observation \ref{obskite} for some $q_k^* \leq_{\bar{\bar{\ell}}^{n+1}} q_k$ we have
	 $(q_k^*)^{(\bar u^k, \bar v^k, \bar w^k, \bar s^k)} = q_*$, and so
	 \begin{equation} \label{i00}
	 	((q_k^*)^{(\bar u^k, \bar v^k, \bar w^k, \bar s^k)})^{\{\varp^\bbS_m\}, (\langle 0 \rangle)} \Vdash \name z_j = i_0,
	 \end{equation}
 	\begin{equation} \label{i10}
 		((q_k^*)^{(\bar u^k, \bar v^k, \bar w^k, \bar s^k)})^{\{\varp^\bbS_m\}, (\langle 1 \rangle)} \Vdash \name z_j = i_1.
 	\end{equation}
	 
	 For a suitable extension $q_+ \leq (q_k^*)^{((\bar u')^k, (\bar v')^k, (\bar w')^k, (\bar s')^k)}$ 
	  \begin{equation} q_+ \Vdash \name z_j = i_* \text{ for some }i_{*} \in \{0,1\}. \end{equation}
		We need the following simple fact, which uses only that if $p$ is a Sacks condition, then when replacing $p$ by $q \leq_1 p$, then $p^{(\langle 0 \rangle)}$ and $p^{(\langle 1 \rangle)}$ can be dealt with independently.
		\begin{fact}
			If $r \in \bbQ'$, $(\bar u, \bar v, \bar w, \bar s), (\bar u', \bar v', \bar w', \bar s') \in \mathbf{seq}_{n+1}(r)$, $r_+ \leq r^{(\bar u', \bar v', \bar w', \bar s')}$ (where we demand also $r_+(\varp^\infty_l)^{(w'_l)} \leq r(\varp^\infty_l)^{(w'_l)}$ for each $l \in \omega$), and $ s_m \neq s_m'$, then there exists $r' \leq_{\bar{\bar{\ell}}} r$ such that
			\begin{itemize}
				\item $(r')^{(\bar u', \bar v', \bar w', \bar s')} = r_+$,
				\item $r'(\varp^\bbS_m)^{(s_m)} = r(\varp^\bbS_m)^{(s_m)}$, in particular
				$$ ((r')^{(\bar u, \bar v, \bar w, \bar s)})^{\{\varp^\bbS_m\}, (\langle 0 \rangle)} \leq (r^{(\bar u, \bar v, \bar w, \bar s)})^{\{\varp^\bbS_m\}, (\langle 0 \rangle)},$$
				$$ ((r')^{(\bar u, \bar v, \bar w, \bar s)})^{\{\varp^\bbS_m\}, (\langle 1 \rangle)} \leq (r^{(\bar u, \bar v, \bar w, \bar s)})^{\{\varp^\bbS_m\}, (\langle 1 \rangle)}.$$
			\end{itemize}
		\end{fact}
		Applying the fact to $q_k^*$ and $q_+$ yields the condition $q_k^{**}$ with
		\begin{equation} (q_k^{**})^{((\bar u')^k, (\bar v')^k, (\bar w')^k, (\bar s')^k)} = q_+, \end{equation}
		and so
		\begin{equation} \label{q**} (q_k^{**})^{((\bar u')^k, (\bar v')^k, (\bar w')^k, (\bar s')^k)} \Vdash \name z_j = i_*. \end{equation}
		Furthermore, the fact gives us that 
		$$ (q_k^{**})^{(\bar u^k, \bar v^k, \bar w^k, \bar s^k)} \leq_{\{\varp^\bbS_m\}, 1}(q_k^*)^{(\bar u^k, \bar v^k, \bar w^k, \bar s^k)},$$
		so by \eqref{i00} and \eqref{i10}
		 \begin{equation} \label{i0}
			((q_k^{**})^{(\bar u^k, \bar v^k, \bar w^k, \bar s^k)})^{\{\varp^\bbS_m\}, (\langle 0 \rangle)} \Vdash \name z_j = i_0,
		\end{equation}
		\begin{equation} \label{i1}
			((q_k^{**})^{(\bar u^k, \bar v^k, \bar w^k, \bar s^k)})^{\{\varp^\bbS_m\}, (\langle 1 \rangle)} \Vdash \name z_j = i_1.
		\end{equation}
		Now since $i_0 \neq i_1$ either $i_0 \neq i_*$, or $i_1 \neq i_*$, w.l.o.g. we can assume that 
			\begin{equation} i_0 \neq i_*.
		\end{equation}
		Finally, appealing to Observation \ref{obskite} again, there exists $q_k^{***} \leq_{\bar{\bar{\ell}}^{n+1}} q_k^{**}$, such that 
		$$(q_k^{***})^{(\bar u^k, \bar v^k, \bar w^k, \bar s^k)} = (q_k^{**})^{(\bar u^k, \bar v^k, \bar w^k, \bar s^k)})^{\{\varp^\bbS_m\}, (\langle 0 \rangle)},$$
		so
		\begin{equation} (q_k^{***})^{(\bar u^k, \bar v^k, \bar w^k, \bar s^k)} \Vdash \name z_j \neq i_*. \end{equation}
		This, together with \eqref{q**} shows that setting $q_{k+1} = q_k^{***}$ works, 
		since by possibly replacing $q_k^{***}$ with a $\leq_{\bar{\bar{\ell}}^{n+1}}$-extension w.l.o.g.\ both $(q_k^{***})^{(\bar u^k, \bar v^k, \bar w^k, \bar s^k)}$ and  $(q_k^{***})^{((\bar u')^k, (\bar v')^k, (\bar w')^k, (\bar s')^k)}$ decide $\name z \rest [0,j]$.

		\begin{enumerate}[label = $(\blacklozenge_{3})(\roman*)$, ref = $(\blacklozenge_{3})(\roman*)$]
		\setcounter{enumi}{2}
		\item \label{caseiii} Case $iii$: $\xi_n = \varp^\iota_m$ for some $\iota \in \{ \mathbf{1}, \infty\}$ (and $m$).
		\end{enumerate}

	\begin{lemma} \label{guszkl0}
		Let $q_* \in \bbQ'$, $\Vdash_{\bbQ'} \name z \in 2^\omega$, $n \in \omega$ be given,
		suppose that $\iota \in \{ \mathbf{1}, \infty\}$, $m \in \omega$ are such that $\xi_n = \varp^\iota_m$.
		\Then \ there exists 
		\begin{enumerate}[label = $(\roman*)$, ref = $(\roman*)$]
			\item $r_* \leq_{\bar{\bar{\ell}}^n} q_*$
			\item $\bar{\mathbf{x}} = \langle \mathbf{x}_{(\bar u, \bar v, \bar w, \bar s)}: \ (\bar u, \bar v, \bar w, \bar s) \in \mathbf{seq}_{\bar{\bar{\ell}}^n}(q_*) \rangle$
			\item $\bar{x} = \langle x^{(\bar u, \bar v, \bar w, \bar s)}: \ (\bar u, \bar v, \bar w, \bar s) \in \mathbf{seq}_{\bar{\bar{\ell}}^{n+1}}(r_*) \rangle$
		\end{enumerate}
		satisfying the following.
				
		\begin{enumerate}[label = $(\boxtimes)$, ref = $(\boxtimes)$]
			\item 	\label{bbf1} for each $(\bar u, \bar v, \bar w, \bar s)$ we have
			\begin{itemize}
				\item  $x^{(\bar u, \bar v, \bar w, \bar s)} \in \ ^{<\omega} 2$, and
				\item  $\mathbf{x}_{(\bar u, \bar v, \bar w, \bar s)\rest n} \in \ \{ \mathbf{mul}, \mathbf{eq},\mathbf{un}\}$,
			\end{itemize}  such that $r_*^{(\bar u, \bar v, \bar w, \bar s)\rest n)} \in D^{\iota, \mathbf{x}_{(\bar u, \bar v, \bar w, \bar s)\rest n}}_m(\name z)$, and
			 whenever $(\bar u', \bar v', \bar w', \bar s') \neq (\bar u, \bar v, \bar w, \bar s)$, and $\mathbf{x}_{(\bar u, \bar v, \bar w, \bar s)\rest n} \neq \mathbf{eq}$, then either
			 \begin{enumerate}[label = $(\boxtimes)_{(\alph*)}$, ref = $(\boxtimes)_{(\alph*)}$]
				\item $\mathbf{x}_{(\bar u, \bar v, \bar w, \bar s)\rest n} = \mathbf{un}$, and then $\iota = \mathbf{1}$ holds, when for the unique $k$ satisfying 
				$$k <\min(\lh(x^{(\bar u, \bar v, \bar w, \bar s)}), \lh(x^{(\bar u', \bar v', \bar w', \bar s')}))$$
				\[ x^{(\bar u, \bar v, \bar w, \bar s)}_k \neq x_k^{(\bar u', \bar v', \bar w', \bar s')},\]
			we have $$x^{(\bar u, \bar v, \bar w, \bar s)} \rest k = x^{(\bar u, \bar v, \bar w, \bar s)} \rest k \neq \mathfrak{s}_k$$ (where $\bar{\mathfrak{s}}$ is from \ref{x1}-\ref{x2}),
			\item \label{mul}  or $\mathbf{x}_{(\bar u, \bar v, \bar w, \bar s)\rest n} = \mathbf{mul}$, and 
				$$k_0< k_1 <\min(\lh(x^{(\bar u, \bar v, \bar w, \bar s)}), \lh(x^{(\bar u', \bar v', \bar w', \bar s')})$$
			\[ x^{(\bar u, \bar v, \bar w, \bar s)}_{k_0} \neq x_{k_0}^{(\bar u', \bar v', \bar w', \bar s')},\]
			\[ x^{(\bar u, \bar v, \bar w, \bar s)}_{k_1} \neq x_{k_1}^{(\bar u', \bar v', \bar w', \bar s')}.\] 	
			\end{enumerate} 
					
		\end{enumerate}	
	
	\end{lemma}
	\begin{observation} \label{elfaop+}
		If $n$, $\iota$, $m$ are as in  the lemma, and $r'$ satisfies
		$$\forall  (\bar u, \bar v, \bar w, \bar s) \in \mathbf{seq}_{\bar{\bar{\ell}}^n}(r'): \
		(r')^{(\bar u, \bar v, \bar w, \bar s)} \in D^{\iota, \mathbf{x}^{(\bar u, \bar v, \bar w, \bar s)}}_m(\name z),$$
		 then the same statements holds for any $r'' \leq_{\bar{\bar{\ell}}^{n+1}} r'$.
	\end{observation}
	Before proving the lemma note that 
		\begin{enumerate}[label = $(\blacklozenge_{3})(iii)_\arabic*$, ref = $(\blacklozenge_{3})(iii)_\arabic*$]
		\item \label{ppppro}
	applying to $q_* = r_n$ it yields the desired condition $r_{n+1} \leq_{\bar{\bar{\ell}}^{n}} r_n$, $\mathbf{x}_{(\bar u, \bar v, \bar w, \bar s)}$'s ($(\bar u, \bar v, \bar w, \bar s) \in \mathbf{seq}_{\bar{\bar{\ell}}^n}(r_{n+1})$, and setting \begin{itemize}
		\item 	$y^{(\bar u, \bar v, \bar w, \bar s)} = \text{ the maximal element in }$
		$$ \{x \in \ ^{\omega>}2 : \ r_{n+1}^{(\bar u, \bar v, \bar w, \bar s)} \Vdash x \in \name z \} (\supseteq x^{(\bar u, \bar v, \bar w, \bar s)}_N)$$
		\end{itemize}
	the requirements in \ref{r0}-\ref{r2} are clearly satisfied: just use the same argument as after \ref{pppro}, therefore finishing the case \ref{caseiii}, and the induction in \ref{nemkez}, too.
		\end{enumerate} 
	
	\begin{PROOF}{Lemma \ref{guszkl0}}(Lemma \ref{guszkl0})
		We are going to construct $r_*$ regardless of the specific value of $\iota \in \{ \mathbf{1}, \infty\}$. We remark that although for $\iota = \infty$ a simpler argument would also suffice, as the case of $\iota = \mathbf{1}$ itself  needs a slightly more involved (and painful) reasoning, it is easier to handle the two together.
		First
		\begin{enumerate}[label = $(\blacktriangle)_{\arabic*}$, ref = $(\blacktriangle)_{\arabic*}$] 
			\item we choose an enumeration $\langle (\bar u^j, \bar v^j, \bar w^j, \bar s^j): j<M \rangle$ of all the possible quadruples $(\bar u, \bar v, \bar w, \bar s)$ from the set $\mathbf{seq}_{\bar{\bar{\ell}}^{n}}(q_*)$.
		\end{enumerate}
		We need the following.
		
		\begin{enumerate}[label = $(\blacktriangle)_{\arabic*}$, ref = $(\blacktriangle)_{\arabic*}$]
			\stepcounter{enumi} 
			\item \label{hulye} We are going to define $q_{**} \leq_{\bar{\bar{\ell}}^n} q_*$, 
			as well as the sequence 
			$$\bar{\mathbf{x}} = \langle  \mathbf{x}_{(\bar u, \bar v, \bar w, \bar s)}: \ (\bar u, \bar v, \bar w, \bar s) \in \mathbf{seq}_{\bar{\bar{\ell}}^n}(q_*) \rangle,$$ with $\mathbf{x}_{(\bar u, \bar v, \bar w, \bar s)} \in  \{$\textbf{un},\textbf{mul},\textbf{eq}$\}$ satisfying the following.
			(For $j<M$ writing sometimes $\mathbf{x}_j$ instead of $\mathbf{x}_{(\bar u^j, \bar v^j, \bar w^j, \bar s^j)}$) we would like the sequences to have the properties as follows:
			\begin{enumerate}
				\item $\bigcap_{j<M} \{q \leq_{\bar{\bar{\ell}}^n} q_{**}: \ q^{(\bar u^j, \bar v^j, \bar w^j, \bar s^j)} \in D_m^{\iota,\mathbf{x}_j}(\name z)\}$ is $\leq_{\bar{\bar{\ell}}^n}$-dense below $q_{**}$, 
				\item moreover,
				whenever $\mathbf{x}_i = \mathbf{mul}$ for some $i<M$, then for every $p \leq_{\bar{\bar{\ell}}^n} q_{**}$
				$$ \begin{array}{l} p \in \bigcap_{j<i} \{q \leq_{\bar{\bar{\ell}}^n} q_{**}: \ q^{(\bar u^j, \bar v^j, \bar w^j, \bar s^j)} \in D_m^{\iota,\mathbf{x}_j}(\name z)\} \\ \Rightarrow \\  p^{(\bar u^i, \bar v^i, \bar w^i, \bar s^i)} \in D_m^{\iota,\mathbf{mul}}(\name z). \end{array}$$
			\end{enumerate}
		\end{enumerate}
		
		\begin{sclaim} \label{sor}
			Suppose that $\iota$, $m$, $n$, $q_*$ are as in the lemma. Then there exist $q_{**} \leq_{\bar{\bar{\ell}}^n} q_*$ and
			$\bar{\mathbf{x}}$ satisfying the requirements in \ref{hulye}.
		\end{sclaim}
		\begin{PROOF}{Subclaim \ref{sor}}
			First we define $q_0$, $\mathbf{x}_0$ as follows.
			\begin{enumerate}[label = $(\intercal)_{\arabic*}$,ref = $(\intercal)_{\arabic*}$]
				\item Set the auxiliary variable $q_{0} = q_*$.
				\begin{enumerate}[label = $\bullet_{\arabic*}$, ref = $\bullet_{\arabic*}$]
					\item First, suppose that the set
					$$\{q \leq_{\bar{\bar{\ell}}^n} q_{0}: \ q^{(\bar u^0, \bar v^0, \bar w^0, \bar s^0)} \in D_m^{\iota,\mathbf{un}}(\name z)\}$$
					is $\leq_{\bar{\bar{\ell}}^n}$-dense below $q_{0}$, in which case set $\mathbf{x}_0 = \mathbf{un}$, $q_1 = q_{0}$.
					\item Otherwise, define $q_{0}' \leq_{\bar{\bar{\ell}}^n}q_{0}$ so that there is no $q \leq_{\bar{\bar{\ell}}^n} q_{0}'$ with 
					$q^{(\bar u^0, \bar v^0, \bar w^0, \bar s^0)} \in D_m^{\iota,\mathbf{un}}(\name z)$.
					\item Second, if the set
					$$\{q \leq_{\bar{\bar{\ell}}^n} q'_{0}: \ q^{(\bar u^0, \bar v^0, \bar w^0, \bar s^0)} \in D_m^{\iota,\mathbf{eq}}(\name z)\}$$
					is  $\leq_{\bar{\bar{\ell}}^n}$-dense below $q'_{0}$, then we let $\mathbf{x}_0 = \mathbf{eq}$, and $q_1 = q'_{0}$.
					\item  If it is not the case, then there is
					$q_{0}'' \leq_{\bar{\bar{\ell}}^n}q_0'$ for which there is no $q \leq_{\bar{\bar{\ell}}^n} q_{0}''$ with 
					$q^{(\bar u^0, \bar v^0, \bar w^0, \bar s^0)} \in D_m^{\iota,\mathbf{eq}}(\name z)$. 
				\end{enumerate}	  
			\end{enumerate}
			Then
			\begin{enumerate}[label = $(\intercal)_{\arabic*}$,ref = $(\intercal)_{\arabic*}$]
				\setcounter{enumi}{1}
				\item \label{k1} set $\mathbf{x}_0 = \mathbf{mul}$, and $q_1 = q''_{0}$,
			\end{enumerate}
			and observe that by Definition \ref{defmul} 
			\begin{enumerate}[label = $(\intercal)_{\arabic*}$,ref = $(\intercal)_{\arabic*}$]
				\setcounter{enumi}{2} 
				\item \label{k2}	for each 	$q \leq_{\bar{\bar{\ell}}^n}q_{0}'' = q_1$ we have
				$q^{(\bar u^0, \bar v^0, \bar w^0, \bar s^0)} \in D_m^{\iota,\mathbf{mul}}(\name z)$.
			\end{enumerate}
			
			\begin{enumerate}[label = $(\intercal)_{\arabic*}$,ref = $(\intercal)_{\arabic*}$]
				\setcounter{enumi}{3} 
				\item 	This way  we are going to define 
				\begin{itemize}
					\item the $\leq_{\bar{\bar{\ell}}^n}$-decreasing sequence 
					$$ q_* = q^0 \geq_{\bar{\bar{\ell}}^n} q^1 \geq_{\bar{\bar{\ell}}^n} \dots \geq_{\bar{\bar{\ell}}^n} q^{M},$$
					\item together with the sequence $\bar{\mathbf{x}}$
				\end{itemize}
				by induction on  $j$ such that for each $i \leq M$
			\end{enumerate}
			\begin{enumerate}[label = $(\intercal)_{4}^{(\alph*)}$,ref = $(\intercal)_{4}^{(\alph*)}$]
				\item \label{k3a} $\bigcap_{j<i} \{q \leq_{\bar{\bar{\ell}}^n} q_{i}: \ q^{(\bar u^j, \bar v^j, \bar w^j, \bar s^j)} \in D_m^{\iota,\mathbf{x}_j}(\name z)\}$ is $\leq_{\bar{\bar{\ell}}^n}$-dense below $q_{i}$, 
				\item \label{k3b}and if $\mathbf{x}_{i-1} = \mathbf{mul}$, then for arbitrary $p \leq_{\bar{\bar{\ell}}^n} q_{i}$:
				$$ \begin{array}{l} p \in \bigcap_{j<i-1} \{q \leq_{\bar{\bar{\ell}}^n} q_i: \ q^{(\bar u^j, \bar v^j, \bar w^j, \bar s^j)} \in D_m^{\iota,\mathbf{x}_j}(\name z)\} \\ \Rightarrow \\  p^{(\bar u^{i-1}, \bar v^{i-1}, \bar w^{i-1}, \bar s^{i-1})} \in D_m^{\iota,\mathbf{mul}}(\name z). \end{array}$$
			\end{enumerate}
			Note that $q_1$ and $\mathbf{x}_0$ clearly satisfy the demands if $\mathbf{x}_0 \in \{\mathbf{un}, \mathbf{eq}\}$, and also when $\mathbf{x}_0= \mathbf{mul}$, for which recall \ref{k1}, \ref{k2}.
			\begin{enumerate}[label = $(\intercal)_{\arabic*}$,ref = $(\intercal)_{\arabic*}$]
				\setcounter{enumi}{4} 
				\item	Suppose that $0<i<M$, and $q_i$, and the $\mathbf{x}_j$'s are already defined for $j<i$. 
				\begin{enumerate}[label = $\bullet_{\arabic*}$, ref = $\bullet_{\arabic*}$]
					\item \label{D*}	Set 
					$$D^* = \bigcap_{j<i} \{q \leq_{\bar{\bar{\ell}}^n} q_i: \ q^{(\bar u^j, \bar v^j, \bar w^j, \bar s^j)} \in D_m^{\iota,\mathbf{x}_j}(\name z)\},$$
					which is $\leq_{\bar{\bar{\ell}}^n}$-dense below $q_{i}$.
					\item \label{l*}	Note that
					$$q \leq_{\bar{\bar{\ell}}^{n+1}} p \ \iff \ (\forall j<M) q^{(\bar u^j, \bar v^j, \bar w^j, \bar s^j)} \leq_{\{\varp^\iota_m\},1} p^{(\bar u^j, \bar v^j, \bar w^j, \bar s^j)}).$$
					Using this observation (recalling \ref{elfaop} from Fact \ref{elfa}) we obtain
					\begin{equation}  \label{*op} (\forall q',q \in \bbQ'): \ (q' \leq_{\bar{\bar{\ell}}^{n+1}} q, \ q \in D^*) \ \to \ (q' \in D^*). \end{equation}
				\end{enumerate}
				Again, following the pattern of the definition of $q_1$ and $\mathbf{x}_0$,
				
				\begin{enumerate}[label = $\bullet_{\arabic*}$, ref = $\bullet_{\arabic*}$]
					\setcounter{enumii}{2}
					\item  if 
					$$\{q \in D^*: \ q \leq_{\bar{\bar{\ell}}^n} q_{i}, \ q^{(\bar u^i, \bar v^i, \bar w^i, \bar s^i)} \in D_m^{\iota,\mathbf{un}}(\name z)\}$$
					is $\leq{\bar{\bar{\ell}}^n}$-dense below $q_i$, then we set $\mathbf{x}_i = \mathbf{un}$, and $q_{i+1} = q_{i}$,
					\item  otherwise let $q_{i}' \leq_{\bar{\bar{\ell}}^n} q_{i}$ be such that for no $q \leq_{\bar{\bar{\ell}}^n} q_{i}'$ do we have $q^{(\bar u^i, \bar v^i, \bar w^i, \bar s^i)} \in D_m^{\iota,\mathbf{un}}(\name z)$.
					\item If the set
					$$\{q \in D^*: \ q \leq_{\bar{\bar{\ell}}^n} q'_{i}: \ q^{(\bar u^i, \bar v^i, \bar w^i, \bar s^i)} \in D_m^{\iota,\mathbf{eq}}(\name z)\}$$
					is  $\leq_{\bar{\bar{\ell}}^n}$-dense below $q'_{i}$, then we let $\mathbf{x}_i = \mathbf{eq}$, and $q_{i+1} = q'_{i}$.
					\item  If it is not the case, then there is $q''_{i} \leq_{\bar{\bar{\ell}}^n} q'_{i}$, such that 
					\[ \forall q \in D^*: \ (q \leq_{\bar{\bar{\ell}}^n} q''_{i}) \ \to \ (q \notin D_m^{\iota,\mathbf{un}}(\name z) \cup D_m^{\iota,\mathbf{eq}}(\name z)).  \]
					But then by $\eqref{*op}$, \ref{D*} and by Definition \ref{defmul} we have that
					\[ \forall q \in D^*: \ (q \leq_{\bar{\bar{\ell}}^n} q''_{i}) \ \to \ (q \in D_m^{\iota,\mathbf{mul}}(\name z)),  \]
					and letting $q_{i+1} = q_{i}''$, $\mathbf{x}_{i} = \mathbf{mul}$ we are done, \ref{k3a}, \ref{k3b} hold for $i_* = i+1$.
				\end{enumerate}
				
			\end{enumerate}
			Finally, letting $q_{**} = q_{M}$, it is easy to check that $q_{**}$, $\bar{\mathbf{x}}$ are as desired.
		\end{PROOF}
		
		\begin{sclaim} \label{fasclelo}
			Suppose that  $\iota \in \{ \mathbf{1}, \infty\}$, $l$, $m$, $n$ are as in Lemma \ref{guszkl0},  $\bar{\mathbf{x}}$, $q_{**}$ given by Subclaim \ref{sor} (i.e.\ satisfying the requirements in \ref{hulye}, w.r.t.\ the fixed enumeration of $\mathbf{seq}_{\bar{\bar{\ell}}^n}(q_{*})= \mathbf{seq}_{\bar{\bar{\ell}}^n}(q_{**})$).
			Then there exists a sequence
			$$q_{**} \geq_{\bar{\bar{\ell}}^n} q_0 \geq_{\bar{\bar{\ell}}^n} q_1 \geq_{\bar{\bar{\ell}}^n} \dots  \geq_{\bar{\bar{\ell}}^n} q_l$$
			such that 
			\begin{enumerate}
				\item  $(\forall k\leq l) \ q_k \in \bigcap_{j<M} \{q \leq_{\bar{\bar{\ell}}^n} q_{**}: \ q^{(\bar u^j, \bar v^j, \bar w^j, \bar s^j)} \in D_m^{\iota,\mathbf{x}_j}(\name z)\},$
				\item  for each $k<l$, $j<M$: 
				$$q_{k+1}^{(\bar u^j, \bar v^j, \bar w^j, \bar s^j)} \leq (q_k^{(\bar u^j, \bar v^j, \bar w^j, \bar s^j)})^{\{\varp^\iota_m\}, (\langle 0 \rangle)},$$
				\item 	and for each $j<M$, if $\mathbf{x}_j = \mathbf{un}$, then there exist $i^j_0 < i^j_1 < \dots < i^j_l$ (and $c^j_k \in \{0,1\}$), such that
				 $$\begin{array}{lll} \forall k\leq l:  & (q_k^{(\bar u^j, \bar v^j, \bar w^j, \bar s^j)})^{\{\varp^\iota_m\}, (\langle 0 \rangle)} & \Vdash \name z_{i^j_k} = c^j_k, \\
				 & (q_k^{(\bar u^j, \bar v^j, \bar w^j, \bar s^j)})^{\{\varp^\iota_m\}, (\langle 1 \rangle)} & \Vdash \name z_{i^j_k} = 1- c^j_k, \end{array}$$
			\end{enumerate}
		\end{sclaim}
		\begin{PROOF}{Subclaim \ref{fasclelo}}
			By \ref{hulye} we can choose $p_0 \leq_{\bar{\bar{\ell}}^n} q_{**}$ with 
			\begin{equation} \label{ae} p_0 \in D_*:= \bigcap_{j<M}\{p \leq_{\bar{\bar{\ell}}^n} q_{**}: \ p^{(\bar u^j, \bar v^j, \bar w^j, \bar s^j)} \in D_m^{\iota,\mathbf{x}_j}(\name z) \}. \end{equation}
			This means that
			\newcounter{elokccc}
			 \begin{enumerate}[label = $(\ddagger)_{\arabic*}$, ref = $(\ddagger)_{\arabic*}$]
				\setcounter{enumi}{\value{elokccc}}
				\item  \label{wen} whenever $j<M$ is such that $\mathbf{x}_j = \mathbf{un}$, then for some $i^j_0 \in \omega$, $c^j_0 \in \{0,1\}$:
				$$\begin{array}{ll}  (p_0^{(\bar u^j, \bar v^j, \bar w^j, \bar s^j)})^{\{\varp^\iota_m\}, (\langle 0 \rangle)} & \Vdash \name z_{i^j_0} = c^j_0, \\
					(p_0^{(\bar u^j, \bar v^j, \bar w^j, \bar s^j)})^{\{\varp^\iota_m\}, (\langle 1 \rangle)} & \Vdash \name z_{i^j_0} = 1- c^j_0. \end{array}$$
				\stepcounter{elokccc}
			\end{enumerate}

			 Now for each such fixed $j$ there is $p' \leq_{\{\varp^\iota_m\}, 1}$ $q_0^{(\bar u^j, \bar v^j, \bar w^j, \bar s^j)}$ so that both $(p')^{\{\varp^\iota_m\}, (\langle 0 \rangle)}$ and $(p')^{\{\varp^\iota_m\}, (\langle 1 \rangle)}$ decide $\name z \rest [0, i^j_0)$ (in fact $q^{(\bar u^j, \bar v^j, \bar w^j, \bar s^j)} \in D_m^{\iota,\mathbf{un}}(\name z)$ implies that then $p'$ decides that initial segment). Therefore,
			 \begin{enumerate}[label = $(\ddagger)_{\arabic*}$, ref = $(\ddagger)_{\arabic*}$]
			 		\setcounter{enumi}{\value{elokccc}}
			 	\item  \label{ex} there exists $q_0 \leq_{\bar{\bar{\ell}}^{n+1}} p_0$, for which 
			 	\begin{equation} \label{aea} \forall j<M: \ q_0^{(\bar u^j, \bar v^j, \bar w^j, \bar s^j)} \parallel \name z \rest [0, i^j_0), \end{equation}
			 	and (automatically by \eqref{ae}  $q_0 \leq_{\bar{\bar{\ell}}^{n+1}} p_0$)
			 	\begin{equation} \label{aeae} q_0 \in D_*. \end{equation}
			 		\stepcounter{elokccc}
			 \end{enumerate}

			 Now let $p_1 \leq_{\bar{\bar{\ell}}^n} q_0$ be such that  
			 \begin{equation} \label{p1} \text{for each }j<M: \ p_1^{(\bar u^j, \bar v^j, \bar w^j, \bar s^j)} \leq (q_0^{(\bar u^j, \bar v^j, \bar w^j, \bar s^j)})^{\{\varp^\iota_m\}, (\langle 0 \rangle)} \end{equation} 
			 (in particular, $p_1 \nleq_{\bar{\bar{\ell}}^{n+1}} q_0$), and $p_1 \in D_*$.
			 Similarly to \ref{wen} and \ref{ex}
		 \begin{enumerate}[label = $(\ddagger)_{\arabic*}$, ref = $(\ddagger)_{\arabic*}$]
			 			\setcounter{enumi}{\value{elokccc}}
			 	\item  whenever $j<M$ is such that $\mathbf{x}_j = \mathbf{un}$, then for some $i^j_1 \in \omega$, $c^j_1 \in \{0,1\}$:
			 	$$\begin{array}{ll}  (p_1^{(\bar u^j, \bar v^j, \bar w^j, \bar s^j)})^{\{\varp^\iota_m\}, (\langle 0 \rangle)} & \Vdash \name z_{i^j_1} = c^j_1, \\
			 		(p_1^{(\bar u^j, \bar v^j, \bar w^j, \bar s^j)})^{\{\varp^\iota_m\}, (\langle 1 \rangle)} & \Vdash \name z_{i^j_1} = 1- c^j_1, \end{array}$$
		 			\stepcounter{elokccc}
			 \end{enumerate}
			 \begin{enumerate}[label = $(\ddagger)_{\arabic*}$, ref = $(\ddagger)_{\arabic*}$]
			 		\setcounter{enumi}{\value{elokccc}}
			 	\item   there exists $q_1 \leq_{\bar{\bar{\ell}}^{n+1}} p_1$, for which 
			 	\begin{equation} \label{aeaea} \forall j<M: \ q_1^{(\bar u^j, \bar v^j, \bar w^j, \bar s^j)} \parallel \name z \rest [0, i^j_1), \end{equation}
			 	and (automatically by  $q_1 \leq_{\bar{\bar{\ell}}^{n+1}} p_1$)
			 	\begin{equation} \label{aeaeae} q_1 \in D_*. \end{equation}
			 		\stepcounter{elokccc}
			 \end{enumerate}
		 	Observe that by \eqref{p1} and \eqref{aea}
		 	\begin{enumerate}[label = $(\ddagger)_{\arabic*}$, ref = $(\ddagger)_{\arabic*}$]
		 		\setcounter{enumi}{\value{elokccc}}
		 		\item  for each $j$: $i^j_0 < i^j_1$.
		 		\stepcounter{elokccc}
		 	\end{enumerate}
	 		Following this pattern, we can define the sequence by induction on $k \leq l$.

		\end{PROOF}

		\begin{sclaim} \label{fascl}
			Suppose that  $\iota \in \{ \mathbf{1}, \infty\}$, $n$, $m$, are as in Lemma \ref{guszkl0},  $\bar{\mathbf{x}}$, $q_{**}$ given by Subclaim \ref{sor} (i.e.\ satisfying the requirements in \ref{hulye}, w.r.t.\ the fixed enumeration of $\mathbf{seq}_{\bar{\bar{\ell}}^n}(q_{*})= \mathbf{seq}_{\bar{\bar{\ell}}^n}(q_{**})$).
			Suppose that 
			$$q_{**} \geq_{\bar{\bar{\ell}}^n} q_0 \geq_{\bar{\bar{\ell}}^n} q_1 \geq_{\bar{\bar{\ell}}^n} \dots  \geq_{\bar{\bar{\ell}}^n} q_l$$
			
			such that  the $q_k$'s are given by Subclaim \ref{fasclelo}, so
			\begin{equation} \label{ellza} (\forall k\leq l) \ q_k \in \bigcap_{j<M} \{q \leq_{\bar{\bar{\ell}}^n} q_{**}: \ q^{(\bar u^j, \bar v^j, \bar w^j, \bar s^j)} \in D_m^{\iota,\mathbf{x}_j}(\name z)\}, \end{equation}
			moreover,  
			\begin{equation} \label{elag} \text{for each }k<l: \ \forall j <M: \ q_{k+1}^{(\bar u^j, \bar v^j, \bar w^j, \bar s^j)} \leq (q_{k}^{(\bar u^j, \bar v^j, \bar w^j, \bar s^j)})^{\{\varp^\iota_m\},(\langle 0 \rangle)} \end{equation} 
			(in particular, $q_k \ngeq_{\bar{\bar{\ell}}^{n+1}} q_{k+1}$).
			
			Then,
			\begin{enumerate}
				\item  there exists $q_{***} \leq_{\bar{\bar{\ell}}^n} q_0$ (in fact, even $q_{***} \leq_{\bar{\bar{\ell}}^{n+1}} q_0$), for which 
				\[ \forall j < M, \ \forall \bar t \in \ ^{k\geq}2:  \ (q_{***}^{(\bar u^j, \bar v^j, \bar w^j, \bar s^j)})^{\{\varp^\iota_m\}, (\bar t)} \in D_m^{\iota,\mathbf{x}_j}(\name z)\},\]
				and if $j$ is such that $\mathbf{x}_j = \mathbf{un}$, then for the sequence  $i^j_0 < i^j_1 < \dots < i^j_l$ from Subclaim \ref{fasclelo}:
				\begin{equation} \label{unifo}
					\begin{array}{rll}
								 \forall \bar t \in \ ^{k\geq}2 \ \exists c, z^*: &  (q_{***}^{(\bar u^j, \bar v^j, \bar w^j, \bar s^j)})^{\{\varp^\iota_m\}, (\bar t \tieconcat \langle 0 \rangle )} & \Vdash  \name z_{i^j_k} = c, \\
								 		 & (q_{***}^{(\bar u^j, \bar v^j, \bar w^j, \bar s^j)})^{\{\varp^\iota_m\}, (\bar t \tieconcat \langle 1 \rangle )} & \Vdash  \name z_{i^j_k} = 1- c,\\
								 &	(q_{***}^{(\bar u^j, \bar v^j, \bar w^j, \bar s^j)})^{\{\varp^\iota_m\}, (\bar t)} & \Vdash  \name z \rest [0,i^j_k) 	  = z^*.
					\end{array}			
				\end{equation}
				\item \label{more}	Moreover, if for each $k$ ``$=$" holds in \eqref{elag}, then $q_{***}$ can be chosen to be $q_0$.
			\end{enumerate}

		\end{sclaim}
		\begin{PROOF}{fascl}
			Condition $\eqref{elag}$ implies that
			we can define the condition $q^*_{l-1} \leq_{\bar{\bar{\ell}}^{n+1}} q_{l-1}$ so that 
			\[ (\forall j<M): \  q_{l}^{(\bar u^j, \bar v^j, \bar w^j, \bar s^j)} = (q_{l-1}^*)^{(\bar u^j, \bar v^j, \bar w^j, \bar s^j)})^{\{\varp^\iota_m\}, (\langle 0 \rangle)}\]
			\newcounter{csicou} \setcounter{csicou}{0}
			\begin{enumerate}[label = $(\circledcirc)_{\arabic*}$ , ref = $(\circledcirc)_{\arabic*}$]
				\setcounter{enumi}{\value{csicou}}
				\item 
				note that replacing $q_{l-1}$ with $q^*_{l-1}$ still preserves 
				$q_{l-1} \geq_{\bar{\bar{\ell}}^n} q_{l}$, but $q_{l-1} \ngeq_{\bar{\bar{\ell}}^*{n+1}} q_{l}$, that is, $q^*_{l-1} \geq_{\bar{\bar{\ell}}^n} q^*_l$, but $q^*_{l-1} \ngeq_{\bar{\bar{\ell}}^{n+1}} q^*_l$. Similarly, for $k = l-1$  \eqref{ellza}  holds recalling that $D_m^{\iota, \mathbf{y}}$'s are closed under $\leq_{\{\varp^\iota_m\},1}$-extensions if $\mathbf{y} \in \{ \mathbf{un}, \mathbf{eq}, \mathbf{mul}\}$ (\ref{elfaop} from Fact \ref{elfa}).
				\stepcounter{csicou}
			\end{enumerate}
			Doing this by downward induction on $k =l-1$, $l-2$, $\dots$, $0$,
			\begin{enumerate}[label = $(\circledcirc)_{\arabic*}$ , ref = $(\circledcirc)_{\arabic*}$]
				\setcounter{enumi}{\value{csicou}}
				\item \label{wlog}
				replacing $q_k$ by $q^*_{k} \leq_{\bar{\bar{\ell}}^{n+1}} q_{k}$ when necessary w.l.o.g.\  we can assume that 
				\[ (\forall k<l)(\forall j<M): \  q_{k+1}^{(\bar u^j, \bar v^j, \bar w^j, \bar s^j)} = ((q_{k})^{(\bar u^j, \bar v^j, \bar w^j, \bar s^j)})^{\{\varp^\iota_m\}, (\langle 0 \rangle)},\]
				\stepcounter{csicou}
			\end{enumerate}
			and introducing the
			\begin{enumerate}[label = $(\circledcirc)_{\arabic*}$ , ref = $(\circledcirc)_{\arabic*}$]
				\setcounter{enumi}{\value{csicou}}
				\item \label{wlog2}
				notation $\vec{0}^k$ for the constant $0$ sequence of length $k$, i.e. $\vec{0}^1 = \langle 0 \rangle$, $\vec{0}^{k+1} = \vec{0}^k \tieconcat \langle 0 \rangle$,
				\[ (\forall k<l)(\forall j<M): \  q_{k}^{(\bar u^j, \bar v^j, \bar w^j, \bar s^j)} = ((q_{0})^{(\bar u^j, \bar v^j, \bar w^j, \bar s^j)})^{\{\varp^\iota_m\}, (\vec{0}^k)}.\]
				
				\stepcounter{csicou}
			\end{enumerate}
			
			Now
			\begin{enumerate}[label = $(\circledcirc)_{\arabic*}$ , ref = $(\circledcirc)_{\arabic*}$] 
				\setcounter{enumi}{\value{csicou}}
				\item we claim that (assuming \ref{wlog}) choosing $q_{***} = q_0$ works.
				\stepcounter{csicou}
			\end{enumerate}
			
			\begin{enumerate}[label = $(\circledcirc)_{\arabic*}$ , ref = $(\circledcirc)_{\arabic*}$]
				\setcounter{enumi}{\value{csicou}}
				\item \label{4} Fix $j<M$, we are going to prove that 
				\[ (\forall \bar t \in \ ^{l\geq}2): \ \ (q_0^{(\bar u^j, \bar v^j, \bar w^j, \bar s^j)})^{\{\varp^\iota_m\},(\bar t)} \in D^{\iota,\mathbf{x}_j}_m(\name z),\]
				\stepcounter{csicou}	
			\end{enumerate}
			\begin{enumerate}[label = $(\circledcirc)_{\arabic*}$ , ref = $(\circledcirc)_{\arabic*}$]
				\setcounter{enumi}{\value{csicou}}
				\item \label{arl}
				First we argue \ref{4} for $j$'s such that $\mathbf{x}_j \in \{\mathbf{un},\mathbf{eq}\}$, and prove \eqref{unifo} as well.
				\stepcounter{csicou}
			\end{enumerate}
			So fix $j$ with $\mathbf{x}_j \in \{\mathbf{un},\mathbf{eq}\}$.
			If $\mathbf{x}_j = \mathbf{un}$, then for each $k \leq l$ for  the natural number $i^j_k$ from Subclaim \ref{fasclelo} (and for some $c_0,c_1$) 
			\begin{equation} \label{negy} \begin{array}{rl} (q_k^{(\bar u^j, \bar v^j, \bar w^j, \bar s^j)})^{\{\varp^\iota_m\},(\langle 0 \rangle)}   \Vdash & \name z_{i^j_k} = c_0, \\ (q_k^{(\bar u^j, \bar v^j, \bar w^j, \bar s^j)})^{\{\varp^\iota_m\},(\langle 1 \rangle)}   \Vdash & \name z_{i^j_k} = c_1, \\
					\text{where }  c_0 \neq c_1, & \end{array} \end{equation}
				where we also have
				$$i^j_0 < i^j_1 < \dots < i^j_l.$$
		
		 (For convenience, if $\mathbf{x}_j = \mathbf{eq}$, then we let $i^j_k = -1$ for each $k \leq l$.)
			Observe that the fact that $(q_k^{(\bar u^j, \bar v^j, \bar w^j, \bar s^j)})^{\{\varp^\iota_m\},(\langle 0 \rangle)} \in D^{\iota,\mathbf{x}_j}_m(\name z)$, where  $\mathbf{x}_j \in \{ \mathbf{un}, \mathbf{eq} \}$ implies that 
			\begin{equation} \label{megy} \begin{array}{l}
					\text{whenever } p \leq_{\{\varp^\iota_m\},1} q_k^{(\bar u^j, \bar v^j, \bar w^j, \bar s^j)}), \text{ and } a \in \{0,1\}: \\
					(i \neq i^j_k) \rightarrow \ [(p^{\{\varp^\iota_m\},(\langle 0 \rangle)} \Vdash \name z_i = a \iff \ p^{\{\varp^\iota_m\},(\langle 1 \rangle)} \Vdash \name z_i = a]. \end{array} \end{equation}

			\begin{enumerate}[label = $(\circledcirc)_{\arabic*}$ , ref = $(\circledcirc)_{\arabic*}$]
				\setcounter{enumi}{\value{csicou}}
				\item \label{cicrc7} Now for any $\bar t \in \ ^{n\geq}2$, if $\mathbf{x}_j = \mathbf{un}$, then set $i^j_* = i^j_*(|\bar t|) = i^j_{|\bar t|}$, otherwise  if $\mathbf{x}_j = \mathbf{eq}$, then set $i^j_* = -1$. It suffices to show that 
				\begin{enumerate}[label = $\circ_{\arabic*}$, ref = $\circ_{\arabic*}$]
					\item \label{fon} whenever $i \in \omega$, $i \neq i^j_*$ and $r \leq_{\{\varp^\iota_m\},1} q_0^{(\bar u^j, \bar v^j, \bar w^j, \bar s^j)})^{\{\varp^\iota_m\}, (\bar t)}$ are such that $r^{\{\varp^\iota_m\}, (\langle 0 \rangle)}$, or $r^{\{\varp^\iota_m \}, (\langle 1 \rangle)}$ decides the value of $\name z_i$, then so does $r$, and
					\item \label{fo2} if $\mathbf{x}_j = \mathbf{un}$, and so $i^j_* \geq 0$, then  
					$$(q_0^{(\bar u^j, \bar v^j, \bar w^j, \bar s^j)})^{\{\varp^\iota_m\}, (\bar t \tieconcat \langle 0 \rangle)}, \text{ and } (q_0^{(\bar u^j, \bar v^j, \bar w^j, \bar s^j)})^{\{\varp^\iota_m\}, (\bar t \tieconcat \langle 1 \rangle)}$$ force different values to $\name z_{i^j_*}$.
				\end{enumerate} 		
			\end{enumerate} 	
			We fix $k\leq l$, and argue  \ref{fon} and \ref{fo2} simultaneously for each $\bar t \in \ ^k2$. Let $r \leq_{\{\varp^m\}, k+1} q_0^{(\bar u^j, \bar v^j, \bar w^j, \bar s^j)}$, (so
			\begin{equation} \label{egy} 
				\begin{array}{rl}
					\forall  \bar t \in \ ^k2: & r^{\{\varp^\iota_m\},(\bar t \tieconcat \langle 0 \rangle)} \leq (q_0^{(\bar u^j, \bar v^j, \bar w^j, \bar s^j)})^{\{\varp^\iota_m\},(\bar t \tieconcat \langle 0 \rangle)} \\ &
					r^{\{\varp^\iota_m\},(\bar t \tieconcat \langle 1 \rangle)} \leq (q_0^{(\bar u^j, \bar v^j, \bar w^j, \bar s^j)})^{\{\varp^\iota_m\},(\bar t \tieconcat \langle 1 \rangle)}) \end{array}
			\end{equation}
			by symmetry it is enough to show that
			\begin{equation} \label{ezkellmegy} \begin{array}{l}
					\forall i \neq i^j_k, \ \forall \bar t \in \ ^k2: \ \  (r^{\{\varp^\iota_m\}, (\bar t \tieconcat \langle 0 \rangle)} \Vdash \name z_i = a)
					\ \ \Leftrightarrow  \ \
					(r^{\{\varp^\iota_m\}, (\bar  t \tieconcat  \langle 1 \rangle)} \Vdash \name z_i = a), \\
				\end{array}
			\end{equation}
			and
			\begin{equation} \label{ezkellnegy} \begin{array}{l}
					\forall \bar t \in \ ^k2: \ \ (r^{\{\varp^\iota_m\},(\bar t \tieconcat \langle 0 \rangle)} \Vdash \name z_{i^j_*} = a)
					\ \ \Leftrightarrow  \ \
					(r^{\{\varp^\iota_m\},(\bar t \tieconcat  \langle 1 \rangle)} \Vdash \name z_{i^j_*} = 1-a). \\
				\end{array}
			\end{equation}
			
			
			\begin{enumerate}[label = $(\circledcirc)^{(2)}_{\arabic*}$ , ref = $(\circledcirc)´^{(2)}_{\arabic*}$]
				\setcounter{enumi}{\value{csicou}} 
				\item \label{as} We claim that for any $\bar t^* \in \ ^{\omega>}2$, $d \leq k$ and  $a \in \{0,1\}$, if $i \neq i^j_d$, then  
				\begin{equation} \label{kellnegy} r^{\{\varp^\iota_m\},(\vec{0}^{d} \tieconcat \langle 1 \rangle \tieconcat \bar t^*)} \Vdash \name z_i = a \ \iff \ r^{\{\varp^\iota_m\},(\vec{0}^{d} \tieconcat \langle 0 \rangle \tieconcat \bar t^*)} \Vdash \name z_i = a, \end{equation}
				and 
				\begin{equation} \label{kellmegy} r^{\{\varp^\iota_m\},(\vec{0}^{d} \tieconcat \langle 1 \rangle \tieconcat \bar t^*)} \Vdash \name z_{i^j_d} = a \ \iff \ r^{\{\varp^\iota_m\},(\vec{0}^{d} \tieconcat \langle 0 \rangle \tieconcat \bar t^*)} \Vdash \name z_{i^j_d} = 1-a, \end{equation}
				\stepcounter{csicou} 
			\end{enumerate}
			Before arguing \ref{as} first we note that it would finish the proof of \ref{fon} and \ref{fo2}: For any $\bar t \in \ ^k2$ and $i \in \omega$ (applying \ref{fon} $2\cdot |\{b <k: \ t_b = 1\}|$-many times) we obtain that 
			$$\forall a \in \{0,1\}: \ r^{\{\varp^\iota_m\},(\vec{0}^{k} \tieconcat \langle 0 \rangle)} \Vdash \name z_{i^j_l} = a \ \iff \ r^{\{\varp^\iota_m\}, (\bar t \tieconcat \langle 0 \rangle )} \Vdash \name z_{i^j_l} = f^{t_l}(a),$$ and 
			$$\forall a \in \{0,1\}: \ r^{\{\varp^\iota_m\},(\vec{0}^{k} \tieconcat \langle 1 \rangle)} \Vdash \name z_{i^j_l} = a \ \iff \ r^{\{\varp^\iota_m\}, (\bar t \tieconcat \langle 1 \rangle)} \Vdash \name z_{i^j_l} = f^{t_l}(a),$$
			where $f(a) = 1-a$, $f^0 = f^2 = f \circ f = $ id, and we mean $0$ under $t_l$, when $l\geq k$.
			But then by \ref{wlog}, \eqref{egy} we have $r^{\{\varp^\iota_m\},(\vec{0}^{k})} \leq_{\{ \varp^\iota_m\},1} (q_k^{(\bar u^j, \bar v^j, \bar w^j, \bar s^j)})$, which means
			\[ r^{\{\varp^\iota_m\},(\vec{0}^{k} \tieconcat \langle 0 \rangle)} \leq (q_k^{(\bar u^j, \bar v^j, \bar w^j, \bar s^j)})^{\{\varp^\iota_m\},(\langle 0 \rangle)}),\]
			\[ r^{\{\varp^\iota_m\},(\vec{0}^{k} \tieconcat \langle 1 \rangle)} \leq (q_k^{(\bar u^j, \bar v^j, \bar w^j, \bar s^j)})^{\{\varp^\iota_m\},(\langle 1 \rangle)}).\]
			Now if $i = i^j_*$, then this together with \eqref{negy} and \eqref{kellnegy} from \ref{as} imply	\eqref{ezkellnegy}. Similarly, for $i \neq i^j_*$ \eqref{megy} and \eqref{kellmegy} from \ref{as} imply	\eqref{ezkellmegy}. Hence it remains to verify \ref{as}.
			
			But \eqref{egy} and \ref{wlog} imply that
			$$ r^{\{\varp^\iota_m\},(\vec{0}^{d} \tieconcat \langle 0 \rangle \tieconcat \bar t^*)} \leq (q_d^{(\bar u^j, \bar v^j, \bar w^j, \bar s^j)})^{\{\varp^\iota_m\},(\langle 0 \rangle)}, $$
			$$ r^{\{\varp^\iota_m\},(\vec{0}^{d} \tieconcat \langle 1 \rangle \tieconcat \bar t^*)} \leq (q_d^{(\bar u^j, \bar v^j, \bar w^j, \bar s^j)})^{\{\varp^\iota_m\},(\langle 1 \rangle}, $$
			and clearly if $i = i^j_d$, then \eqref{negy} implies \eqref{kellnegy}, while for $i \neq  i^j_d$, then \eqref{kellmegy} follows from \eqref{megy}.

			\begin{enumerate}[label = $(\circledcirc)_{\arabic*}$ , ref = $(\circledcirc)_{\arabic*}$]
				\setcounter{enumi}{\value{csicou}}
				\item Now assuming that 	we have \ref{arl}, we prove that (by induction on $j$):
				\[ \forall j<M, \ \forall\bar t\in \ ^{n\geq}2: \ \mathbf{x}_m = \mathbf{mul} \ \rightarrow \  (q_0^{(\bar u^j, \bar v^j, \bar w^j, \bar s^j)})^{\{ \varp^\iota_m\}, (\bar t)} \in D^{\iota,\mathbf{x}_m}_m(\name z).  \]
				\stepcounter{csicou}
			\end{enumerate}	
			Assume that $j_*<M$ is such that $\mathbf{x}_{j_*} = \mathbf{mul}$, 
			\begin{equation} \label{<j*} \forall j < j_*, \ \forall \bar t \in \ ^{l\geq}2:  \ (q_0^{(\bar u^j, \bar v^j, \bar w^j, \bar s^j)})^{\{\varp^\iota_m\}, (\bar t)} \in D_m^{\iota,\mathbf{x}_j}(\name z).\end{equation}
			Fix $\bar t \in \ ^{l\geq}2$, and suppose on the contrary, that 
			$$(q_0^{(\bar u^{j_*}, \bar v^{j_*}, \bar w^{j_*}, \bar s^{j_*})})^{\{\varp^\iota_m\}, (\bar t)} \notin D_m^{\iota,\mathbf{mul}}(\name z).$$
			Recalling Definition \ref{defmul} for some $p \leq_{\{\varp^\iota_m\}, 1}  (q_0^{(\bar u^{j_*}, \bar v^{j_*}, \bar w^{j_*}, \bar s^{j_*})})^{\{\varp^\iota_m\}, (\bar t)}$ we have 
			$p \in D_m^{\iota,\mathbf{un}}(\name z) \cup D_m^{\iota,\mathbf{eq}}(\name z)$.
			Then there is a condition $p' \leq_{\bar{\bar{\ell}}^n} q_0$, for which 
			\begin{equation} \label{rr*} (p')^{(\bar u^{j_*}, \bar v^{j_*}, \bar w^{j_*}, \bar s^{j_*})}  \ \in\ (D_m^{\iota,\mathbf{un}}(\name z) \cup D_m^{\iota,\mathbf{eq}}(\name z)), \end{equation}
			and
			$$\forall j <M: \ (p')^{(\bar u^{j}, \bar v^{j}, \bar w^{j}, \bar s^{j})} \leq_{\{\varp^\iota_m\},1} (q_0^{(\bar u^{j}, \bar v^{j}, \bar w^{j}, \bar s^{j})})^{\{\varp^\iota_m\}, (\bar t)}, $$
			so
			\begin{equation} \label{emo} (p')^{(\bar u^{j_*}, \bar v^{j_*}, \bar w^{j_*}, \bar s^{j_*})} \in  D_m^{\iota,\mathbf{un}}(\name z) \cup D_m^{\iota,\mathbf{eq}}(\name z), \end{equation}
			and recalling \ref{elfaop} from Fact \ref{elfa} we could infer from \ref{<j*} that 
			\begin{equation} \label{r<j*} \forall j < j_*:  \ ((p')^{(\bar u^j, \bar v^j, \bar w^j, \bar s^j)}) \in D_m^{\iota,\mathbf{x}_j}(\name z).\end{equation}
			But now, since $p' \leq_{\bar{\bar{\ell}}^n} q_0 \leq_{\bar{\bar{\ell}}^n} q_{**}$, $\mathbf{x}_{j_*} = \mathbf{mul}$, necessarily 
			$$(p')^{(\bar u^{j_*}, \bar v^{j_*}, \bar w^{j_*}, \bar s^{j_*})}\in D_m^{\iota,\mathbf{mul}}(\name z),$$ 
			contradicting \eqref{rr*} (as $D^{\iota,\textbf{un}}_{m}(\name z)$, $D^{\iota,\textbf{mul}}_{m}(\name z)$, $D^{\iota,\textbf{eq}}_{m}(\name z)$ are pairwise disjoint by obvious reasons (\ref{triv}).
		\end{PROOF}
		
		\begin{sclaim} \label{mul+}
			Suppose that $\iota \in \{\mathbf{1}, \infty\}$, let the condition $q$ be  in $ D^{\iota,\mathbf{mul}}_m(\name z)$.
			Then some $q' \leq_{\{\varp^\iota_m\},1} q $ satisfies the following:
			
			There exist $i_* \neq i_{**} \in \omega$, $c_*,c_{**} \in \{0,1\}$, such that
			\[ \begin{array}{rrclrrcl} (q')^{\{\varp^\iota_m\}, (\langle 0 \rangle)} \Vdash & \name z_{i_*}& = &c_*, &  (q')^{\{\varp^\iota_m\}, (\langle 1 \rangle)} \Vdash & \name z_{i_*} &= &1- c_* \\
				(q')^{\{\varp^\iota_m\}, (\langle 0 \rangle)} \Vdash & \name z_{i_{**}} & = & c_{**}, &  (q')^{\{\varp^\iota_m\}, (\langle 1 \rangle)} \Vdash & \name z_{i_{**}} & = & 1- c_{**}.	
			\end{array}
			\]
			Moreover, $q'$ can be chosen so that both $(q')^{\{\varp^\iota_m\}, (\langle 0 \rangle)}$ and $(q')^{\{\varp^\iota_m\}, (\langle 1 \rangle)}$ decide the first $\max(i_*,i_{**})+1$-many digits of $\name z$.
		\end{sclaim}
		\begin{PROOF}{Subclaim \ref{mul+}}
			For each $i \in \omega$ there exists $q_+ \leq_{\{\varp^\iota_m\},1} q$, such that both $q_+^{\{\varp^\iota_m\},(\langle 0 \rangle)}$ and $q_+^{\{\varp^\iota_m\},(\langle 1 \rangle)}$ decide $\name z_i$. Since $q \notin D^{\iota,\mathbf{eq}}_m(\name z)$, for some $i_* \in \omega$ and $q_+ \leq_{\{\varp^\iota_m\},1} q$   the conditions $q_+^{\{\varp^\iota_m\},(\langle 0 \rangle)}$ and $q_+^{\{\varp^\iota_m\},(\langle 1 \rangle)}$ decide about $\name z_{i_*}$ differently. Since $q_+ \notin D^{\iota,\mathbf{eq}}_m(\name z)$ (again by $q_+ \leq_{\{\varp^\iota_m\},1} q$ and Definition \ref{defmul}) there exists $i_{**} \neq i_*$, such that $q' \leq_{\{\varp^\iota_m\},1} q_+$, and $(q')^{\{\varp^\iota_m\},(\langle 0 \rangle)}$ and $(q')^{\{\varp^\iota_m\},(\langle 1 \rangle)}$ force different values to $\name z_{i_{**}}$.
		\end{PROOF}
		\begin{sclaim} \label{ekl}
			Assume that  $\iota \in \{ \mathbf{1}, \infty \}$, $n$, $m$, $q_*$ are as in Lemma \ref{guszkl0}, $q_{**}$, $\bar{\mathbf{x}}^\iota$ are as in \ref{hulye}, moreover, there is no $j<M$ for which $\mathbf{x}_j = $$\mathbf{un}$.
			Then if $r_{**} \leq_{\bar{\bar{\ell}}^n} q_{**}$ is such that 
			$$r_{**} \in \bigcap_{j<M} \{q \leq_{\bar{\bar{\ell}}^n} q_{**}: \ q^{(\bar u^j, \bar v^j, \bar w^j, \bar s^j)} \in D_m^{\iota,\mathbf{x}_j}(\name z)\},$$ 
			then there is an $r_{*} \leq{\bar{\bar{\ell}}^{n+1}} r_{**}$ that satisfies the requirements in \ref{bbf1} of Lemma \ref{guszkl0}.
		\end{sclaim}
		\begin{PROOF}{Subclaim \ref{ekl}}
			 Note that all the requirements except \ref{mul} hold for $r_{**}$, and so for any $r_* \leq{\bar{\bar{\ell}}^{n+1}} r_{**}$ too by Observation \ref{obskite}. Now we only have to appeal to Subclaim \ref{mul+} $|\{j<M: \ \mathbf{x}_j = \mathbf{mul}\}|$-many times.
		\end{PROOF}
		
		\begin{sclaim} \label{ekl1.5}
			Assume that  $\iota = \infty $, $n$, $m$, $q_*$ are as in Lemma \ref{guszkl0}, $q_{**} \leq_{\bar{\bar{\ell}}^n} q_*$, $\bar{\mathbf{x}}$ are as in \ref{hulye}, and suppose that there exists $j<M$ with $\mathbf{x}_j = \mathbf{un}$.
					Then there  exists $r_* \leq_{\bar{\bar{\ell}}^n} q_{**}$,  and $\bar{\mathbf{x}}'$ that satisfy \ref{bbf1} of Lemma \ref{guszkl0} (with some $\bar x$), where $$ \begin{array}{rlcl}\forall j<M: & \mathbf{x}_j \in \{ \mathbf{mul}, \mathbf{eq}\} & \to & \mathbf{x}_j = \mathbf{x}_j', \\
					&	 \mathbf{x}_j = \mathbf{un} & \to & \mathbf{x}'_j = \mathbf{mul}. \end{array}
					$$
		\end{sclaim}
		\begin{PROOF}{Subclaim \ref{ekl1.5}}
			Proceed first similarly to the proof of Subclaim \ref{ekl}, and appeal to Subclaim \ref{mul+} $|\{j<M: \ \mathbf{x}_j = \mathbf{mul}\}|$-many times, and so for some $p_{*} \leq_{\bar{\bar{\ell}}^n} q_{**}$ we have that
			\begin{enumerate}[label = $\boxdot_{(\alph*)}$, ref  = $\boxdot_{(\alph*)}$] 
				\item \label{D} $ \forall j < M: \  p_{*}^{(\bar u^j, \bar v^j, \bar w^j, \bar s^j)} \in D_m^{\iota,\mathbf{x}_j}(\name z)$,
				\item \label{ij*} if $\mathbf{x}_j = \mathbf{mul}$, then for some $i_j^* \neq i_j^{**}$ the conditions
				$$(p_{*}^{(\bar u^j, \bar v^j, \bar w^j, \bar s^j)})^{(\{\varp^\infty_m\}, (\langle 0 \rangle))}, (p_{*}^{(\bar u^j, \bar v^j, \bar w^j, \bar s^j)})^{(\{\varp^\infty_m\}, (\langle 1 \rangle))}$$
				decide differently about $\name z_{i_j^*}$, as well as about $\name z_{i_j^{**}}$. Moreover, both conditions decide $\name z \rest [0, \max(i_j^*, i_j^{**})+1]$.
			\end{enumerate}
			Now pick $p_{**} \leq_{\bar{\bar{\ell}}^n} p_*$, so that
			
			\begin{enumerate}[label = $\bigodot_{(\alph*)}$, ref  = $\bigodot_{(\alph*)}$] 
				\item 	$\forall j<M: \ p_{**}^{(\bar u^j, \bar v^j, \bar w^j, \bar s^j)} \leq (p_{*}^{(\bar u^j, \bar v^j, \bar w^j, \bar s^j)})^{(\{\varp^\infty_m\}, (\langle 0 \rangle)},$
				\item \label{D**}	and 
				$ \forall j < M: \  p_{**}^{(\bar u^j, \bar v^j, \bar w^j, \bar s^j)} \in D_m^{\iota,\mathbf{x}_j}(\name z).$
			\end{enumerate}
			So (by replacing $p_*$ with a $\leq_{\bar{\bar{\ell}}^{n+1}}$-extension of it) w.l.o.g.\ we can assume that 
			\begin{equation} \label{elagh} \forall j<M: \ p_{**}^{(\bar u^j, \bar v^j, \bar w^j, \bar s^j)} = (p_{*}^{(\bar u^j, \bar v^j, \bar w^j, \bar s^j)})^{(\{\varp^\infty_m\}, (\langle 0 \rangle)}. \end{equation}
			Define $r_* \leq_{\bar{\bar{\ell}}^n} p_*$ so that
			\begin{equation}\label{r*0} \forall j<M: \ (r_{*}^{(\bar u^j, \bar v^j, \bar w^j, \bar s^j)})^{\{\varp^\infty_m\}, (\langle 0 \rangle)} = (p_*^{(\bar u^j, \bar v^j, \bar w^j, \bar s^j)})^{\{\varp^\infty_m\}, (\langle 0, 0 \rangle)}, \end{equation}
			and similarly,
			\begin{equation} \label{r*1} \forall j<M: \ (r_{*}^{(\bar u^j, \bar v^j, \bar w^j, \bar s^j)})^{\{\varp^\infty_m\}, (\langle 1 \rangle)} = (p_*^{(\bar u^j, \bar v^j, \bar w^j, \bar s^j)})^{\{\varp^\infty_m\}, (\langle 1, 1 \rangle)}. \end{equation}
			
			Now if $j<M$ is such that $\mathbf{x}_j = \mathbf{mul}$, then for $a \in \{0,1\}$
				$$(r_{*}^{(\bar u^j, \bar v^j, \bar w^j, \bar s^j)})^{\{\varp^\infty_m\}, (\langle a \rangle)} \leq (p_*^{(\bar u^j, \bar v^j, \bar w^j, \bar s^j)})^{\{\varp^\infty_m\}, (\langle a \rangle)}, $$
			so $(r_{*}^{(\bar u^j, \bar v^j, \bar w^j, \bar s^j)})^{\{\varp^\infty_m\}, (\langle a \rangle)}$ ($a \in \{0,1\}$) decide differently about $\name z_{i^*_j}$ and $\name z_{i^{**}_j}$ (by \eqref{ij*}).
			
			If $j<M$ is such that $\mathbf{x}_j = \mathbf{eq}$, then it follows from \eqref{elagh}, \ref{D}, \ref{D**} and \ref{more} from Subclaim \ref{fascl}, that 
			$$\forall  t_0 \in \ \{0,1\}: \ (p_{*}^{(\bar u^j, \bar v^j, \bar w^j, \bar s^j)})^{\{\varp^\infty_m\}, (\langle t_0\rangle)} \in D^{\iota,\mathbf{eq}}_m(\name z),$$
			but then a similar straightforward  calculation shows that for no $i$ and $p' \leq_{\{\varp^\infty_m\}, 2} p_*^{(\bar u^j, \bar v^j, \bar w^j, \bar s^j)}$, no $\bar t_0 \neq \bar t_1 \in \ ^2\{0,1\}$ do $(p')^{\{\varp^\infty_m\}, (\bar t_0)}$ and $(p')^{\{\varp^\infty_m\}, (\bar t_1)}$ decide differently about $\name z_i$. This clearly implies that $r_{*}^{(\bar u^j, \bar v^j, \bar w^j, \bar s^j)}) \in D^{\iota,\mathbf{eq}}_m(\name z)$.
			
			Finally, if $j$ is such that $\mathbf{x}_j = \mathbf{un}$, then we argue that
			$r_{*}^{(\bar u^j, \bar v^j, \bar w^j, \bar s^j)} \in D^{\iota,\mathbf{mul}}_m(\name z)$.
			 Again, \eqref{elagh}, \ref{D}, \ref{D**} and \ref{more} from Subclaim \ref{fascl} together imply that 
			 $$\forall  t_0 \in \ \{0,1\}: \ (p_{*}^{(\bar u^j, \bar v^j, \bar w^j, \bar s^j)})^{\{\varp^\infty_m\}, (\langle t_0\rangle)} \in D^{\iota,\mathbf{un}}_m(\name z),$$
			 and similarly to the argument in \ref{arl} in Subclaim \ref{fascl} there are $i_0 < i_1$ and $c_0,c_1 \in \{0,1\}$, such that
			 $$ \forall  \bar t = \langle t_0, t_1 \rangle \in \ ^2 \{0,1\}: \ (p_{*}^{(\bar u^j, \bar v^j, \bar w^j, \bar s^j)})^{\{\varp^\infty_m\}, (\langle t_0, t_1\rangle)} \Vdash \name z_{i_0} = f^{t_0}(c_0) \ \wedge z_{i_0} = f^{t_1}(c_1),$$
			 (where $f(c) = 1-c$, $f^0 =$ id).
			 From this we obtain
			  $$  \begin{array}{rcl} (p_{*}^{(\bar u^j, \bar v^j, \bar w^j, \bar s^j)})^{\{\varp^\infty_m\}, (\langle 0,0 \rangle)} & \Vdash & \name z_{i_0} = c_0 \ \wedge \name z_{i_1} = c_1, \\
			   (p_{*}^{(\bar u^j, \bar v^j, \bar w^j, \bar s^j)})^{\{\varp^\infty_m\}, (\langle 1,1 \rangle)} & \Vdash & \name z_{i_0} = 1- c_0 \ \wedge \name z_{i_1} = 1- c_1, \end{array}$$
			  so recalling \eqref{r*0}, \eqref{r*1} clearly $r_*^{(\bar u^j, \bar v^j, \bar w^j, \bar s^j)} \in D^{\infty,\mathbf{mul}}_m(\name z)$. (And since we can always $\leq_{\bar{\bar{\ell}}^{n+1}}$-extend $r_*$ so that $(r_*^{(\bar u^j, \bar v^j, \bar w^j, \bar s^j)})^{\{\varp^\infty\},(\langle a \rangle)}$ ($a \in \{0,1\}$) decides $\name z \rest [0, \max(i_0,i_1)+1)$, which gives the $x^{(\bar u, \bar v, \bar w, \bar s)}$'s for the two  $(\bar u, \bar v, \bar w, \bar s) \in \mathbf{seq}_{\bar{\bar{\ell}}^{n+1}}(r_*)$ for which
			  $$ (\bar u^j, \bar v^j, \bar w^j, \bar s^j) \sqsubseteq (\bar u, \bar v, \bar w, \bar s).)$$
		\end{PROOF}

		\begin{sclaim} \label{ekl2}
			Assume that $\iota = \mathbf{1}$, $q_*$, $\name z$, $n$, $m$ are as in Lemma \ref{guszkl0},   $q_{**} \leq_{\bar{\bar{\ell}}^n} q_*$, $\bar{\mathbf{x}}$ are as in \ref{hulye}, $l \in \omega$, moreover, 
			\begin{equation} \label{rossz} 0<|\{j<M: \ \mathbf{x}_j = \mathbf{un}\}|  < 2^l. \end{equation} 
			Suppose that $q_{***} \leq_{\bar{\bar{\ell}}^n}  q_{**}$ is given by applying Subclaims \ref{fasclelo} and \ref{fascl} to $q_{**}$ and $l$.
			Then there exists $\bar t \in \ ^{l}2$ for which some $r_{*} \leq_{\bar{\bar{\ell}}^{n+1}} r_{\bar t}$  satisfies the requirements in \ref{bbf1}, (where $r_{\bar t}$ is 
			defined by the equality
			\[ \forall j <M: \ r_{\bar t}^{(\bar u^j, \bar v^j, \bar w^j, \bar s^j)} = (q_{***}^{(\bar u^j, \bar v^j, \bar w^j, \bar s^j)})^{\{\varp^\iota_m\}, (\bar t)}). \]
			
		\end{sclaim}
		\begin{PROOF}{Subclaim \ref{ekl2}}
			
		Fixing $j<M$ so that $\mathbf{x}_j = \mathbf{un}$,
		$$\forall \bar t \in \ ^l 2: \ (q_{***}^{(\bar u^j, \bar v^j, \bar w^j, \bar s^j)})^{\{\varp^\iota_m\}, (\bar t)}) \in D^{\mathbf{1},\mathbf{un}}_m(\name z), $$
		
		and by \eqref{unifo}
		 \begin{enumerate}
		 	\item 	 there are  natural numbers 
		 	$$i^j_0 < i^j_1 < \dots < i^j_l,$$
		 	such that
		 	$$	\forall \bar t \in \ ^{l} 2, \ \forall a \in \{0,1\}: \ (q_{***}^{(\bar u^j, \bar v^j, \bar w^j, \bar s^j)})^{\{\varp^\iota_m\}, (\bar t \tieconcat \langle a \rangle)}  \parallel \name z \rest [0,i^j_l+1).$$
		 	\begin{itemize}
		 		\item $	\forall \bar t \in \ ^{l} 2, \ \exists \bar z^j_{\bar t} \in \ ^{i^j_l}2: \  (q_{***}^{(\bar u^j, \bar v^j, \bar w^j, \bar s^j)})^{\{\varp^\iota_m\}, (\bar t)})  \Vdash \bar z^j_{\bar t} \subseteq \name z$,
		 		\item if $\bar t \in \ ^{l\geq } 2, \ i \neq i^j_{|\bar t|}$, $r \leq_{\{ \varp^{\mathbf{1}}_m \},1} (q_{***}^{(\bar u^j, \bar v^j, \bar w^j, \bar s^j)})^{\{\varp^\iota_m\}, (\bar t)})$, then 
		 		$$ r^{\{\varp^\iota_m\}, (\bar t \tieconcat \langle 0 \rangle)} \Vdash \name z_i = a \ \iff \ r^{\{\varp^\iota_m\}, (\bar t \tieconcat \langle 1 \rangle)} \Vdash \name z_i = a \ (\forall a \in \{0,1\}),$$
		 		\item if $\bar t \in \ ^{l\geq } 2$, $r \leq_{\{ \varp^{\mathbf{1}}_m \},1} (q_{***}^{(\bar u^j, \bar v^j, \bar w^j, \bar s^j)})^{\{\varp^\iota_m\}, (\bar t)})$, $a \in \{0,1\}$, then 
		 		$$ r^{\{\varp^\iota_m\}, (\bar t \tieconcat \langle 0 \rangle)} \Vdash \name z_{i^j_{|\bar t|}} = a \ \iff \ r^{\{\varp^\iota_m\}, (\bar t \tieconcat \langle 1 \rangle)} \Vdash \name z_{i^j_{|\bar t|}} = 1- a.$$
		 	\end{itemize}
		 \end{enumerate}
	 	\begin{enumerate}
			\item[(2)] \label{kul} Observe that if $j<M$, $\bar t \neq \bar t'$, then $\bar z^j_{\bar t} \neq \bar z^j_{\bar t'}$.
		\end{enumerate}

	
	Now for \ref{bbf1} 
		\begin{enumerate}
			\item[(3)] 	it suffices to choose $\bar t \in \ ^l 2$, and set $r_* = r_{\bar t}$ so that whenever $j<M$ is such that $\mathbf{x}_j = \mathbf{un}$, then 
			\begin{equation} \label{nems} \bar z^j_{\bar t} \neq \mathfrak{s}_{i^j_l}, \end{equation}
			which is shown by the following:
			 since for any $j$ with $\mathbf{x}_j = \mathbf{eq}$ obviously $r_{\bar t}^{(\bar u^j, \bar v^j, \bar w^j, \bar s^j)} \in D^{\mathbf{1},\mathbf{eq}}_m$, while if $\mathbf{x}_j = \mathbf{mul}$, then we can  replace $r_*$ with some $r_{**} \leq_{\{\varp^\mathbf{1}_m\},1} r_* = r_{\bar t}$ given by Subclaim \ref{mul+}, preserving 
			$$r_{**} \in \bigcap_{j<M}\{p \in \bbQ': \ p^{(\bar u^j, \bar v^j, \bar w^j, \bar s^j)} \in D_m^{\iota,\mathbf{x}_j}(\name z) \}.$$
			
		\end{enumerate}
	But for each $j<M$ with $\mathbf{x}_j = \mathbf{un}$ there is at most one $\bar t \in \ ^{i^j_l} 2$ that does not satisfy \eqref{nems}, so by \eqref{rossz} there exists a sequence $\bar t$ that is suitable for our demands.

		\end{PROOF}
	\end{PROOF}

	\end{PROOF}
	It is only left to argue \ref{e3}, that will complete the proof of Theorem $\ref{main}$.
	
	So fix 
	\begin{itemize}
		\item 	$\alpha < \lambda_\mathbf{0}$, and a $\bbG_{\mathbf{0}}(\bar{\mathfrak{s}})$-independent tree
		$$T_\mathbf{0} \in V_\mathbf{0} \cap \sP(^{\omega>}2),$$
		where $V_\mathbf{0} = V^{\bbQ^{\mathbf{0}}_{\lambda_\mathbf{0} \setminus \{ \alpha \}} \times \bbQ^{\mathbf{1}} \times \bbQ^{\mathbf{\infty}} \times \bbQ^{\bbS}}$,
		\item $\beta \in [\lambda_\mathbf{0}, \lambda_\mathbf{1})$, and a $\bbG_{\mathbf{1}}$-independent tree
		$$T_\mathbf{1} \in V_\mathbf{1} \cap \sP(^{\omega>}2),$$
		where $V_\mathbf{1} = V^{\bbQ^{\mathbf{0}} \times \bbQ^{\mathbf{1}}_{\lambda_\mathbf{1} \setminus \lambda_{\mathbf{0}} \setminus \{ \beta \}} \times \bbQ^{\mathbf{\infty}} \times \bbQ^{\bbS}}$,
		\item $\gamma \in [\lambda_\mathbf{1}, \lambda_\mathbf{\infty})$, and an $E_0$-independent tree
		$$T_\mathbf{\infty} \in V_\infty  \cap \sP(^{\omega>}2),$$
		where 		$V_\infty = V^{\bbQ^{\mathbf{0}} \times \bbQ^{\mathbf{1}} \times \bbQ^{\mathbf{\infty}}_{\lambda_\mathbf{\infty} \setminus \lambda_{\mathbf{1}} \setminus \{ \gamma \}}  \times \bbQ^{\bbS}}$,
	\end{itemize}
	and we shall check that the generic real in question is not in $[T_\iota]$ ($\iota \in \{ \mathbf{0}, \mathbf{1}, \infty\}$).
	Assume on the contrary (i.e. $\neg$ \ref{e3}), let $p_\iota \in \bbP^\iota$ be such that
	\begin{equation} \label{releme} V_\iota \models \ ``p \Vdash_{\bbP^\iota} \name r \in [T_\iota]", \end{equation}	
	where $\name r$ is the generic real given by $\bbP^\iota$
	(note that $\bbP^\iota \in V$, and so we have to carefully manipulate $p$ when working in $V_\iota$ as there are more reals in that model than in $V$).
	
	By Definition $\ref{def1}$ and $p \in \bbP^\iota$ 
	\begin{itemize}
		\item if $\iota = \mathbf{0}$, then there exists $j \in \omega$ with $p_{2j} = C_{2j}$, $p_{2j+1} = C_{2j+1}$. W.l.o.g.\ we can assume that $|p_0| = |p_1| = \dots = |p_{2j-1}| = 1$, and if $p_i = \{ t_i \}$ ($i<2j$), then let 
		$$\bar t^* = \bar t_0 \tieconcat \bar t_1 \tieconcat \dots \tieconcat \bar t_{2j-1},$$
		\item if $\iota = \mathbf{1}$, then 
		 there exists $j \in \omega$ with $p_{j} = \{0,1\}$, and $|p_0| = |p_1| = \dots = |p_{2j-1}| = 1$, and if $p_i = \{ a_i \}$ ($i<j$), then let 
		$$\bar t^* = \langle a_i: \ i<j \rangle,$$
			\item if $\iota = \mathbf{\infty}$, then 
		there exists $j \in \omega$ with 
		 \begin{equation} \label{tj'} p_{j} = \{\bar t'_j, \bar t_j''\}, \end{equation}
		  and $|p_0| = |p_1| = \dots = |p_{2j-1}| = 1$, and if $p_i = \{ \bar t_i  \}$ ($i<j$), then let 
		$$\bar t^* = \bar t_0 \tieconcat \bar t_1 \tieconcat \dots \tieconcat \bar t_{j-1},$$
	\end{itemize}
	Now
		\begin{itemize}
		\item if $\iota = \mathbf{0}$, then using \ref{D1} pick $t_{2j} \in p_{2j}$ so that $$\bar t^{**} := \bar t^* \tieconcat \bar  t_{2j} = \mathfrak{s}_k$$
		 for some $k \in \omega$. Letting $p' \in \bbP^\mathbf{0}$ denote a condition for which $p' \leq p$, $p'_{2j} = \{t_{2j} \}$,
			\item if $\iota = \mathbf{0}$, or $\infty$, then
			 $$\bar t^{**} := \bar t^*,$$
			 and  let $p'= p$.
		\end{itemize}
		Clearly
	$$p' \Vdash \name r \in [\bar t_{**}],$$
	so we can assume that $[t_{**}] \cap [T_\iota] \neq \emptyset$.
	Consider
	\begin{itemize}
		\item 	 the sets
		$$T_{\bar t_{**}(0)} = \{ \bar t \in \ ^{\omega\geq}2: \ \bar t_{**} \tieconcat \langle 0 \rangle \tieconcat \bar t \in T\},$$
		and 
		$$T_{\bar t_{**}(1)} = \{ \bar t \in \ ^{\omega\geq}2: \ \bar t_{**} \tieconcat \langle 1 \rangle \tieconcat \bar t \in T\},$$
		if $\iota = \mathbf{0}$ or $\mathbf{1}$,
		\item while if $\iota = \infty$, then let
			$$T_{\bar t_{**}(0)} = \{ \bar t \in \ ^{\omega\geq}2: \ \bar t_{**} \tieconcat \bar t_j' \tieconcat \bar t \in T\},$$
		and 
		$$T_{\bar t_{**}(1)} = \{ t \in \ ^{\omega\geq}2: \ \bar t_{**} \tieconcat \bar t_j'' \tieconcat \bar t \in T\},$$
		(where $p'_j = \{\bar t_j', \bar t_j'' \}$ recalling \eqref{tj'}).
	\end{itemize}
Now as $[T_\mathbf{0}]$ ($[T_\mathbf{1}]$, $[T_\mathbf{\infty}]$, resp.) is 
	$\bbG_{\mathbf{0}}(\bar{\mathfrak{s}})$-independent ($\bbG_\mathbf{1}$-, $E_0$-independent, resp.) compact set for which $t_{**} \in T$, there must be $k \in \omega$ such that
	the sets  $T_{\bar t_{**}(0)} \cap \ ^k2$, and $T_{\bar t_{**}(1)} \cap \ ^k2$ are disjoint.
	
	Now by further extending $p'$ if necessary we can assume that
		\begin{itemize}
		\item (if $\iota = \mathbf{0}$) 
	 $|p'_{2j+2}| = |p'_{2j+3}| = \dots = |p'_{2j+k+1}| = 1$, and if $p'_{2j+2+i} = \{ \bar  t_{2j+2+i} \}$ ($i<k$), then  the sequence $$\bar t_{***} = \bar t_{2j+2} \tieconcat \bar t_{2j+3} \tieconcat \dots \tieconcat \bar t_{2j+k+1} \in \ ^{\omega\geq}2$$
	is obviously of length $\geq k$.
		\item (if $\iota = \mathbf{1}$) 
		$|p'_{j+2}| = |p'_{j+3}| = \dots = |p'_{j+k+1}| = 1$, and if $p'_{j+2+i} = \{ a_{j+2+i} \}$ ($i<k$), then  the sequence $$\bar t_{***} = \langle a_{j+2+i}: \ i<k \rangle \in \ ^{\omega\geq}2$$
		is of length $k$,
			\item (if $\iota = \mathbf{\infty}$) 
		$|p'_{j+2}| = |p'_{j+3}| = \dots = |p'_{j+k+1}| = 1$, and if $p'_{2j+2+i} = \{ \bar \bar t_{2j+2+i} \}$ ($i<k$), then  the sequence $$\bar t_{***} = \bar t_{2j+2} \tieconcat \bar t_{2j+3} \tieconcat \dots \tieconcat \bar t_{2j+k+1} \in \ ^{\omega\geq}2$$
		is obviously of length $\geq k$.
	\end{itemize}
	If $\iota \in \{ \mathbf{0}, \mathbf{1} \}$, then
	let $a \in \{0,1 \}$ be such that $t_{***} \rest k \in T_{\bar t_{**}(a)}$.
	Our observation above means that $t_{***} \rest k \notin T_{\bar t_{**}(1-a)}$, thus
	\begin{equation} \label{notT} 
	 			\bar t_{**} \tieconcat \langle 1-a \rangle \tieconcat \bar t_{***} \notin T_\iota.
	\end{equation}
	Extend $p'$ to $p'' \in \bbQ$ such that $p''_{2j+1} = \{ \langle 1-a \rangle \}$ (if $\iota = \mathbf{0}$), or .$p''_{j+1} = \{ 1-a \}$ (if $\iota = \mathbf{1}$),
	and then
	$$ p'' \Vdash \name r \in [\bar t_{**} \tieconcat \langle 1-a \rangle \tieconcat \bar t_{***}]$$
	which together with $\eqref{notT}$	contradicts  $\eqref{releme}$.
	We can also reach the same contradiction in the case $\iota = \infty$, just working with $p_j = \{ \bar t_j', \bar t_j''\}$ instead of $\{0,1\}$.

\end{PROOF}

\begin{problem}
	Is it true, that in the final model there is a partition of the Cantor space into $\lambda_\mathbf{0}$-many $\mathbb{G}_{0}(\bar{\mathfrak{s}})$-independent Borel sets (while the other assertions from Theorem $\ref{main}$ still hold)? Is it consistent that there is a partition of $2^\omega$ into $\lambda$-many $\mathbb{G}_{0}(\bar{\mathfrak{s}})$-independent Borel sets, where $ \cov(\cM) < \lambda < 2^{\aleph_0}$, and less than $\lambda$-many (or just  not more than $\cov(\cM)$-many) $\mathbb{G}_{0}(\bar{\mathfrak{s}})$-independent Borel sets do not cover $2^\omega$? 
	What can we say about the corresponding invariant of $\mathbb{G}_{1}$, or $E_0$?
\end{problem} 

\begin{problem}
	Define the graph 
	$$ \bbG_n = \{ (x,y) \in \ [ ^\omega 2]^2: \  |\{j \in \omega: \ x_j \neq y_j \} | \leq n \}$$
	for $n \in \omega$ fixed.
	Can we separate $\cov(I_{\bbG_n})$ and $\cov(I_{\bbG_{n+1}})$?
	Can we separate infinitely many $\cov(T_{\bbG_j})$'s?
\end{problem} 

\section{Acknowledgement}
The first named author is grateful for  M. Gaspar for some valueable discussions. We thank Shimon Garti for the helpful suggestions and remarks, which improved the paper.

\end{PROOF}

\bibliographystyle{amsalpha}
\bibliography{shlhetal,F1988}

\end{document}